\documentclass[reqno,10pt]{amsart}
\usepackage{amssymb,amsmath,amscd, amsthm,stmaryrd,verbatim,
	mathrsfs,tikz,color , amsfonts, enumerate, esvect
}
\usepackage{enumitem}
\usepackage{braket}
\usetikzlibrary{patterns}
\usetikzlibrary{er}
\usepackage{genyoungtabtikz}

\usepackage[colorlinks,linkcolor=blue,urlcolor=cyan,citecolor=blue]{hyperref}

\numberwithin{equation}{section}
\allowdisplaybreaks

\usepackage[style=alphabetic,sorting=nyt]{biblatex}
\usepackage{fullpage}
\usepackage[textwidth=0.85in, textsize=tiny]{todonotes}
\newcommand{\sgn}{\operatorname{sgn}}

\theoremstyle{plain}
\newtheorem{Thm}{Theorem}[section]
\newtheorem{Prop}[Thm]{Proposition}
\newtheorem{Lem}[Thm]{Lemma}
\newtheorem{Cor}[Thm]{Corollary}

\newtheorem*{thm*}{Theorem}
\theoremstyle{definition}
\newtheorem{Defff}[Thm]{Definition}
\newtheorem{ex}[Thm]{Example}
\newtheorem{Rem}[Thm]{Remark}

\newcommand{\lam}{\lambda}

\newcommand{\al}{\alpha}

\newcommand{\ep}{\varepsilon}

\newcommand{\Ann}{\mathrm{Ann}}
\newcommand{\eq}{\begin{equation}}
\newcommand{\en}{\end{equation}}
\newcommand{\beqna}[1]{\begin{eqnarray}\label{#1}}
\newcommand{\eeqna}{\end{eqnarray}}
\newcommand{\beqn}[1]{\begin{equation}\label{#1}}
\newcommand{\eeqn}{\end{equation}}

\newcommand{\mf}[1]{\mathfrak{#1}}
\renewcommand{\subset}{\subseteq}

\newcommand{\bil}[2]{\langle{#1},{#2}^{\vee} \rangle }
\newcommand{\hs}{ \mathfrak{h}^*}

\newcommand{\gkd}{\operatorname{GKdim}}

\newcommand{\nodd}{}

\newcommand{\hobox}[3]{\draw (0+#1,0-#2) rectangle (1+#1,-1-#2)++(-0.5,+0.5) node {$ #3$};}

\newcommand{\domscale}{0.51}
\newcommand{\sh}{\operatorname{sh}}
\newcommand{\ev}{^{\mathrm{ev}}}
\newcommand{\od}{^{\mathrm{odd}}}

\usepackage[centertableaux]{ytableau}
\colorlet{srcol}{black!15}

\makeatletter
\pgfkeys{/ytableau/options,
  noframe/.default = false,
  noframe/.is choice,
  noframe/true/.code = {%
    \global\let\vrule@YT=\vrule@none@YT
    \global\let\hrule@YT=\hrule@none@YT
  },
  noframe/false/.code = {%
    \global\let\vrule@YT=\vrule@normal@YT
    \global\let\hrule@YT=\hrule@normal@YT
  },
  noframe/on/.style = {noframe/true},
  noframe/off/.style = {noframe/false},
}
\makeatother
\ytableausetup{noframe=off,boxsize=1.7em}
\let\ytb=\ytableaushort

\newcommand{\tytb}[1]{{\tiny\ytb{#1}}}

\makeatletter

\newcommand{\trivial}[2][]{\if\relax\detokenize{#1}\relax{\color{red}\vspace{0em} $[$  #2 $]$}\else
\ifx#1h
\ifcsname showtrivial\endcsname{\color{orange} \vspace{0em}  $[$ #2 $]$}\fi\else{\colorWrong argument!}\fi\fi\ignorespaces}

\def\Univ{{\mathfrak{U}}}

\def\Irr{\mathrm{Irr}}
\def\barP{\bar{P}}
\def\floor#1{\lfloor{#1}\rfloor}
\def\bfdd{{\mathbf{d}}}
\def\cO{{\mathcal{O}}}
\def\Spr{{\mathrm{Springer}}}

\def\abs#1{\left|#1\right|}
\def\Ind{\mathrm{Ind}}

\def\cO{\mathcal{O}}
\def\fhh{\mathfrak{h}}
\def\fgg{\mathfrak{g}}
\def\so{\mathfrak{so}}
\def\sp{\mathfrak{sp}}
\def\bfdd{\mathbf{d}}
\def\bfff{\mathbf{f}}
\def\bfdp{\mathbf{p}}
\def\bfda{\mathbf{a}}
\def\bfpp{\mathbf{p}}
\def\bfqq{\mathbf{q}}
\def\bfpq{\mathbf{q}}

\def\csqcup{\stackrel{c}{\sqcup}}

\def\cupcol{{{\stackrel{\smash{c}}{\sqcup}\,}}}
\def\cuprow{{{\stackrel{\smash{r}}{\sqcup}\,}}}

\begin{document}
	
 \title[Annihilator variety]{On the   annihilator variety of a highest weight module for
		classical Lie algebras}
 
 \author{Zhanqiang Bai}
\address[Bai]{School of Mathematical Sciences, Soochow University, Suzhou 215006, P. R. China}
\email{zqbai@suda.edu.cn}

\author{Jia-Jun Ma}
	\address[Ma]{Department of Mathematics, School of Mathematical Sciences, Xiamen University, Xiamen 361005, P. R. China}
	\address[Ma]{Department of Mathematics, Xiamen University Malaysia campus, Malaysia}
	\email{hoxide@gmail.com}
	
\author {Yutong Wang}
		\address[Wang]{School of Mathematical Sciences,  Soochow University, Suzhou 215006, P. R. China}
	\email{wyt123qq@hotmail.com}
	
	
	\date{\today}

	\begin{abstract}
		Let $\mathfrak{g}$ be a classical complex  simple Lie algebra. Let $L(\lambda)$ be a highest weight module of $\mathfrak{g}$ with highest weight $\lambda-\rho$, where $\rho$ is half the sum of positive roots. The associated variety of the annihilator ideal of $L(\lambda)$  is called  the annihilator variety of $L(\lambda)$.
  It is known that the annihilator variety of any highest weight module $L(\lambda)$ is the Zariski closure of a nilpotent orbit in $\mathfrak{g}^*$. But in general, this nilpotent orbit is not easy to describe for a given highest weight module $L(\lambda)$.
		In this paper, we will give some simple  formulas to  characterize this unique nilpotent orbit appearing in the   annihilator variety of a highest weight module for	classical Lie algebras. Our formulas are given by introducing two  algorithms, i.e., bipartition algorithm and partition algorithm. To get a special or metaplectic special partition from a domino type partition, we define the H-algorithm based on the Robinson-Schensted insertion algorithm. By using this H-algorithm, we can easily determine this nilpotent orbit from the information of $\lambda$.

		
	\end{abstract}

 \maketitle

	 \tableofcontents
	\section{Introduction}

	Let $\mathfrak{g}$ be a simple complex Lie algebra and $\mathfrak{h}$ be a Cartan subalgebra. 
 A two-sided ideal in
	$U(\mathfrak{g})$ is called primitive if it is the annihilator of an irreducible  representation of $\mathfrak{g}$. The classification of primitive ideals is a very important problem for the representation theory of semisimple complex Lie algebras. This problem is now completed thanks to the work of Duflo \cite{Du}, Borho-Jantzen \cite{BJ77}, Barbasch-Vogan \cite{BarV82,BV83} and Joseph \cite{J80-1,J80-2,J81-1}. Some  review of this problem can be found in \cite{Jo82} and \cite{Jo83ICM}.
 
 For a $U(\mathfrak{g})$-module $M$, let $I(M)=\Ann(M)$ be its annihilator ideal in  $U(\mathfrak{g})$ and $J(M)$ be the corresponding graded ideal in $ S(\mathfrak{g})=\text{gr} (U(\mathfrak{g}))$. The zero set of $J(M)$ in the dual vector space $\mathfrak{g}^*$ of $\mathfrak{g}$ is called {\it the annihilator variety} of $M$, which is also called {\it the associated variety} of $I(M)$. We usually denote it by $V(\Ann(M))$.  The study of associated varieties of primitive ideals or annihilator varieties of highest weight modules is a very important problem. 
	 Borho-Brylinski \cite{BoB1} proved that the associated variety of a primitive ideal $I$ with a fixed regular integral infinitesimal
	character is the Zariski closure of the nilpotent orbit in $\mathfrak{g}^*$ attached to $I$, via the Springer correspondence. Joseph \cite{Jo85} extended this result to a primitive ideal with a general infinitesimal  character.
	He mainly used Springer correspondence and truncated induction functor,
	which will be recalled in our paper. The study of  annihilator varieties of  irreducible representations is closely related with many research fields, such as representations of Weyl groups, Kazhdan-Lusztig cells, representations of Lie groups and $W$-algebras.  See for example \cite{BoB1,BarV82,BMSZ-typeC,BMSZ-ABCD,BMSZ-counting,GS,GSK,GT,Mc96,Lo15,LY}.

 But in general,  the nilpotent orbit appearing in the annihilator variety of a highest weight module is not  easy to describe for classical Lie algebras.
	In this paper, we want to give a simple description of the nilpotent orbit appearing in the annihilator variety of a highest weight module by generalizing the Robinson-Schensted insertion algorithm used in \cite{BaiX,BXX}. In Joseph's proof of the irreducibility of the annihilator variety of a highest weight module, this corresponding nilpotent orbit is not explicitly determined since the complexity of truncated induction for integral Weyl groups,  Springer correspondence and the numeral of a very even nilpotent orbit.  
 In our first algorithm, we will use bipartitions (i.e., some pairs of partitions) to obtain this nilpotent orbit. This algorithm will use the R-S algorithm to produce some pairs of partitions (i.e., bipartitions), which will give us the partition of the desired nilpotent orbit.
 In our second algorithm,  we will directly use partitions to obtain this nilpotent orbit. In this algorithm, we will use the R-S algorithm to produce some  partitions (Young tableaux), which will give us the partition of the desired nilpotent orbit after some collapses of partitions. 
 


 \subsection{Partition algorithm for type $A$}
	When $\lambda\in\hs$,  we write $\lambda=(\lambda_1,\dots,\lambda_n)=\sum_{i=1}^n\lam_i\ep_i$, where $\lam_i\in \mathbb{C}$ and $\{\ep_i \mid 1\leq i\leq n\}$ is the canonical basis of the Euclidean space $\mathbb{R}^n$. We associate to $ \lambda $ a set $ S(\lambda) $ of some Young tableaux as follows. Let $ \lambda_Y: \lambda_{i_1}, \lambda_{i_2}, \dots, \lambda_{i_r} $  be a maximal subsequence of $ \lambda_1,\lambda_2,\dots,\lambda_n $ such that $ \lambda_{i_k} $, $ 1\leq k\leq r $ are congruent to each other by $ \mathbb{Z} $. Then the Young tableau $P(\lambda_Y)$ associated to the subsequence $ \lambda_Y $ by using R-S  algorithm is a Young tableau in $ S(\lambda) $.

	\begin{Thm}[Theorem \ref{typeA}]
		Let $\mathfrak{g}=\mathfrak{sl}(n, \mathbb{C})$.	Suppose $\lambda\in \mathfrak{h}^*$. Then
		$$V(\Ann (L(\lambda)))=\overline{\mathcal{O}}_{p(\lambda)},$$
		where $ p(\lambda)$ is the partition of the Young tableau $$P(\lambda)=\underset{P(\lambda_{Y})\in S(\lambda)}{\stackrel{c}{\sqcup}} P(\lambda_{Y}).$$

	\end{Thm}
	Here $P{\stackrel{c}{\sqcup}}Q$ denotes a new Young tableau whose columns are the union of columns of the Young tableaux $P$ and $Q$.
	
	\begin{ex}Let $\mathfrak{g}=\mathfrak{sl}(7, \mathbb{C})$.
		Suppose $\lambda=(7, 5.5, 0.5, 6.5, 3,1,2)$. Then we can take $\lambda_{Y_1}=(7, 3,1,2)$ and $\lambda_{Y_2}=( 5.5,0.5, 6.5)$. So we get
		\[
		P(\lambda_{Y_1})=
		\begin{tikzpicture}[scale=\domscale+0.1,baseline=-28pt]
		\hobox{0}{0}{1}
		\hobox{1}{0}{2}
		\hobox{0}{1}{3}
		\hobox{0}{2}{7}
		\end{tikzpicture}\quad \text{and }\quad
		P(\lambda_{Y_2})=
		\begin{tikzpicture}[scale=\domscale+0.1,baseline=-20pt]
		\hobox{0}{0}{0.5}
		\hobox{1}{0}{6.5}
		\hobox{0}{1}{5.5}
		\end{tikzpicture}.
		\]
		Therefore, $p(\lambda_{Y_1})=[2, 1, 1]$,  $p(\lambda_{Y_2})=[2,1]$ and $p(\lambda)=p(\lambda_{Y_1}){\stackrel{c}{\sqcup}}p(\lambda_{Y_2}) =[4,2,1]$. So the annihilator variety is $V(\Ann (L(\lambda)))=\overline{\mathcal{O}}_{[4,2,1]}.$
	\end{ex}
	
\subsection{Bipartition algorithm}	
	
	For the types $B_n, C_n$ and $D_n$, let
	$[\lambda] $  be the set of  maximal subsequences $ x $ of  $ \lambda $ such that any two entries of $ x $ have an integral  difference or sum. In this case, we set $ [\lambda]_1 $ (resp. $ [\lambda]_2 $) to be the subset of $ [\lambda] $ consisting of sequences with  all entries belonging to $ \mathbb{Z} $ (resp. $ \frac12+\mathbb{Z} $).
	We set $[\lambda]_{1,2}=[\lambda]_1\cup [\lambda]_2, \quad [\lambda]_3=[\lambda]\setminus[\lambda]_{1,2}$. 	Since there is at most one element in $[\lambda]_1 $ and $[\lambda]_2 $, we denote them by  $(\lambda)_{(0)}$ and $(\lambda)_{(\frac{1}{2})}$.

	Let  $ x=(\lambda_{i_1}, \lambda_{i_2},\dots \lambda_{i_r})\in[\lambda] $. We define
	$${x}^-=(\lambda_{i_1},\lambda_{i_2},\dots,\lambda_{i_{r-1}}, \lambda_{i_{r}},-\lambda_{i_{r}},-\lambda_{i_{r-1}},\dots,-\lambda_{i_{2}},-\lambda_{i_{1}}).$$


	Let  $ x=(\lambda_{i_1}, \lambda_{i_2},\dots \lambda_{i_r})\in[\lambda]_3 $. Let $  y=(\lambda_{j_1}, \lambda_{j_2},\dots, \lambda_{j_p}) $ be the maximal subsequence of $ x $ such that $ j_1=i_1 $ and the difference of any two entries of $ y$ is an integer. Let $ z= (\lambda_{k_1}, \lambda_{k_2},\dots, \lambda_{k_q}) $ be the subsequence obtained by deleting $ y$ from $ x $, which is possible empty.
	Define
	$$  \tilde{x}=(\lambda_{j_1}, \lambda_{j_2},\dots, \lambda_{j_p}, -\lambda_{k_q}, -\lambda_{k_{q-1}},\dots ,-\lambda_{k_1}).  $$
	
	Note that from   $ x=(\lambda_{i_1}, \lambda_{i_2},\dots \lambda_{i_r})\in[\lambda]_{1,2} $, we have a Young tableau $P(x^-)$ by using the R-S  algorithm.	From $P(x^-)$, we have a partition  $p(x^-)$, which will give us a $B$-symbol (or $D$-symbol), see \S \ref{Weyl}. A $B$-symbol (or $D$-symbol) will give us a pair of partitions (i.e., a bipartition) of type $B$ (resp. type $D$), see \S \ref{j-induction}. From a bipartition, we can get a partition for a nilpotent orbit. The detailed process can be found in Section \ref{BC}. By using the induction operator $j$, from a partition $p(\tilde{x})$  coming from the Young tableau $P(\tilde{x})$ for $x\in[\lambda]_{3}$, we can get a bipartition of type $B$ (or type $D$). See equation \ref{A-B} and \ref{A-D}.
	
	We call the above process {\it bipartition algorithm}.

	\begin{Thm}[Theorem \ref{type-B}]
		Let $\mathfrak{g}=\mathfrak{so}(2n+1, \mathbb{C})$.	Suppose $\lambda\in \mathfrak{h}^*$ and $[\lambda]=(\lambda)_{(0)}\cup (\lambda)_{(\frac{1}{2})}\cup [\lambda]_3$ with $[\lambda]_3=\{{\lambda}_{{Y}_1},\dots,{\lambda}_{{Y}_m}\}$.  Then
		$$V(\Ann (L(\lambda)))=\overline{\mathcal{O}}_{p_{_B}(\lambda)},$$
		where $ p_{_B}(\lambda)$ is the partition obtained from the following bipartition $$({\bf d}_{0}{\stackrel{c}{\sqcup}}{\bf d}_{00}{\stackrel{c}{\sqcup}}_i{\bf d}_i, {\bf f}_{0}{\stackrel{c}{\sqcup}}{\bf f}_{00}{\stackrel{c}{\sqcup}}_i{\bf f}_i ).$$
		Here $ ({\bf d}_{0},{\bf f}_{0})$ is the $B$-type bipartition obtained from $p((\lambda)_{(0)}^-)$, $({\bf d}_{00},{\bf f}_{00})$ is the $B$-type bipartition obtained from $p((\lambda)_{(\frac{1}{2})}^-)$ and $ ({\bf d}_{i},{\bf f}_{i})$ is the $B$-type bipartition obtained from the partition $p(\tilde{\lambda}_{{Y}_i})$.

	\end{Thm}
	From $(\lambda)_{(\frac{1}{2})}$, we have a partition $p((\lambda)_{(\frac{1}{2})}^-)$, which will give us a special partition of type $D$. This special type $D$ partition will give us a $C$-type metaplectic special partition by computing its dual of $D$-collapse of its dual partition. See Proposition \ref{D-C}. This $C$-type metaplectic special partition will give us a bipartition called  $C$-type bipartition. The definition of  collapse of partitions will be given in \S \ref{j-induction}.
	
	\begin{Thm}[Theorem \ref{type-C}]
		Let $\mathfrak{g}=\mathfrak{sp}(n, \mathbb{C})$.	Suppose $\lambda\in \mathfrak{h}^*$ and $[\lambda]=(\lambda)_{(0)}\cup (\lambda)_{(\frac{1}{2})}\cup [\lambda]_3$ with $[\lambda]_3=\{{\lambda}_{{Y}_1},\dots,{\lambda}_{{Y}_m}\}$.  Then
		$$V(\Ann (L(\lambda)))=\overline{\mathcal{O}}_{p_{_C}(\lambda)},$$
		where $ p_{_C}(\lambda)$ is the $C$-type partition obtained from the following bipartition $$({\bf d}_{0}{\stackrel{c}{\sqcup}}{\bf d}_{00}{\stackrel{c}{\sqcup}}_i{\bf d}_i, {\bf f}_{0}{\stackrel{c}{\sqcup}}{\bf f}_{00}{\stackrel{c}{\sqcup}}_i{\bf f}_i ).$$
		Here $ ({\bf d}_{0},{\bf f}_{0})$ is the $B$-type bipartition obtained from $p((\lambda)_{(0)}^-)$, $({\bf d}_{00},{\bf f}_{00})$ is the $C$-type bipartition obtained from $p((\lambda)_{(\frac{1}{2})}^-)$ and $ ({\bf d}_{i},{\bf f}_{i})$ is the $B$-type bipartition obtained from the partition $p(\tilde{\lambda}_{{Y}_i})$.

	\end{Thm}

	\begin{ex}\label{Cexa}
		Let $\mathfrak{g}=\mathfrak{sp}(15, \mathbb{C})$.
		Suppose $$\lambda=(-3.5,2.5,1.5,3,6,9,-7,2,5,-8,-4,1,0.6,2.6,-0.6).$$ Then we can take $$(\lambda)_{(0)}=(3,6,9,-7,2,5,-8,-4,1),$$  $$(\lambda)_{(\frac{1}{2})}=(-3.5,2.5,1.5),$$ and $$\lambda_{Y_1}=(0.6,2.6,-0.6).$$
		
		For $(\lambda)_{(0)}$, the Young tableau $T((\lambda)_{(0)}^-)$ will give us a partition $p((\lambda)_{(0)}^-)=[5,4,3,3,2,1]$. Then we
		have a  strictly increasing sequence $$
		(0,2,4,6,7,9,11)=(0,1,2,3,3,4,5)+(0,1,2,3,4,5,6).$$
		It will give us a
		$B$-symbol $\Lambda= \begin{pmatrix}
		0~1~2~3\\3~4~5
		\end{pmatrix} $. The corresponding special symbol is $\Lambda^s=\begin{pmatrix}
		0~2~3~5\\1~3~4
		\end{pmatrix}$.  The corresponding bipartition is $({\bf d}_{0},{\bf f}_{0})=([2,1,1],[2,2,1])$.

		For $(\lambda)_{(\frac{1}{2})}$,  the Young tableau $P((\lambda)_{(\frac{1}{2})}^-)$ will give us a partition $p((\lambda)_{(\frac{1}{2})}^-)=[3,1,1,1]$. Then we
		have a  strictly increasing sequence $$
		(0,2,3,4,7)=(0,1,1,1,3)+(0,1,2,3,4).$$
		It will give us a
		$B$-symbol $\Lambda= \begin{pmatrix}
		0~1~2\\1~3
		\end{pmatrix} $, and a $D$-symbol $\Lambda= \begin{pmatrix}
		0~1~2\\0~2~4		\end{pmatrix} \sim \begin{pmatrix}
		0~1\\1~3		\end{pmatrix} $.  The corresponding partition of type $D$ is the partition ${\bf d}_b=[3,1,1,1]$.
		From the  $C$-type metaplectic special partition $(({\bf d}^*_b)_D)^*=[4,1,1]$, we can get a $C$-type bipartition $({\bf d}_{00},{\bf f}_{00})=([2,1],[0])$.
		
		For $\lambda_{Y_1}$, we have $\tilde{\lambda}_{Y_1}=(0.6,2.6,0.6)$.  The Young tableau $P(\tilde{\lambda}_{Y_1})$ will give us a partition $p(\tilde{\lambda}_{Y_1})=[2,1]$. Its dual is $p(\tilde{\lambda}_{Y_1})^*=[2,1]$. This will give us a $B$-type bipartition $ ({\bf d}_{1},{\bf f}_{1})=([2],[1])$.

		Thus we can get a bipartition $({\bf d}_{0}{\stackrel{c}{\sqcup}}{\bf d}_{00}{\stackrel{c}{\sqcup}}_i{\bf d}_i, {\bf f}_{0}{\stackrel{c}{\sqcup}}{\bf f}_{00}{\stackrel{c}{\sqcup}}_i{\bf f}_i )=([6,2,1],[3,2,1])$.
		The corresponding  $C$-type partition is ${\bf d}=[12,6,4,4,2,2]$. Thus we have $V(\Ann (L(\lambda)))=\overline{\mathcal{O}}_{[12,6,4,4,2,2]}$.
	\end{ex}
	

	\begin{Thm}[Theorem \ref{typeD}]
  Let $\mathfrak{g}=\mathfrak{so}(2n, \mathbb{C})$.	Suppose $\lambda\in \mathfrak{h}^*$ and $[\lambda]=(\lambda)_{(0)}\cup (\lambda)_{(\frac{1}{2})}\cup [\lambda]_3$ with $[\lambda]_3=\{{\lambda}_{{Y}_1},\dots,{\lambda}_{{Y}_m}\}$.  Then
		$$V(\Ann (L(\lambda)))=\overline{\mathcal{O}}_{p_{_D}(\lambda)},$$
		where $ p_{_D}(\lambda)$ is the partition obtained from the following unordered bipartition $$\{{\bf d}_{0}{\stackrel{c}{\sqcup}}{\bf d}_{00}{\stackrel{c}{\sqcup}}_i{\bf d}_i, {\bf f}_{0}{\stackrel{c}{\sqcup}}{\bf f}_{00}{\stackrel{c}{\sqcup}}_i{\bf f}_i \}.$$
		Here $ \{{\bf d}_{0},{\bf f}_{0}\}$ is the $D$-type unordered bipartition obtained from $(\lambda)_{(0)}$, $\{{\bf d}_{00},{\bf f}_{00}\}$ is the $D$-type unordered bipartition obtained from $(\lambda)_{(\frac{1}{2})}$ and $ \{{\bf d}_{i},{\bf f}_{i}\}$ is the $D$-type unordered bipartition obtained from $\tilde{\lambda}_{{Y}_i}$. 
  When $n=2m'$ is even and ${\mathcal{O}}_{p_{_D}(\lambda)}$ is a very even orbit, it will be type I  if $k(\lambda)\equiv 0 
  ~(\mathrm{mod}~ {2})$ and type II if $k(\lambda)\equiv 1 
  ~(\mathrm{mod}~ {2})$. Here
	 we use $k(\lambda)$ to denote the number of very even unordered bipartitions with numeral II in the set of  $\{\{{\bf d}_{0},{\bf f}_{0}\}, \{{\bf d}_{00},{\bf f}_{00}\},  \{{\bf d}_{i},{\bf f}_{i}\}\mid  1\leq i\leq m'  \}$. 
	\end{Thm}

 \subsection{Partition algorithm}
	By using Sommers's recipe \cite{Sommers}, we also have another algorithm, called {\it partition algorithm}. We use $X$ to denote the corresponding type of Lie algebras, i.e., $X=B,C$ or $D$.

	\begin{Thm}[Theorem \ref{partition-algorithm-} and Corollary \ref{nemural-d}]
		Suppose $\mathfrak{g} = \mathfrak{so}(2n+1, \mathbb{C}), \mathfrak{sp}(n, \mathbb{C})$ or $\mathfrak{so}(2n, \mathbb{C})$, $\lambda\in \mathfrak{h}^*$ and
		$[\lambda]=(\lambda)_{(0)} \cup (\lambda)_{(\frac{1}{2})}\cup [\lambda]_3$
		with $[\lambda]_3=\{{\lambda}_{{Y}_1},\dots,{\lambda}_{{Y}_m}\}$. Let
		\begin{enumerate}
			\item $\mathbf{p}_{0}$ be the $X$-type special partition associated to
			$(\lambda)_{(0)}$;
			\item ${\mathbf p}_{\frac{1}{2}}$ be the $C$-type special  partition (for type $B$) or $C$-type metaplectic special partition (for types $C$ and $D$) associated to
			$(\lambda)_{(\frac{1}{2})}$;
			\item $\mathbf{p}_{i}$ be the $A$-type partition associated to
			$\tilde{\lambda}_{Y_i}$.
		\end{enumerate}
		Let $\mathbf{p}_{\lambda}$ be the $X$-collapse of
		\begin{equation}
		\mathbf{d}_{\lambda} := \mathbf{p}_{0} {\stackrel{c}{\sqcup}} {\mathbf p}_{\frac{1}{2}} {\stackrel{c}{\sqcup}}({{\stackrel{c}{\sqcup}}}_i 2 \mathbf{p}_i).
		\end{equation}
		Then we have
		\[
		V(\Ann (L(\lambda)))=\overline{\mathcal{O}}_{\mathbf{p}_{\lambda}}.
		\]
 When $n=2m'$ is even and ${\mathcal{O}}_{\mathbf{p}_{\lambda}}$ is a very even orbit, its numeral is determined by the number $k(\lambda)$ of very even orbits with numeral II in the set of very even orbits of type $D$:    $\{\mathcal{O}_{\bfpp_{0}}, \mathcal{O}_ {{\bf p}_{\frac{1}{2}}},  \mathcal{O}_{2\bfpp_i}| 1\leq i\leq m  \}$. 
So $\mathcal{O}_{{\bf p}_{\lambda}}$ will be type I if $k(\lambda)\equiv 0 
  ~(\mathrm{mod}~ {2})$ and type II if $k(\lambda)\equiv 1 
  ~(\mathrm{mod}~ {2})$.
 \end{Thm}	
	
	In  the beginning, $\mathbf{p}_{0}$ and ${\mathbf p}_{\frac{1}{2}}$ are obtained by using special symbols. By inspecting the algorithms used in \cite{BXX}, we find some new algorithms without using symbols in \S \ref{refi}. These algorithms are called {\it Hollow diagram algorithms} (or simply H-algorithms).  So the partition type of this nilpotent orbit ${\mathcal{O}}_{\mathbf{p}_{\lambda}}$ can be determined very easily by successively using R-S algorithm,  H-algorithm and collapse of partitions.

	
	

	\begin{ex}	In the above Example \ref{Cexa},
		for $\mathfrak{g}=\mathfrak{sp}(15, \mathbb{C})$ and $$\lambda=(-3.5,2.5,1.5,3,6,9,-7,2,5,-8,-4,1,0.6,2.6,-0.6).$$ Then from $(\lambda)_{(0)}=(3,6,9,-7,2,5,-8,-4,1),$
		we have $p((\lambda)_{(0)}^-)=[5,4,3,3,2,1]$, which is not a special partition of type $C$. Then we
		have a  strictly increasing sequence $$
		(0,2,4,6,7,9,11)=(0,1,2,3,3,4,5)+(0,1,2,3,4,5,6).$$
		It will give us a
		$B$-symbol $\Lambda= \begin{pmatrix}
		0~1~2~3\\3~4~5
		\end{pmatrix} $. The corresponding special symbol is $\Lambda^s=\begin{pmatrix}
		0~2~3~5\\1~3~4
		\end{pmatrix}$. The corresponding special partition of type $C$ is the partition ${\bf p}_0=[4,4,3,3,2,2]$.  If we use the H-algorithm of type $C$,  we will have

\[p((\lambda)_{(0)}^-)=
\tytb{EOEOE,OEOE,EOE,OEO,EO,O} \to
\tytb{{\none[1]}\none O\none O\none,{\none[2]}O\none O,\none\none O,\none O\none O,{\none[3]}\none O,{\none[4]}O}
\to
\tytb{{\none[1]}\none O\none O\none,{\none[2]}O\none OE,\none\none O,\none O\none O,{\none[3]}\none O,{\none[4]}OE}
\to
\tytb{{\none[1]}EOEO,{\none[2]}OEOE,\none EOE,\none OEO,{\none[3]}EO,{\none[4]}OE}={\bf p}_0.
\]
Thus we can also get the same partition ${\bf p}_0=[4,4,3,3,2,2]$.

		From $(\lambda)_{(\frac{1}{2})}=(-3.5,2.5,1.5)$,  the Young tableau $P((\lambda)_{(\frac{1}{2})}^-)$ will give us a partition $p((\lambda)_{(\frac{1}{2})}^-)=[3,1,1,1]$. Then we
		have a  strictly increasing sequence $$
		(0,2,3,4,7)=(0,1,1,1,3)+(0,1,2,3,4).$$
		It will give us a
		$B$-symbol $\Lambda= \begin{pmatrix}
		0~1~2\\1~3
		\end{pmatrix} $, and a special $D$-symbol $\Lambda= \begin{pmatrix}
		0~1~2\\0~2~4		\end{pmatrix} \sim \begin{pmatrix}
		0~1\\1~3		\end{pmatrix} $.  The corresponding special partition of type $D$ is the partition ${\bf d}_b=[3,1,1,1]$.
		So the $C$-type metaplectic special partition is ${\bf p}_{\frac{1}{2}}=(({\bf d}^*_b)_D)^*=[4,1,1]$. If we use the H-algorithm of metaplectic type, we will have
  \begin{align*}
    p((\lambda)_{(\frac{1}{2})}^-)=\tytb{EOE,O,E,O} \to
    \tytb{{\none[1]}E\none E, \none,  \none E} \to
\tytb{{\none[1]}E\none EO, \none,  \none E}
\to
\tytb{{\none[1]}EO EO, \none O,  \none E}
={\bf p}_{\frac{1}{2}}.
\end{align*}
  Thus    we can also get the same  partition $[4,1,1]$.
		
		For $\lambda_{Y_1}=(0.6,2.6,-0.6)$, we have $\tilde{\lambda}_{Y_1}=(0.6,2.6,0.6)$.  The Young tableau $P(\tilde{\lambda}_{Y_1})$ will give us a partition $\mathbf{p}_{1}=p(\tilde{\lambda}_{Y_1})=[2,1]$.
		
		Thus
		\begin{equation*}
		\mathbf{d}_{\lambda} = \mathbf{p}_{0} {\stackrel{c}{\sqcup}} {\mathbf p}_{\frac{1}{2}} {\stackrel{c}{\sqcup}}  2\mathbf{p}_1=[12,7,4,3,2,2].
		\end{equation*}
		And $\mathbf{p}_{\lambda}=(\mathbf{d}_{\lambda})_C=[12,6,4,4,2,2]$.

		Thus we have the same result  $V(\Ann (L(\lambda)))=\overline{\mathcal{O}}_{[12,6,4,4,2,2]}$.
		
	\end{ex}	
	
Finally we design the following web page for people to use our algorithm: \begin{center}
\textcolor{blue}{http://liealgebra.slashblade.top/lie/classification}.   
\end{center}

	This paper is organized as follows. In \S\ref{subsec:gkh}, we give some necessary preliminaries about annihilator varieties of highest weight modules and integral Weyl groups. In \S\ref{j-induction}, we recall some properties of the $j$-induction operator.
 In \S \ref{Weyl}, we recall the Springer correspondence between representations of Weyl groups and  nilpotent orbits.  Then the combinatorial characterizations of the annihilator varieties (nilpotent orbits) of  all the simple highest weight modules for classical Lie algebras are obtained in \S\ref{type-A}, \S\ref{BC}, \S\ref{Dchar} and \S\ref{refi}.

	\section{Notation and Preliminary Results}\label{subsec:gkh}
	Let $\mathfrak{g}$ be a simple complex Lie algebra and $\mathfrak{h}$ be a Cartan subalgebra.
	Let $\Phi^+\subset\Phi$ be the set of positive roots determined by a Borel subalgebra $\mathfrak{b}$ of $\mathfrak{g}$. So we have a Cartan decomposition $\mathfrak{g}=\mathfrak{n}\oplus\mathfrak{h}\oplus\mathfrak{n}^-$. Denote by $\Delta$ the set of simple roots in $\Phi^+$.
	Choose a subset $I\subset \Delta$ and  it generates a subsystem
	$\Phi_I\subset\Phi$.
	Let $\mathfrak{p}_I$ be the standard parabolic subalgebra corresponding to $I$ with Levi decomposition $\mathfrak{p}_I=\mathfrak{ l}_I\oplus\mathfrak{u}_I$. We frequently drop the
	subscript $I$ if there is no confusion.
	
	Let $F(\lambda)$ be a finite-dimensional irreducible $\mathfrak{l}$-module with highest weight $\lambda -\rho\in\mathfrak{h}^*$. It can also be viewed as a
	$\mathfrak{p}$-module with trivial $\mathfrak{u}$-action. The {\it generalized Verma module} $N_I(\lambda)$ is defined by
	\[
	N_I(\lambda):=U(\mathfrak{g})\otimes_{U(\mathfrak{p})}F(\lambda).
	\]
	The simple quotient of $N_I(\lambda)$ is denoted by $L(\lambda)$, which is a highest weight module with highest weight $\lambda-\rho$.

	
	
 \subsection{Annihilator varieties of highest weight modules}
	
	Let $M$ be a finite generated $U(\mathfrak{g})$-module. Fix a finite dimensional generating subspace $M_0$ of $M$. Let $U_{n}(\mathfrak{g})$ be the standard filtration of $U(\mathfrak{g})$. Set $M_n=U_n(\mathfrak{g})\cdot M_0$ and
	\(
	\text{gr} (M)=\bigoplus\limits_{n=0}^{\infty} \text{gr}_n M,
	\)
	where $\text{gr}_n M=M_n/{M_{n-1}}$. Thus $\text{gr}(M)$ is a graded module of $\text{gr}(U(\mathfrak{g}))\simeq S(\mathfrak{g})$.
	We use $\Ann (M)$ denote the two-sided ideal of $U(\mathfrak{g})$
	consisting of elements that act by zero on $M$.
	
	
	
	\begin{Defff}
		The  \textit{associated variety} of $M$ is defined by
		\begin{equation*}
		V(M):=\{X\in \mathfrak{g}^* \mid f(X)=0 \text{ for all~} f\in \operatorname{Ann}_{S(\mathfrak{g})}(\operatorname{gr} M)\}.
		\end{equation*}
	\end{Defff}
	
	The above  definition is independent of the choice of $M_0$ (e.g., \cite{NOT}).
	
	\begin{Defff} Let $\mathfrak{g}$ be a finite-dimensional semisimple Lie algebra. Let $I$ be a two-sided ideal in $U(\mathfrak{g})$. Then $\text{gr}(U(\mathfrak{g})/I)\simeq S(\mathfrak{g})/\text{gr}I$ is a graded $S(\mathfrak{g})$-module. Its annihilator is $\text{gr}I$. We define its associated variety by
		$$V(I):=V(U(\mathfrak{g})/I)=\{X\in \mathfrak{g}^* \mid p(X)=0\ \mbox{for all $p\in {\text{gr}}I$}\}.
		$$
	\end{Defff}
	
	Following \cite{GSK}, $V(\Ann (M))$ is called the \textit{annihilator variety} of the $U(\mathfrak{g})$-module $M$.


	We use $\stackrel{R}{\sim}$ (resp. $\stackrel{L}{\sim}$ and $\stackrel{LR}{\sim}$ ) to denote the Kazhdan-Lusztig right (resp. left and double) cell equivalence relation \cite{KL}.

	Let $G$ be a connected semisimple finite dimensional complex algebraic group with Lie algebra $\mathfrak{g}$. Let $W$ be the Weyl group of $\mathfrak{g}$.  We use $L_w$ to denote the simple highest weight $\mathfrak{g}$-module of highest weight $-w\rho-\rho$ with  $w\in W$.  We denote $I_w=\Ann(L_w)$.  Borho-Brylinski \cite{BoB1} proved that the annihilator variety of $L_w$ is irreducible, i.e., it is the closure of a single nilpotent orbit. We denote   $V(I_w)= V(\Ann (L_w))=\overline{\mathcal{O}}_w$.

	From \cite{BoB1} or \cite{Ta}, we know that there exists a bijection between the double cells in the Weyl group $W$ and special nilpotent orbits. In other words, $\mathcal{O}_w=\mathcal{O}_y$ if and only if $w\stackrel{LR}{\sim}y$.




By using Springer's correspondence and $j$-induction operator, Joseph extended the result in \cite{BoB1} to arbitrary infinitesimal character.	We have the following result.
	
	\begin{Prop}[ \cite{Jo85}]
		Let $\mathfrak{g}$ be a reductive Lie algebra and $I$ be a primitive ideal in $U(\mathfrak{g})$.Then $V(I)$ is the closure of a single nilpotent coadjoint orbit $\mathcal{O}_I$ in $\mathfrak{g}^*$. In particular, for a highest weight module $L(\lambda)$, we denote $V(\Ann (L(\lambda)))=\overline{\mathcal{O}}_{Ann(L(\lambda))}$.
	\end{Prop}

\subsection{Robinson-Schensted insertion algorithm}\label{R-S} We recall the Robinson-Schensted insertion algorithm which will be used in our paper. Some details can be found in \cite{BaiX}.

For  a totally ordered set $ \Gamma $, we  denote by $ \mathrm{Seq}_n (\Gamma)$ the set of sequences $ x=(x_1,x_2,\dots, x_n) $   of length $ n $ with $ x_i\in\Gamma $. In our paper, we usually take $\Gamma$ to be $\mathbb{Z}$ or a coset of $\mathbb{Z}$ in $\mathbb{C}$.
	Then we have a  Young tableau $P(x)$ obtained by applying the following Robinson-Schensted insertion algorithm  to $x\in \mathrm{Seq}_n (\Gamma)$. 
 \begin{Defff}[Robinson-Schensted insertion algorithm]
For an element  $ x \in  \mathrm{Seq}_n (\Gamma)$, we write  $x=(x_1,\dots,x_n)$. We associate to $x $ a  Young tableau  $ P(x) $ as follows. Let $ P_0 $ be an empty Young tableau. Assume that we have constructed Young tableau $ P_k $ associated to $ (x_1,\dots,x_k) $, $ 0\leq k<n $. Then $ P_{k+1} $ is obtained by adding $ x_{k+1} $ to $ P_k $ as follows. Firstly we add $ x_{k+1} $ to the first row of $ P_k $ by replacing the leftmost entry $ x_j $ in the first row which is \textit{strictly} bigger than $ x_{k+1} $.  (If there is no such an entry $ x_j $, we just add a box with entry $x_{k+1}  $ to the right side of the first row, and end this process). Then add this $ x_j $ to the next row as the same way of adding $x_{k+1} $ to the first row.  Finally we put $P(x)=P_n$.

\end{Defff}

We use $p(x)=(p_1,\dots, p_k)$ to denote the shape of $P(x)$, where $p_i$ is the number of boxes in the $i$-th row of  $P(x)$.
When $\sum\limits_{1\leq i\leq k} p_i=N$, $p(x)$ will be a partition of $N$ and we still denote this partition by $p(x)=[p_1,\dots, p_k]$.

In general, Robinson-Schensted insertion algorithm is abbreviated to R-S algorithm. 

\begin{ex}
    Suppose $x=(-5,-6,-4,2, -2,-4)$.  Usually  we write $x=(-5,-6,-4,2, -2,-4')$ and regard $-4<-4'$. Then from the R-S algorithm, we have
$$
\begin{tikzpicture}[scale=\domscale+0.1,baseline=-13pt]
		\hobox{0}{0}{-5}
		\end{tikzpicture}\stackrel{-6}{\to} \begin{tikzpicture}[scale=\domscale+0.1,baseline=-18pt]
		\hobox{0}{0}{-6}
		\hobox{0}{1}{-5}
		\end{tikzpicture}\stackrel{-4}{\to} \begin{tikzpicture}[scale=\domscale+0.1,baseline=-18pt]
		\hobox{0}{0}{-6}
		\hobox{1}{0}{-4}
		\hobox{0}{1}{-5}
				\end{tikzpicture}\stackrel{2}{\to}
    \begin{tikzpicture}[scale=\domscale+0.1,baseline=-18pt]
		\hobox{0}{0}{-6}
		\hobox{1}{0}{-4}
  \hobox{2}{0}{2}
		\hobox{0}{1}{-5}
		\end{tikzpicture}\stackrel{-2}{\to}
  \begin{tikzpicture}[scale=\domscale+0.1,baseline=-18pt]
		\hobox{0}{0}{-6}
		\hobox{1}{0}{-4}
  \hobox{2}{0}{-2}
		\hobox{0}{1}{-5}
  \hobox{1}{1}{2}
		\end{tikzpicture}\stackrel{-4'}{\to}
   \begin{tikzpicture}[scale=\domscale+0.1,baseline=-18pt]
		\hobox{0}{0}{-6}
		\hobox{1}{0}{-4}
  \hobox{2}{0}{-4'}
		\hobox{0}{1}{-5}
  \hobox{1}{1}{-2}
  \hobox{0}{2}{2}
		\end{tikzpicture}=P(x).$$
 Thus we have $p(x)=[3,2,1]$, which is a partition of $6$.
           
           

\end{ex}

\subsection{Integral Weyl group}\label{integral2}
	Let $(-, -)$ be the standard bilinear form on $\mathfrak{h}^*$. For $\lambda\in\mathfrak{h}^*$, define
	\begin{equation*}
	\Phi_{[\lambda]}:=\{\alpha\in\Phi\mid\langle\lambda, \alpha^\vee\rangle\in\mathbb{Z}\},
	\end{equation*}
	where $\langle\lambda, \alpha^\vee\rangle=2(\lambda, \alpha)/(\alpha, \alpha)$. Set
	\[
	W_{[\lambda]}:=\{w\in W\mid w\lambda-\lambda\in \mathbb{Z}\Phi\}.
	\]
	Then $\Phi_{[\lambda]}$ is a root system with Weyl group $W_{[\lambda]}$
	(e.g., \cite{Hum08}).
	Let $\Delta_{[\lambda]}$ be the simple system of $\Phi_{[\lambda]}$.  Set $J=\{\alpha\in\Delta_{[\lambda]}\mid\langle\lambda, \alpha^\vee\rangle=0\}$. Denote by $W_J$ the Weyl group generated by reflections $s_\alpha$ with $\alpha\in J$. Let $\ell_{[\lambda]}$ be the length function on $W_{[\lambda]}$. Thus $\ell_{[\lambda]}=\ell$ when $\lambda$ is integral. Put
	\begin{equation*}\label{ceq1}
	W_{[\lambda]}^J:=\{w\in W_{[\lambda]}\mid \ell_{[\lambda]}(ws_\alpha)=\ell_{[\lambda]}(w)+1\ \mbox{for all}\ \alpha\in J\}.
	\end{equation*}
	Thus $W_{[\lambda]}^J$ consists of the shortest representatives of the cosets $wW_J$ with $ w\in W_{[\lambda]} $. When $\lambda$ is integral, we simply write $W^J:=W_{[\lambda]}^J$ .
	
	A weight $ \mu\in\hs $ is called \textit{anti-dominant} if $ \bil{\mu}{\al} \notin\mathbb{Z}_{>0}$ for all $ \al\in\Phi^+ $. For any $\lambda\in\mathfrak{h}^*$, there exists a unique anti-dominant weight $\mu\in\hs$ and a unique $w\in W_{[\mu]}^J$ such that $\lambda=w\mu$.
\begin{Prop}[{\cite[Proposition 3.5]{Hum08}}]\label{anti}
    Let $\lambda\in \mathfrak{h}^*$, with corresponding root system $\Phi_{[\lambda]}$ and Weyl group $W_{[\lambda]}$. Let $\Delta_{[\lambda]}$ be the simple system of $\Phi_{[\lambda]} \cap \Phi^+$  in $\Phi_{[\lambda]}$.  Then $\lambda$ is antidominant if and only if one of the
following three equivalent conditions holds:
\begin{enumerate}
    \item $\langle\lambda, \alpha^\vee\rangle\leq 0$ for all $\alpha \in \Delta_{[\lambda]}$;
    \item $\lambda\leq s_{\alpha}\lambda$ for all $\alpha \in \Delta_{[\lambda]}$;
    \item  $\lambda\leq w\lambda$ for all $w \in W_{[\lambda]}$.
\end{enumerate}
Therefore there is a unique antidominant weight in the orbit $W_{[\lambda]}\lambda$.
\end{Prop}

	Let $ \lambda$ be a regular element
	in $\mathfrak{h}^*$ and anti-dominant with respect to $\Phi^+ $. Then $\lambda$ is
	dominant with respect to $\Phi_{[\lambda]}^+$.
	For each $w\in  W_{[\lambda]}$, we define
	$$a(w)= |\Phi^+|-\gkd L(w\lambda).$$
	
	Let $Q$ be the $\mathfrak{h}$ root lattice of $\mathfrak{g}$. We define $$[\lambda]_{R}=\lambda+Q\in \mathfrak{h}^*.$$
	For
	each $w\in W$, we can attach a polynomial $p_w$ such that $$p_w(\mu)=\mathrm{rank}(U(\mathfrak{g})/\Ann (L(w\mu)) )$$ where $\mu\in[\lambda]_R$
	is dominant. $p_w$  is called the {\it Goldie-rank
		polynomial} attached to the primitive ideal $\Ann (L(w\lambda))$.
	
	We take a double cell $\mathcal{C}^{LR}$  in $W_{[\lambda]}$ and a set of representatives $\{w_1, w_2, \dots , w_k \}$ of the left cells in $\mathcal{C}^{LR}$.
	Then we have the following proposition.
	
	\begin{Prop}[{\cite{J80-1,J80-2,J81-1,BarV82}}]\label{intehral-W}
		The set $\{p_{w_i}| 1\leq i\leq k\}$  forms a basis of a special representation $\pi_w$ of $W_{[\lambda]}$  realized in $S^{a(\pi_w)}(\mathfrak{h})$. Here $a(\pi_w)=a(w)$ is the minimal degree $m$ such that $\pi_w$ occurs in $S^{m}(\mathfrak{h})$, which is called the fake degree
		of $\pi_w$.	There is a one-to-one correspondence between the set of double cells in $W_{[\lambda]}$
		and the set  of special representations of $W_{[\lambda]}$.
	\end{Prop}
	
	If we use $\mathcal{C}^{LR}_w$ to denote the double cell in $W_{[\lambda]}$ containing the element $w$, then the correspondence in the above proposition is mapping $\mathcal{C}^{LR}_w$ to $\pi_w$.

\section{Representations of classical Weyl groups and the $j$-induction}\label{j-induction}
In this section, we recall some properties of the $j$-induction operator. Some details can be found in \cite{Ca85,Lus84,BMSZ-counting}.

Suppose $W$ is a Weyl group attached to a root system $\Phi$ (with a fixed subset of simple roots $\Delta$) in a finite dimensional real vector space $V$. Then $W$ acts on the space $P_k(V)$ of degree $k$ homogeneous polynomials on $V$.  
After Joseph \cite{J80-2}, a representation $\sigma\in \Irr(W)$ is called \emph{univalent} if it
occurs with multiplicity one in  $P_{b(\sigma)}(V)$ where $b(\sigma)$ is the minimal degree  such that $\sigma$ occurs in $P_{b(\sigma)}(V)$.
The number $b(\sigma)$ is called the \emph{fake degree} of $\sigma$.
Let $\Univ(W)\subset \Irr(W)$ be the set of univalent representations. 
Under the Springer correspondence, the image of the trivial local system on a nilpotent orbit is always univalent \cite[Corollarie~4]{BM}. 

Suppose $W'$ is a subgroup of $W$ generated by reflections in a subroot system of $\Phi$.   
The $j$-induction from $W'$ to $W$ is a well defined map 
\[
\begin{array}{rcl}
j_{W'}^{W} \colon \Univ(W')& \rightarrow &\Univ(W) \\
\sigma' & \mapsto & j_{W'}^{W}\sigma'
\end{array}
\]
where $j_{W'}^{W}\sigma'$ is the representation generated by the $\sigma'$ 
isotypical component in $P_{b(\sigma')}(V)$. Furthermore, $b(j_{W'}^{W}\sigma') =b(\sigma')$ and $j_{W'}^{W}\sigma'$ occurs in $P_{b(\sigma')}(V)$ with multiplicity one.  
The definitions of ``univalent'' and ``$j$-induction'' are independent of the choice of $V$. 

The $j$-induction satisfies the induction by stage. 
\begin{Lem}[{See \cite[Theorem~11.2.4]{Ca85}}]\label{j.1}
Let $W''\subset W'$ be two subgroups of $W$ generated by two sub root systems (not necessary to be parabolic subgroups). 
Then 
\[
j^W_{W''}=j^W_{W'}\circ j^{W'}_{W''}.
\]
\end{Lem}
Some details can be found in \cite[\S 11.2]{Ca85} and \cite[Chapter~4]{Lus84}.

 
	Let $\lambda$ be a regular element
	in $\mathfrak{h}^*$ and anti-dominant with respect to $\Phi^+ $. Then $\lambda$ is
	dominant with respect to $\Phi_{[\lambda]}^+$.
 From a representation of the integral Weyl group $W_{[\lambda]}$, we can get a representation of the Weyl group $W$ by using the $j$-induction operator. 

 We recall the construction of $\pi_w$ in \S \ref{integral2}, which is a representation of $W_{[\lambda]}$. Then we have the following result.
\begin{Prop}[{\cite[Theorem~3.10]{Jo85}}]\label{AV}
    Let $w\in W_{[\lambda]}$. 
    Under the $W$ action, $\pi_w$ generates an
    irreducible $W$-module, with the same fake degree $a(w)$.  This irreducible representation of $W$ is called the $j$-induction (or  truncated induction) of $\pi_w$, to be denoted by $\tilde{\pi}_w=j_{W_{[\lambda]}}^W(\pi_w)$. This $W$-module $\tilde{\pi}_w$ corresponds to  a nilpotent orbit $\mathcal{O}_{\tilde{\pi}_w}$
    with trivial local system via the Springer correspondence. Furthermore,
$$V(\Ann (L(w\lambda)))=\overline{\mathcal{O}}_{\tilde{\pi}_w}.$$
\end{Prop}

\subsection{Representations of $S_n$}
    We view $S_n$ as the Weyl group attached to the root system of type $A_{n-1}$ and let $\sgn$ denote the sign representation of $S_n$.
	For $W_n$, let $\varepsilon$ denote the unique non-trivial character, which is trivial on $S_n$. Note that $\varepsilon$ is also the restriction of the sign representation of $S_{2n}$ on $W_n$.
	
    We use a multiset $\bfdd=\set{d_1,d_2,\dots,d_k}$ of positive integers to denote a partition of $|{\bfdd}|:=d_1+\dots+d_k$. Alternatively, for a sequence $d_1\geq d_2\geq \dots\geq d_k$ of non-negative integers,  
    write $\bfdd=[d_1,d_2,\dots,d_k]$ for the partition $\set{d_i| 0<d_i, 0\leq i \leq k}$. 
    A partition is identified with the Young diagram such that $\bfdd$ is the multi-set of the lengths of non-empty rows in the Young diagram.  By abuse of notation, 
    we also use $\bfdd$ to denote the Young diagram corresponding to the partition $\bfdd$.
    As usual ${\bf d}^*=\set{d^*_1,d^*_2,\dots,d^*_q}$ is the dual partition of $\bfdd$ corresponding to the transpose of the Young diagram $\bfdd$.
    Suppose $\bfff:=[f_1,f_2,\dots, f_l]$ is a partition, let $\bfdd \cuprow \bfff$ denote that partition $\bfdd \cup \bfff$ and $\bfdd\cupcol \bfff$ denote the partition given by $(\bfdd\cupcol\bfff)^* = \bfdd^*\cuprow \bfff^*$.

From \cite[Theorem 10.1.1]{CM},  we  can identify the set of Young diagrams of total size $n$ with the set of complex nilpotent orbits such that each row corresponds to a Jordan block.  
    Consequently,  the set of Young diagrams is also identified with the set  $\Irr(S_{n})$ via the Springer correspondence. We adopt the normalization of the correspondence such that the trivial orbit (the Young diagram having a single column) corresponds to the sign representation of $S_n$. 
 
    We recall  the the Springer
	correspondence in terms of Macdonald's construction of irreducible representations of $S_{n}$ via $j$-induction in \cite{Mc72}: 
 Let  $\cO_{\bfdd}$ be a nilpotent orbit corresponds to the partition $\bfdd$, then we have 
	\begin{equation}\label{Apartition}
	\Spr(\cO_{\bfdd}) = \pi_{\bfdd} := j_{S_{\bfdd^*}}^{S_{n}} \sgn,
	\end{equation}
	where 
    \[
    S_{\bfdd^*}:=S_{d^*_1}\times S_{d^*_2} \times\dots \times S_{d^*_k}
    \]
    is the parabolic subgroup in $S_n$ attached to the sub root system of type $A_{d_1^*-1}\times \dots \times A_{d_k^*-1}$.

    Consequently, we have the following. 
    \begin{Lem}\label{j.2}
Let $\bf d$ and $\bf f$ be two partitions. 
Then we have
\[
   j^{S_{\abs{\bfdd}+\abs{\bfff}}}_{S_{\abs{\bfdd}}\times S_{\abs{\bfff}}}({\pi}_{\bfdd}\otimes {\pi}_{\bfff}) = \pi_{{\bfdd}\cupcol{\bfff}}.
\]
\end{Lem}
	


\subsection{Representations of $W_n$} 

    In the following, we identify the Weyl group $W_n$ of the root system of type $B_n$ with $S_n \ltimes \set{\pm 1}^n$, where $S_n$ is the parabolic subgroup of $W_n$ attached to the parabolic sub root system of type $A_{n-1}$ in $B_n$. Let $\varepsilon_n$ be the unique quadratic character of $W_n$ that is trivial on $S_n$. We identify the Coxeter group $W'_n$ of type $D_n$ with the kernel of $\varepsilon_n$. In that way, we made a choice of parabolic sub root system of type $A_{n-1}$ for $W'_n$.
    For a partition $\bfdd = [d_1,\dots, d_k]$ of $n$, the product group 
    \[
    W_{\bfdd} := W_{d_1}\times  \dots\times W_{d_k}
    \]
    is naturally identified with the subgroup of $W_n$ such that $W_{\bfdd} \cap S_n = S_{\bfdd}$.

    We now recall the well-known parameterization of the set of irreducible representations of $W_n$ by bipartitions of $n$, see \cite[\S 10.1]{CM}.
    By abuse of notation, we identify a $S_n$-module with its pull back to $W_n$ via the natural quotient map $W_n\to S_n$.  
    For a pair $(\bfdd,\bfff)$ of partitions
    such that $\abs{\bfdd}+\abs{\bfff} = n$ (i.e. a bipartition of $n$), 
    define
    \[
    \pi_{(\bfdd,\bfff)} := \Ind_{W_{\abs{\bfdd}}\times W_{\abs{\bfff}}}^{W_n} \pi_{\bfdd}\otimes (\pi_{\bfff}\otimes \varepsilon_{\abs{\bfff}}). 
    \]

    Suppose $\bfdd = [d_1,d_2,\dots, d_k]$ and $\bfff = [f_1,f_2,\dots, f_l]$. Let  
     \[
     W_{{({\bfdd},{\bfff})}} := W'_{d_1} \times \dots \times W'_{d_k} 
     \times W_{f_1}\times \dots W_{f_l} 
     \]
     be the subgroup of $W_n$  generated by  the subsystem  of type $D_{d_1}\times D_{d_2}\times \dots \times D_{d_k}\times B_{f_1}\times B_{f_2}\times \dots \times B_{f_l}$. 
    By Lusztig \cite[\S 4.5]{Lus84}, we have 
    \[
    \pi_{(\bfdd,\bfff)}  = j_{W_{{(\bfdd^*,\bfff^*)}}}^{W_{n}} \sgn
    \]
    where $\sgn$ denote the sign representation of the Weyl group $W_{(\bfdd^*,\bfff^*)}$. 
    In particular, all representations of $W_n$ are univalent and 
    \begin{equation}\label{eq:fdegree}
     b(\pi_{(\bfdd,\bfff)}) = \sum_{c \in \bfdd^*\cup \bfff^*} c^2 - \abs{\bfdd}. 
    \end{equation}
    By induction by stage, we have the following lemma.
    
\begin{Lem}\label{j.3}
The $j$-induction has the following properties:

\begin{enumerate}\label{induction-b}
\item  
   Let $({\bf d}_1,{\bf f}_1)$ and $({\bf d}_2,{\bf f}_2)$ be two bipartitions. 
   Then  we have	 \[
   j^{W_{\abs{\bfdd_1}+\abs{\bfdd_2}+\abs{\bfff_1}+ \abs{\bfff_2}}}_{{W_{\abs{\bfdd_1}+\abs{\bfff_1}}}\times W_{\abs{\bfff_1}+\abs{\bfff_2}}}({\pi}_{({\bfdd}_1,{\bfff}_1)}\otimes {\pi}_{({\bfdd}_2,{\bfff}_2)})=\pi_{({\bfdd}_1\cupcol{\bfdd}_2, {\bfff}_1\cupcol{\bfff}_2 )}.
   \]
    \item 
 We have   
 \[j^{W_n}_{S_n} \sgn=
 \begin{cases}
	\pi_{([k]^*,[k]^*)} &\text{ if }n=2k,\\
	\pi_{([k+1]^*,[k]^*)} &\text{ if }n=2k+1.
	\end{cases}
\]	
\item 
	Let $\bfdd$ be a partition with $\bfdd^*=[d^*_1,d^*_2,\dots,d^*_k]$. 
    Then  
    \begin{equation}\label{A-B}
	j^{W_n}_{S_n}\pi_\bfdd=
	\pi_{({\bfpp},{\bfqq})},
	\end{equation}
	where ${\bf p}^*=[\floor{\frac{d^*_{1}+1}{2}},\dots,\floor{\frac{d^*_{k}+1}{2}}]$ and ${\bf q}^*=[\floor{\frac{d^*_{1}}{2}},\dots,\floor{\frac{d^*_{k}}{2}}]$.
\end{enumerate}
\end{Lem}
\begin{proof}
Property $(2)$ follows from \cite[ (4.5.4)]{Lus84} and \cite[(5.1)]{Lu79}.
 $(1)$ and $(3)$ follow from induction by stages. 
\end{proof}

\trivial[h]{
We know $
	\pi_{\bf d} = j^{S_{n}}_{\prod\limits_{i}{S_{d^*_{i}}}} \sgn.
	$
}




\subsection{Representations of $W'_n$}For type $D_n$, we take the simple root system $\Delta=\{\alpha_1=e_1-e_2,\dots,\alpha_{n-1}=e_{n-1}-e_n, \alpha_n=e_{n-1}+e_n\}$.
 We use $S_n\subset W'_n$ to denote the subgroup generated by reflections of the first $n-1$ simple roots and $S'_n\subset W'_n$ to denote the subgroup generated by reflections of the first $n-2$ simple roots	and $\alpha_n$. 
The irreducible representations of $W'_n$ is obtained by the restriction of that of $W_n$ and therefore is parameterized by an unordered pair of partitions with an addition label:
Let $\set{\bfdd,\bfff}$ be an unordered pair of partitions such that $\abs{\bfdd}+\abs{\bfff}=n$ (i.e. an unordered bipartition of $n$). Then 
$\pi_{(\bfdd,\bfff)}|_{{W'_n}} = \pi_{(\bfff,\bfdd)}|_{W'_n}$. 
When $\bfdd\neq \bfff$, $\pi_{\set{\bfdd,\bfff}}:=\pi_{(\bfdd,\bfff)}|_{W'_n}$ is irreducible. 
When $\bfdd =  \bfff$, 
\[
\pi_{(\bfdd,\bfdd)}|_{W'_n} = \pi_{\set{\bfdd,\bfdd}}^I \oplus \pi_{\set{\bfdd,\bfdd}}^{II},
\]
where $\pi_{\set{\bfdd,\bfdd}}^{I}$ and $\pi_{\set{\bfdd,\bfdd}}^{II}$ are non-isomorphic irreducible $W'_n$-representations and we adopt the convention such that 
\begin{equation}\label{eq:Dlabel}
\pi_{\set{\bfdd,\bfdd}}^{I} = j_{S_n}^{W'_n} \pi_{\bfdd\cuprow \bfdd}, 
\text{~and~}
 \pi_{\set{\bfdd,\bfdd}}^{II} = j_{S'_n}^{W'_n} \pi_{\bfdd\cuprow \bfdd}.
\end{equation}

Now we would like to describe the set of univalent representations of $W'_n$. 
\begin{Lem}\label{lem:univ}
    The set of univalent representations of $W'_n$ equals to 
    \[
    \Set{\pi_{\set{\bfdd,\bfdd}}^{I}, \pi_{\set{\bfdd,\bfdd}}^{II} | {2\abs{\bfdd}=n}}\cup
    \Set{\pi_{\set{\bfdd,\bfff}}| {\abs{\bfdd}+\abs{\bfff} = n\text{ and } \abs{\bfdd}\neq \abs{\bfff}}}.
    \]
\end{Lem}
\begin{proof}
   The representation $\pi_{\set{\bfdd,\bfdd}}^{I}$  is univalent since it is obtained by $j$-induction, and $\pi_{\set{\bfdd,\bfdd}}^{II}$ is univalent since it is the conjugation of  $\pi_{\set{\bfdd,\bfdd}}^{I}$  by an element in $W_n\setminus W'_n$.  
   
    Now we suppose that $\bfdd\neq \bfff$. Note that $\pi_{\set{\bfdd,\bfff}}$ occurs (always with multiplicity one) in  $\pi_{(\bfdd',\bfff')}|_{W'_n}$ if and only if $(\bfdd,\bfff) = (\bfdd',\bfff')$ or $(\bfdd,\bfff) = (\bfff',\bfdd')$. 
   Recall that all representations of $W_n$ are univalent. 
    By \eqref{eq:fdegree}, $\abs{\bfdd}\neq \abs{\bfff}$ if and only if the fake degrees $b(\pi_{(\bfdd,\bfff)})\neq b(\pi_{(\bfff,\bfdd)})$. 
    Let $b = \min(b(\pi_{(\bfdd,\bfff)}), b(\pi_{(\bfff,\bfdd)}))$.
   Then  
    $\pi_{\set{\bfdd,\bfff}}$ occurs in $P_{b}(V)$ with multiplicity one if  $b(\pi_{(\bfdd,\bfff)})\neq b(\pi_{(\bfff,\bfdd)})$ and multiplicity two otherwise. 
    
    This finished the proof. 
\end{proof}
Now it is clear that 
\[
 b(\pi_{\set{\bfdd,\bfdd}}^{I}) = b(\pi_{\set{\bfdd,\bfdd}}^{II})=b(\pi_{(\bfdd,\bfdd)})
\quad\text{ and }\quad b(\pi_{\set{\bfdd,\bfff}}) = b(\pi_{(\bfdd,\bfff)})\text{ when }\abs{\bfdd}\geq \abs{\bfff}.
\]


Up to induction by stage, the following proposition covers all the possible cases of the $j$-induction involving  Weyl groups of type $D$. It will be useful later. 
\begin{Prop}\label{d-induciton} 
Let $m$ be a positive integer, and $X, Y \in \set{I,II}$.  
Let $\bfdd,\bfff,\bfdd'$ and $\bfff'$ be partitions such that 
     $\abs{\bfdd}>\abs{\bfff}$ and $\abs{\bfdd'}>\abs{\bfff'}$. 
 Then the $j$-induction has the following properties:
       \begin{enumerate}
           \item  $j_{S_m}^{W'_m}  \sgn  =\begin{cases}
	\pi^I_{\set{[k]^*,[k]^*}} &\text{ if }m=2k,\\
	\pi_{\set{[k+1]^*,[k]^*}} &\text{ if }m=2k+1.
	\end{cases}$
 
\item 
 $j_{S'_m}^{W'_m}  \sgn  =\begin{cases}
	\pi^{II}_{\set{[k]^*,[k]^*}} &\text{ if }m=2k,\\
	\pi_{\set{[k+1]^*,[k]^*}} &\text{ if }m=2k+1.
	\end{cases}$
  \item  $\pi_{(\bfdd,\bfff)} =  j_{W'_{\abs{\bfdd}+\abs{\bfff}}}^{W_{\abs{\bfdd}+\abs{\bfff}}}
         \pi_{\set{\bfdd,\bfff}}$.
      
     \item 
	Let $\bfdd$ be a partition of $n$ with $\bfdd^*=[d^*_1,d^*_2,\dots,d^*_k]$. 
    Then  
    \begin{equation}\label{A-D}
	j^{W'_n}_{S_n}\pi_\bfdd=
	\pi_{\set{\bfdp,\bfqq}},
	\end{equation}
	where  ${\bf p}^*=[\floor{\frac{d^*_{1}}{2}},\dots,\floor{\frac{d^*_{k}}{2}}]$ and ${\bf q}^*=[\floor{\frac{d^*_{1}+1}{2}},\dots,\floor{\frac{d^*_{k}+1}{2}}]$.

        \item  $\pi_{\set{\bfdd\cupcol \bfdd',\bfff\cupcol \bfff'}}  =  j_{W'_{\abs{\bfdd}+\abs{\bfff}}\times W'_{\abs{\bfdd'}+\abs{\bfff'}}}^{W'_{\abs{\bfdd}+\abs{\bfff}+\abs{\bfdd'}+\abs{\bfff'}}}
        ( \pi_{\set{\bfdd,\bfff}}\otimes 
         \pi_{\set{\bfdd',\bfff'}})$.
         \item $\pi_{\set{\bfdd\cupcol \bfdd',\bfff\cupcol \bfdd'}}  =
        j_{W'_{\abs{\bfdd}+\abs{\bfff}}\times W'_{2\abs{\bfdd'}}}^{W'_{\abs{\bfdd}+\abs{\bfff}+2\abs{\bfdd'}}}
         (\pi_{\set{\bfdd,\bfff}}\otimes 
         \pi_{\set{\bfdd',\bfdd'}}^X)$.
         \item $\pi_{\set{\bfdd\cupcol \bfdd',\bfdd\cupcol \bfdd'}}^Z  =
        j_{W'_{2\abs{\bfdd}}\times W'_{2\abs{\bfdd'}}}^{W'_{2\abs{\bfdd}+2\abs{\bfdd'}}}
         (\pi_{\set{\bfdd,\bfdd}}^X\otimes 
         \pi_{\set{\bfdd',\bfdd'}}^Y)$
         (Here $Z = I$ if $X=Y$ and $Z=II$ otherwise).
          \item $j_{S_m\times W'_{2\abs{\bfdd}}}^{W'_{m+2\abs{\bfdd}}} \sgn \otimes 
         \pi_{\set{\bfdd,\bfdd}}^X  = 
         \begin{cases}
	\pi^X_{\set{[k]^*\cupcol \bfdd,[k]^*\cupcol \bfdd}} &\text{ if }m=2k,\\
	\pi_{\set{[k+1]^*\cupcol \bfdd,[k]^*\cupcol \bfdd}} &\text{ if }m=2k+1.
	\end{cases}$
  \item $j_{S'_m\times W'_{2\abs{\bfdd}}}^{W'_{m+2\abs{\bfdd}}} \sgn \otimes 
         \pi_{\set{\bfdd,\bfdd}}^X  = 
         \begin{cases}
	\pi^Y_{\set{[k]^*\cupcol \bfdd,[k]^*\cupcol \bfdd}} &\text{ if }m=2k,\\
	\pi_{\set{[k+1]^*\cupcol \bfdd,[k]^*\cupcol \bfdd}} &\text{ if }m=2k+1.
	\end{cases}$
  (Here $X\neq Y$).
       \end{enumerate}
      
\end{Prop}
\begin{proof}
The bipartitions in the above formulas can be determined by the fake degree analysis as in  the proof of Lemma~\ref{lem:univ} and the $j$-induction by stage to $W_m$. 
The labels are matched due to \eqref{eq:Dlabel} and our fixed choices of the embedding of various subgroups.  

\trivial[h]{
To save the notation, let $d =\abs{d}$, etc. 
The most tricky formula is the second last one. 
When $X=Y=I$, we have $j_{W'_{2d}\times W'_{2d'}}^{W'_{2d+2d'}} \pi_{\bfdd,\bfdd}^I\otimes \pi_{\bfdd',\bfdd'}^I
= j_{S_{2d}\times S_{2d'}}^{W'_{2d+2d'}} \pi_{\bfdd\cuprow \bfdd} \otimes \pi_{\bfdd'\cuprow \bfdd'} = 
j_{S_{2d+2d'}}^{W'_{2d+2d'}} \pi_{(\bfdd\cupcol \bfdd')\cuprow (\bfdd\cupcol \bfdd')} =  \pi_{\bfdd\cupcol \bfdd'}^I$.

When $X=I,Y=II$, we have $\pi^{II}_{\bfdd',\bfdd'}$ is the conjugation of $\pi^I_{\bfdd',\bfdd'}$ by an element in $W_{2d'}\setminus W'_{2d'}$. 
So the $j$-induction is also the conjugation of the results when $(X,Y)=(I,I)$, which implies $Z=II$.  
Similarly, the situation happens for $X=II, Y=I$. 

When $(X,Y)=(II,II)$, one conjugate twice, one by an element in $W_{2d}\setminus W'_{2d}$ another by an element in $W_{2d'}\setminus W'_{2d'}$. 
Since the product of these two elements is in $W'_{2d+2d'}$, we conclude that $Z=I$. 
}
\end{proof}

\trivial[h]{	
	\begin{Prop}[\cite{Lus84, Ca85}]\label{j}
		The $j$-induction has the following properties:
		\begin{itemize}
			\item[(1)] If given two parabolic subgroups (generated by reflections) $W''\subset W'\subset W$, then $j^W_{W''}=j^W_{W'}\circ j^{W'}_{W''}$;
			\item[(2)] If $\bf d$ and $\bf f$ are two partitions corresponding to two
			representations of $S_m$ and $S_n$ with $m=|{\bf d}|$ and $n=|{\bf f}|$, then
			$$j^{S_{m+n}}_{S_m\times S_n}({\pi}_{\bf d}\otimes {\pi}_{\bf f})=\pi_{{\bf d}{\stackrel{c}{\sqcup}}{\bf f}},$$
			where
			${\bf d}{\stackrel{c}{\sqcup}}{\bf f}$ is the partition such that
			$$\{c_j({\bf d}{\stackrel{c}{\sqcup}}{\bf f})|j\in\mathbb{N}^+ \}=\{d^*_i|i\in\mathbb{N}^+ \}\cup \{f^*_k|k\in\mathbb{N}^+ \}$$
			as multisets;
			\item[(3)]If $({\bf d}_1,{\bf f}_1)$ and $({\bf d}_2,{\bf f}_2)$ are two bipartitions corresponding to irreducible representations of $W_n$ (or $W'_n$) and $W_m$ (or $W'_m$),  (here $n=|{\bf d}_1|+|{\bf f}_1|$ and $m=|{\bf d}_2|+|{\bf f}_2|$),  then  we have	 $$j^{W_{n+m}}_{{W_{n}}\times W_m}({\pi}_{({\bf d}_1,{\bf f}_1)}\otimes {\pi}_{({\bf d}_2,{\bf f}_2)})=\pi_{({\bf d}_1{\stackrel{c}{\sqcup}}{\bf d}_2, {\bf f}_1{\stackrel{c}{\sqcup}}{\bf f}_2 )}=j^{W_{n+m}}_{W'_{({\bf d}_1{\stackrel{c}{\sqcup}}{\bf d}_2, {\bf f}_1{\stackrel{c}{\sqcup}}{\bf f}_2 )}}\sgn.$$
			
   When ${\bf d}_i\ne {\bf f}_i$ and ${\bf d}_1{\stackrel{c}{\sqcup}}{\bf d}_2= {\bf f}_1{\stackrel{c}{\sqcup}}{\bf f}_2$, we have	$$j^{W'_{n+m}}_{{W'_{n}}\times W'_m}({\pi}_{({\bf d}_1,{\bf f}_1)}\otimes {\pi}_{({\bf d}_2,{\bf f}_2)})=\pi^{?}_{({\bf d}_1{\stackrel{c}{\sqcup}}{\bf d}_2, {\bf f}_1{\stackrel{c}{\sqcup}}{\bf f}_2 )}=?$$

   When ${\bf d}_i\ne {\bf f}_i$ and ${\bf d}_1{\stackrel{c}{\sqcup}}{\bf d}_2\ne {\bf f}_1{\stackrel{c}{\sqcup}}{\bf f}_2$, we have	 $$j^{W'_{n+m}}_{{W'_{n}}\times W'_m}({\pi}_{({\bf d}_1,{\bf f}_1)}\otimes {\pi}_{({\bf d}_2,{\bf f}_2)})=\pi_{({\bf d}_1{\stackrel{c}{\sqcup}}{\bf d}_2, {\bf f}_1{\stackrel{c}{\sqcup}}{\bf f}_2 )}=j^{W'_{n+m}}_{W'_{({\bf d}_1{\stackrel{c}{\sqcup}}{\bf d}_2, {\bf f}_1{\stackrel{c}{\sqcup}}{\bf f}_2 )}}\sgn;$$
	\item[(4)]Suppose $({\bf d}_1,{\bf f}_1)$ and $({\bf d}_2,{\bf f}_2)$ are two bipartitions corresponding to irreducible representations of  $W'_n$ and $W'_m$ (here $n=|{\bf d}_1|+|{\bf f}_1|$ and $m=|{\bf d}_2|+|{\bf f}_2|$). If ${\bf d}_i= {\bf f}_i$, we will  have	$$j^{W'_{n+m}}_{{W'_{n}}\times W'_m}({\pi}^X_{({\bf d}_1,{\bf f}_1)}\otimes {\pi}^Y_{({\bf d}_2,{\bf f}_2)})=\pi^Z_{({\bf d}_1{\stackrel{c}{\sqcup}}{\bf d}_2, {\bf f}_1{\stackrel{c}{\sqcup}}{\bf f}_2 )}=j^{W'_{n+m}}_{W'_{({\bf d}_1{\stackrel{c}{\sqcup}}{\bf d}_2, {\bf f}_1{\stackrel{c}{\sqcup}}{\bf f}_2 )}}\sgn,$$
    where $X, Y, Z\in \{I,II\}$ denote the numerals of the representations and $Z=I$ if $X=Y$, otherwise it will be $II$.
		\end{itemize}
		
	\end{Prop}
}	

	

%
%

\subsection{Springer correspondence}

        Nilpotent orbits in $\fgg$ are parameterized by partitions (and with additional labels in type $D$). 
	Let $\bfdd$ be a partition of $2n$
	in types $C,D$ and $2n+1$ in type $B$. 
        We recall Sommer's formulation of Shoji's results on the explicit Springer correspondence in \cite{Sh79}. 
 
        First separate the parts of $\bfdd^*$
	into its odd parts $2\alpha_1 + 1 \geq 2 \alpha_2 + 1 \geq
	2\alpha_3 + 1 \dots
	\geq 2 \alpha_r + 1$
	and its even parts $2 \beta_1 \geq 2 \beta_2 \geq 2 \beta_3 \dots
	\geq 2 \beta_s$ where $\alpha_i$ and $\beta_j$ are positive integers.
        Define 
        \[
        \bfpp^* = \begin{cases}
        \set{\alpha_{2k}+1|1\leq 2k\leq r} \cup \set{\beta_{2k}|1\leq 2k\leq s},
        & \text{ in type $B$},\\
        \set{\alpha_{2k+1}+1|1\leq 2k+1\leq r} \cup \set{\beta_{2k}|1\leq 2k\leq s},
        & \text{ in type $C$},\\
        \set{\alpha_{2k+1}+1|1\leq 2k+1\leq r} \cup \set{\beta_{2k+1}|1\leq 2k+1\leq s},
        & \text{ in type $D$},\\
        \end{cases}
        \]
        and 
        \[
        \bfqq^* = \begin{cases}
        \set{\alpha_{2k+1}|1\leq 2k+1\leq r} \cup \set{\beta_{2k+1}|1\leq 2k+1\leq s},
        & \text{ in type $B$},\\
        \set{\alpha_{2k}|1\leq 2k\leq r} \cup \set{\beta_{2k+1}|1\leq 2k+1\leq s},
        & \text{ in type $C$},\\
        \set{\alpha_{2k}|1\leq 2k\leq r} \cup \set{\beta_{2k}|1\leq 2k\leq s},
        & \text{ in type $D$}.\\
        \end{cases}
        \]


\begin{Prop}[{\cite[Lemma 8]{Sommers}}] \label{som}
The Springer correspondence can be given as follows:
 \begin{enumerate}
 \item 
  In types $B$ and $C$, let $\cO$ be the nilpotent orbit corresponding to $\bfdd$. 
  Then $\pi_{(\bfpp,\bfqq)}$ is the representation associated to the trivial local system on $\cO$ under the Springer correspondence. 
 \item If $\fgg$ is of type $D$ and $\bfdd$ is not very even, let $\cO$ be the nilpotent orbit corresponding to $\bfdd$.
  Then $\pi_{\set{\bfpp,\bfqq}}$ is the representation associated to the trivial local system on $\cO$ under the Springer correspondence. 
 \item If $\fgg$ is of type $D$ and $\bfdd$ is very even, 
    there are two nilpotent orbits $\cO^I$ and $\cO^{II}$ having the Jordan block $\bfdd$ where $\cO^I$ is the orbit which could be obtained by orbit induction from the Levi subgroup attached to $S_{\abs{\bfdd}/2}$. 
  Then $\pi_{\set{\bfqq,\bfqq}}^I$ and $\pi_{\set{\bfqq,\bfqq}}^{II}$ are the representations associated to the trivial local system on $\cO^I$ and $\cO^{II}$ respectively. 
 \end{enumerate}
\end{Prop}

If  we identify partitions, Young diagrams with nilpotent orbits, and
		bipartitions with $W_{n}$-representations, we may use $E_{\bfdd}$ or $E_{\mathcal{O}}$ to denote the Springer representation in the above Proposition \ref{som}.


	Given two partitions ${\bf d}=[d_1,\dots,d_N]$ and ${\bf f}=[f_1,\dots,f_k]$ of $N$, we say that 	${\bf d}$ {\it dominates } 	${\bf f}$ if the following condition holds:
	\begin{equation}\label{domi}
	    \sum\limits_{1\leq j\leq k}d_j\geq \sum\limits_{1\leq j\leq k}f_j
	\end{equation} for $1\leq k\leq N$.
	\begin{Defff}[Collapse]	
		Let ${\bf d}=[d_1,\dots,d_k]$ be a partition of $2n$. There is a unique largest partition of $2n$ of type $D_n$ dominated by ${\bf d}$. If  ${\bf d}$ is not a partition of type $D_n$, then one of its even parts must occur with odd multiplicity. Let $q$ be the largest such part. Then replace the last occurrence of $q$ in ${\bf d}$ by $q-1$ and the first subsequent part $r$ strictly less than $q-1$ by $r+1$. Repeat this process until a partition of type $D_n$ is obtained. This new partition of type $D_n$
		is called the {\it $D$-collapse} of 	${\bf d}$, and we denote it by ${\bf d}_D$. Similarly there are $B$-collapse and $C$-collapse of ${\bf d}$.
	\end{Defff}	
		
		

	
	Some more properties for the  collapse  of partitions can be found in \cite{CM}.

	\begin{Prop}[{\cite[{Proposition 8.3}]{BMSZ-counting}}]\label{D-C}	
		Let $E_{\bf d}\in \mathrm{Irr}(W'_n)$ be a special representation corresponding to a special nilpotent orbit with partition ${\bf d}$, where   ${\bf d}=[d_1,\dots,d_t]$ is a partition of type $D_n$. The special representation $E_{\bf d}\otimes \sgn$ will correspond to the partition  $({\bf d}^*)_D $. And $j^{W_n}_{W'_n}E_{\bf d}$ will correspond to the partition  $(({\bf d}^*)_D)^* $ under the Springer correspondence of type $C$.
	\end{Prop}
	In \cite{BMSZ-typeC}, the partition  $(({\bf d}^*)_D)^* $ is called  {\it metaplectic special} since its dual partition is of type $D$.

\section{Symbols and nilpotent orbits}\label{Weyl}
	In this section, we recall the procedure of obtaining the nilpotent orbit $\mathcal{O}_w$ from a given Weyl group element $w$. Some details can be found in \cite{BarV82,Lus84,Ca85,CM}.
	
	For $i, j\in\mathbb{Z}$ with $i<j$, set $[i,j]= \{i,i+1,\dots, j-1, j\} $. Recall that we map $w\in W$ for classical types to the sequence $ w=(w(1),\dots,w(n)) $. Denote by $s_i$ the involution $ w\in W$ such that $ w(i)=i+1 $, $ w(i+1)=i$ and $w(k)=k$ for all $k\in[1, i-1]\cup[i+2, n]$.
	
	For  a totally ordered set $ \Gamma $, we  denote by $ \mathrm{Seq}_n (\Gamma)$ the set of sequences $ x=(x_1,x_2,\dots, x_n) $   of length $ n $ with $ x_i\in\Gamma $.
	
	For $ x=(x_1,x_2,\dots,x_n)\in \mathrm{Seq}_n (\Gamma)$, set
	\begin{equation*}
	\begin{aligned}
	{x}^-=&(x_1,x_2,\dots,x_{n-1}, x_n,-x_n,-x_{n-1},\dots,-x_2,-x_1),\\
	{}^-{x}=&(-x_n,-x_{n-1},\dots, -x_2,-x_1,x_1,x_2,\dots, x_{n-1}, x_n).
	\end{aligned}
	\end{equation*}

	\subsection{$B$-symbol}\label{B}

 When the root system $\Phi=B_n$ or $C_n$ with $n>1$, the Weyl group is $W=W_n$, where $ W_n $ is the group consisting of permutations $ w $ of the set $[-n, n]$
	such that $ w(-i)=-w(i) $ for all $ i \in[1,n]$. Let $ t\in W_n$ be the element with $ t=(-1,2, \dots, n):=(1,-1)$. Then $ (W_n , T_n)$ is a Coxeter system, with $T_n=\{t,s_1,\dots s_{n-1}\} $.

	Now we recall the Lusztig's symbols \cite{lusztig1977symbol}. Let
	\[
	\begin{pmatrix}
	\lambda_1~\lambda_2~\dots~\lambda_{m+1}\\\mu_1~\mu_2~\dots~\mu_m
	\end{pmatrix}, m\geq 0
	\]
	be a tableau of nonnegative integers such that entries in each row are  strictly increasing. Define an equivalence relation on the set of all such tableaux via
	\begin{equation*}\label{eq:equiv}
	\begin{pmatrix}
	\lambda_1~\lambda_2~\dots~\lambda_{m+1}\\\mu_1~\mu_2~\dots\mu_m
	\end{pmatrix} \sim
	\begin{pmatrix}
	0~\lambda_1+1~\lambda_2+1~\dots~\lambda_{m+1}+1\\0~\mu_1+1~\mu_2+1~\dots\mu_m+1
	\end{pmatrix} .
	\end{equation*}
	Denote by $ \Sigma_B $ the set of equivalence classes under this relation $ \sim $. Use the same notation $ \Lambda=\begin{pmatrix}
	\lambda_1~\lambda_2~\dots~\lambda_{m+1}\\\mu_1~\mu_2~\dots\mu_m
	\end{pmatrix}  \in\Sigma_B$ to denote its equivalence class, called a \textit{$ B $-symbol}.
	
	For $ w\in W_n $, from \cite{BXX}, we have a Young tableau  $ P({}^-{ w}) $, which is  obtained by applying the R-S algorithm to the sequence\[
	{}^-{w}=(-w(n),-w(n-1),\dots, -w(1), w(1), \dots, w(n-1), w(n)),
	\]
	and $p({}^-{w})=\sh(P({}^-{w}))$, which is a partition of $2n$. By \cite[Prop.17]{BarV82}, there is a symbol
	$\Lambda =\begin{pmatrix}
	\lambda_1~\lambda_2~\dots~\lambda_{m+1}\\\mu_1~\mu_2~\dots\mu_m
	\end{pmatrix} \in\Sigma_B $ such that, as multisets
	\begin{equation*}
	\{2\lambda_i,2 \mu_j +1\mid  i\leq  m+1, j\leq m  \} =\{ p_k+2m+1-k \mid  k\leq  2m+1\},
	\end{equation*}
	where $p=p({}^-{w})=[p_1,p_2,\dots]$ with $p_{2m+2}=0$. Here we set $p_l=0$ when $l>2n$. It gives a well-defined map\begin{equation}\label{eq:mapsymb}
	Symb_B:W_n\to \Sigma_B.
	\end{equation}

	\subsection{$D$-symbol}\label{d-weyl}
	When the root system $\Phi=D_n$ with $n>3$, the Weyl group is $W=W_n'$, which consists of those elements $ w$ of  $W_ n$  such that  the number of negative integers in $ \{w(1), w(2), \dots, w(n) \}   $ is even.  Let $ u=ts_1t=(1,-2) $ and $ T'_n=\{u,s_1,s_2,\dots,$  $s_{n-1}\} $. Thus $ (W_n', T'_n) $ is a Coxeter system.
	
	
	
	Here we recall Lusztig's $D$-symbols \cite{lusztig1977symbol}. Let
	\[ \begin{pmatrix}
	\lambda_1~\lambda_2~\dots~\lambda_{m}\\\mu_1~\mu_2~\dots~\mu_m
	\end{pmatrix}, m\geq 0\]
	be a tableau of nonnegative integers such that entries in each row are  strictly increasing. Define an equivalence relation on the set of all such tableaux via
	\begin{equation*}\label{eq:equiv-d}
	\begin{pmatrix}
	\mu_1~\mu_2~\dots\mu_m\\
	\lambda_1~\lambda_2~\dots~\lambda_{m}
	\end{pmatrix} \sim
	\begin{pmatrix}
	\lambda_1~\lambda_2~\dots~\lambda_{m}\\\mu_1~\mu_2~\dots\mu_m
	\end{pmatrix}
	\sim
	\begin{pmatrix}
	0~\lambda_1+1~\lambda_2+1~\dots~\lambda_{m}+1\\0~\mu_1+1~\mu_2+1~\dots\mu_m+1
	\end{pmatrix} .
	\end{equation*}
	Denote by $ \Sigma_D $ the set of equivalence classes under this relation $ \sim $. Use the same notation $ \Lambda=\begin{pmatrix}
	\lambda_1~\lambda_2~\dots~\lambda_{m}\\\mu_1~\mu_2~\dots\mu_m
	\end{pmatrix}  \in\Sigma_D$ to denote its equivalence class,  called a \textit{$ D$-symbol}. Define a map
	\begin{equation}\label{eq:mapbd}
	d:\Sigma_B\to \Sigma_D ,
	\end{equation}
	\[
	\begin{pmatrix}
	\lambda_1~\lambda_2~\dots~\lambda_{m+1}\\\mu_1~\mu_2~\dots\mu_m
	\end{pmatrix}  \mapsto \begin{pmatrix}
	\lambda_1&\lambda_2&\dots&\lambda_{m+1}\\0 &\mu_1+1&\dots&\mu_m+1
	\end{pmatrix}.
	\]

	Let \begin{equation}\label{eq:mapsymb-d}
	Symb_D: W'_n\to \Sigma_D
	\end{equation}
	be  the restriction to $ W_n' $ of the composition $ d\circ Symb_B $ of maps from  \eqref{eq:mapsymb} and \eqref{eq:mapbd}.

	We have the following proposition.
	
	\begin{Prop}[{\cite[\S 4.5 ]{Lus84}}]
		The following map
		\begin{align*}
		\phi_B: \Sigma_B \rightarrow \mathrm{Irr}(W_n)
		\end{align*}
		which maps a $ B $-symbol to an irreducible representation of $W_n$ is a one-to-one correspondence.
		
		The following map
		\begin{align*}
		\phi_D: \Sigma_D \rightarrow \mathrm{Irr}(W'_n)
		\end{align*}
		which maps a $ D $-symbol to an irreducible representation of $W'_n$ is a one-to-one correspondence  (except that when the $D$-symbol corresponds to a very even partition, there are two irreducible representations of  $W'_n$ corresponding to it).
	\end{Prop}
	

	Recall that  the elements of $\mathrm{Irr}(W_n)$ are parameterized by  the bipartitions $(\bf{ d},\bf{f})$   with $|{\bf d}|+|{\bf f}|=n$. So there is a one-to-one correspondence between $B$-symbols and  bipartitions. It can be described as follows. Let $\Lambda =\begin{pmatrix}
	\lambda_1~\lambda_2~\dots~\lambda_{m+1}\\\mu_1~\mu_2~\dots\mu_m
	\end{pmatrix} \in\Sigma_B $. We subtract $i-1$ from the $i$-th element  of the top row and  bottom row.
	Then we get a pair $(\bf{ d},\bf{f})$ corresponding to the new top row and bottom row.


Recall that the elements of $\mathrm{Irr}(W'_n)$ are parameterized by the unordered bipartitions $\{ \bf{ d},\bf{f}\}$  with $|{\bf d}|+|{\bf f}|=n$, except that if $n=2m'$ is even and $\bf{ d}=\bf{f}$, then the unordered bipartition  $\{ \bf{ d},\bf{d}\}$ corresponds to two representations,  $\pi^I_{\{\bf{d},\bf{d}\}}$ and $\pi^{II}_{\{\bf{d},\bf{d}\}}$.	So there is a one-to-one correspondence between D-symbols and  unordered bipartitions 	except that if $n=2m'$ is even and $\bf{ d}=\bf{f}$. It can be described as follows. Let $\Lambda =\begin{pmatrix}
	\lambda_1~\lambda_2~\dots~\lambda_{m}\\\mu_1~\mu_2~\dots\mu_m
	\end{pmatrix} \in\Sigma_D $. We subtract $i-1$ from the $i$-th element  of the top row and  bottom row.
	Then we get an unordered bipartition  $\{ \bf{ d},\bf{f}\}$ corresponding to the new top row and bottom row.

	\subsection{Special orbits}\label{special-symbol}
	
	In the sense of Lusztig \cite{Lus84}, $\mathrm{Irr}(W)$ is a disjoint union of families (i.e., double cells).
	
	\begin{Prop}[{\cite[Theorem 18 ]{BarV82}}]\label{dcell}
		Let $W$ be the Weyl group of type $B_n$ or $D_n$.
		Then two irreducible representations $\Lambda$ and $\Lambda'$ of $W$ belong to the same double cell if and only if the symbol $\Lambda$ is a permutation of $\Lambda'$.
		
	\end{Prop}
	
	A $B$-symbol $\Lambda =\begin{pmatrix}
	\lambda_1~\lambda_2~\dots~\lambda_{m+1}\\\mu_1~\mu_2~\dots\mu_m
	\end{pmatrix} $ is called \textit{special} if
	$$\lambda_1\leq \mu_1\leq \lambda_2\leq \mu_2\leq \dots\leq \mu_m\leq \lambda_{m+1}.$$

	A $D$-symbol $\Lambda =\begin{pmatrix}
	\lambda_1~\lambda_2~\dots~\lambda_{m}\\\mu_1~\mu_2~\dots\mu_m
	\end{pmatrix} $ is called \textit{special} if
	$$\lambda_1\leq \mu_1\leq \lambda_2\leq \mu_2\leq \dots\leq \lambda_m\leq \mu_m \text{~or~} \mu_1\leq \lambda_1 \leq  \mu_2\leq \lambda_2\leq \dots\leq \mu_m \leq \lambda_m.$$

	For a given symbol $\Lambda$, it is clear that there is a unique special symbol $\Lambda^s$, such that they are in the same double cell. An irreducible representation of $W$ is called \textit{special} if and only if the corresponding symbol is special.
	
	Suppose $\mathfrak{g}$ is of type $B_n$. Let ${\mathcal{O}}_{\bf{d}}$ be a nilpotent orbit of $\mathfrak{g}$ with partition ${\bf d}=[d_1,\dots,d_k]$ of $2n+1$, then $k$ is odd. From \cite[\S 10.1]{CM}, there is a $B$-symbol
	$\Lambda =\begin{pmatrix}
	\lambda_1~\lambda_2~\dots~\lambda_{m+1}\\\mu_1~\mu_2~\dots\mu_m
	\end{pmatrix} $ such that, as multisets
	\begin{equation*}
	\{2\lambda_i+1,2 \mu_j\mid  i\leq  m+1, j\leq m  \} =\{ d_l+k-l \mid  l\leq  k\},
	\end{equation*}
	where $2m+1=k$.

	Suppose $\mathfrak{g}$ is of type $C_n$. Let ${\mathcal{O}}_{\bf{d}}$ be a nilpotent orbit of $\mathfrak{g}$ with partition ${\bf d}=[d_1,\dots,d_k]$ of $2n$, then we may suppose $k$ is odd (otherwise we add $0$ as the last part of ${\bf d}$). From \cite[\S 10.1]{CM}, there is a $B$-symbol
	$\Lambda =\begin{pmatrix}
	\lambda_1~\lambda_2~\dots~\lambda_{m+1}\\\mu_1~\mu_2~\dots\mu_m
	\end{pmatrix} $ such that, as multisets
	\begin{equation*}
	\{2\lambda_i,2 \mu_j+1\mid  i\leq  m+1, j\leq m  \} =\{ d_l+k-l \mid  l\leq  k\},
	\end{equation*}
	where $2m+1=k$.

	Suppose $\mathfrak{g}$ is of type $D_n$. Let ${\mathcal{O}}_{\bf{d}}$ be a nilpotent orbit of $\mathfrak{g}$ with partition ${\bf d}=[d_1,\dots,d_k]$ of $2n$, then $k$ is even. From \cite[\S 10.1]{CM}, 
 there is a $D$-symbol
	$\Lambda =\begin{pmatrix}
	\lambda_1~\lambda_2~\dots~\lambda_{m}\\\mu_1~\mu_2~\dots\mu_m
	\end{pmatrix} $ such that, as multisets
	\begin{equation*}\label{eq:symb}
	\{2\lambda_i+1,2 \mu_j\mid  i\leq  m, j\leq m  \} =\{ d_l+k-l \mid  l\leq  k\},
	\end{equation*}
	where $2m=k$.

	From \cite[\S 10.1 ]{CM}, we know that a symbol is special if and only if the corresponding nilpotent orbit is special.
	Since  the special symbols and special partitions (nilpotent orbits) are in one-to-one correspondence, we can use $\mathcal{O}_{\Lambda}$ to denote the special nilpotent orbit
	$\mathcal{O}_{\bf d}$ where $\Lambda$ is the special symbol of the special partition ${\bf d}$.
	
	In general, we have the following procedure to find the special nilpotent orbit for a given $w\in W$. From $w$, we can get a
	Young tableau $ P({}^-{ w}) $, then we can get a symbol $\Lambda$. By some permutation, we can get a special symbol $\Lambda^s$, then we can get the corresponding special nilpotent orbit $\mathcal{O}_w=\mathcal{O}_{\Lambda^s}$.
	
	\begin{ex}Let $ y=(-1,-2,3,4,5) \in W_5$. Then we have a Young tableau
\[P({}^-{y}) =\begin{tikzpicture}[scale=\domscale+0.1,baseline=-20pt]
		\hobox{0}{0}{-5}
		\hobox{1}{0}{-4}
		\hobox{2}{0}{-3}
		\hobox{3}{0}{-2}
  \hobox{4}{0}{3}
  \hobox{5}{0}{4}
  \hobox{6}{0}{5}
		\hobox{0}{1}{-1}
		\hobox{0}{2}{1}
             \hobox{0}{3}{2}
		\end{tikzpicture},\]
  
		a partition $ [7,1,1,1,0] $ and a  strictly increasing sequence \[
		(0,2,3,4,11).
		\] Thus we have a $B$-symbol $\Lambda= \begin{pmatrix}
		0~1~2\\1~5
		\end{pmatrix} $, which is not special. But from $\Lambda$, we have a unique special $B$-symbol $\Lambda^s= \begin{pmatrix}
		0~1~5\\1~2
		\end{pmatrix} $. Finally we can get ${\bf d}=[7,1,1,1,1]$ being a partition of special nilpotent orbit of type $B_5$, or ${\bf d}=[6,2,1,1]$ being a partition of special  nilpotent orbit of type $C_5$.
	\end{ex}

	\section{The  annihilator varieties for type $A_{n-1}$ }\label{type-A}

	Now we suppose  $\mathfrak{g}=\mathfrak{sl}(n, \mathbb{C})$.
	The Weyl group of $\mathfrak{g}$ is  $S_n$. 
  From Sagan \cite{Sagan} or Bai-Xie \cite[Lemma 4.1]{BaiX}, we know that there is a bijection between the Kazhdan-Lusztig right cells in the symmetric group $S_n$ and the Young tableaux through the famous Robinson-Schensted algorithm.
	 For any element $w\in S_n$, we use $P(w)$ to denote the corresponding Young tableau for any $w\in S_n$.	Suppose the shape of $P(w)$ is $p(w)=[p_1,p_2,\dots,p_N]$ which is a decreasing sequence and $\sum p_i=n$. From \cite{CM}, we know that  the nilpotent orbits of $\mathfrak{sl}(n, \mathbb{C})$ are in one-to-one  correspondence with the partitions of $n$. 
 Also we have $\mathcal{O}_w=\mathcal{O}_{p(w)}$.

	Let  $ \mathcal{O} $ be the BGG category consists of $ \mathfrak{g} $-modules which are  semisimple as $ \mathfrak{h} $-modules, finitely generated as $ U(\mathfrak{g}) $-modules and locally $ \mathfrak{n} $-finite. A \textit{block} of $ \mathcal{O} $ is an indecomposable summand of $ \mathcal{O} $ as an abelian subcategory. Let $ \mathcal{O} _\lambda $  be the block containing  the simple module $ L(\lambda) $.

	For two integral weights $ \lambda,\mu \in \mathfrak{h}^*$, we set $ \gamma=\mu-\lambda $. Then we can find a weight $ \bar{\gamma} \in W\gamma$ such that $ \langle {\bar{\gamma}},{\alpha} \rangle \geq 0$ for all $ \alpha \in\Delta^+$. Then $ L(\bar{\gamma}) $ is finite dimensional. The Jantzen's \textit{translation functor} (see Jantzen \cite{Ja79} or Humphreys \cite{Hum08}): $$  T_\lambda^\mu : \mathcal{O}_\lambda \rightarrow \mathcal{O}_\mu $$ is   an exact functor given by $ T_\lambda^\mu(M):=\text{ pr}_\mu (L(\bar{\gamma})\otimes M) $, where $ M\in   \mathcal{O}_\lambda $ and $ \text{ pr}_\mu $ is the natural projection $ \mathcal{O} \rightarrow \mathcal{O}_\mu $.

	From Borho-Brylinski-MacPherson \cite[Lemma 5.2]{BBM}, we know that the associated variety of an irreducible $\mathfrak{g}$-module is invariant under Jantzen's translation functor. In particular, we have the following proposition from Bai-Xie \cite[Corollary 3.3]{BaiX}.

	\begin{Prop}\label{trans}
		For any integral weight $\lambda$, we have
		\[
		T_{-w\rho-\rho}^\lambda (L_w)=L(\lambda)
		\]
		and $V(L(\lambda))=V(L_w)$, where $ w\in W $ is the unique element of minimal length such that $ w^{-1}\lambda $ is antidominant.
	\end{Prop}

	From this proposition we have $V(L(\lambda))=V(L_w)$ and $V(\Ann (L(\lambda)))=V(\Ann (L_w))=V(I_w)=\overline{\mathcal{O}}_w=\overline{\mathcal{O}}_{p(w)}$  by \cite{Jo84}.

	For an integral weight $\lambda$,  we use $p(\lambda)$ to denote the partition corresponding to the Young tableau $P(\lambda)$.   From  \cite[Lemma 4.5]{BaiX}, we know that $p(\lambda)$ and $p(w)$ are the same partitions, where $ w\in W $ is the unique element of minimal length such that $ w^{-1}\lambda $ is antidominant.
	
	So we have the following theorem.

	\begin{Thm}\label{spa}
		Let $\mathfrak{g}=\mathfrak{sl}(n, \mathbb{C})$.	Suppose $\lambda\in \mathfrak{h}^*$ is an integral weight. Then
		$$V(\Ann (L(\lambda)))=V(\Ann (L_w))=V(I_w)=\overline{\mathcal{O}}_w=\overline{\mathcal{O}}_{p(w)}=\overline{\mathcal{O}}_{p(\lambda)},$$ where $ w\in W $ is the unique element of minimal length such that $ w^{-1}\lambda $ is antidominant.

	\end{Thm}
	
	\begin{ex}Let $\mathfrak{g}=\mathfrak{sl}(6, \mathbb{C})$.
		Suppose $\lambda=(2, 4, 6, 9, -29, 8)$. Then
		\[
		P(\lambda)=
		\small{\begin{tikzpicture}[scale=\domscale+0.1,baseline=-20pt]
		\hobox{0}{0}{-29}
		\hobox{1}{0}{4}
		\hobox{2}{0}{6}
		\hobox{3}{0}{8}
		\hobox{0}{1}{2}
		\hobox{1}{1}{9}
		\end{tikzpicture}}.
		\]
		Therefore, $p(\lambda)=[4, 2]$ and the annihilator variety of $L(\lambda)$ is $V(\Ann (L(\lambda)))=\overline{\mathcal{O}}_{[4,2]}.$
	\end{ex}
	
	When $\lambda\in\hs$ is non-integral,  we write $\lambda=(\lambda_1,\dots,\lambda_n)$. Then we associate to $ \lambda $ a set $ S(\lambda) $ of some Young tableaux as follows. Let $ \lambda_Y: \lambda_{i_1}, \lambda_{i_2}, \dots, \lambda_{i_r} $  be a maximal subsequence of $ \lambda_1,\lambda_2,\dots,\lambda_n $ such that $ \lambda_{i_k} $, $ 1\leq k\leq r $ are congruent to each other by $ \mathbb{Z} $. Then the Young tableau $P(\lambda_Y)$ associated to the subsequence $ \lambda_Y $ using R-S  algorithm is a Young tableau in $ S(\lambda) $.
	
	Given two Young tableaux $P(\lambda_{Y_1})$ and $P(\lambda_{Y_2})$, write $P(\lambda_{Y_1}){\stackrel{c}{\sqcup}} P(\lambda_{Y_2})$ for
	the Young tableau whose multiset of nonzero column lengths equals the union of
	those of $P(\lambda_{Y_1})$ and $P(\lambda_{Y_2})$. Also write $2P(\lambda_{Y}) =P(\lambda_{Y}){\stackrel{c}{\sqcup}} P(\lambda_{Y})$. Correspondingly the shape of $P(\lambda_{Y_1}){\stackrel{c}{\sqcup}} P(\lambda_{Y_2})$ is $\sh(P(\lambda_{Y_1}){\stackrel{c}{\sqcup}} P(\lambda_{Y_2}))=p(\lambda_{Y_1}){\stackrel{c}{\sqcup}} p(\lambda_{Y_2})$, and $\sh(2P(\lambda_{Y}))=2p(\lambda_{Y})$.

	Then we  have the following theorem.
	
	\begin{Thm}\label{typeA}
		Let $\mathfrak{g}=\mathfrak{sl}(n, \mathbb{C})$.	Suppose $\lambda\in \mathfrak{h}^*$. Then
		$$V(\Ann (L(\lambda)))=\overline{\mathcal{O}}_{p(\lambda)},$$
		where $ p(\lambda)$ is the partition of the Young tableau $$P(\lambda)=\underset{P(\lambda_{Y})\in S(\lambda)}{\stackrel{c}{\sqcup}} P(\lambda_{Y}).$$

	\end{Thm}
	
	\begin{proof}
		For any $\lambda\in\mathfrak{h}^*$, there exists a unique anti-dominant weight $\mu\in\hs$ and a unique $w\in W_{[\mu]}^J$ such that $\lambda=w\mu$.
		Denote $\tilde{\pi}_w=j_{W_{[\lambda]}}^W(\pi_w)$. From Proposition \ref{AV}, we have $$V(\Ann (L(w\mu)))=\overline{\mathcal{O}}_{\tilde{\pi}_w}.$$
		
		When $\lambda$ is integral,  we know that $p(\lambda)$ and $p(w)$ have the same partitions from  \cite[Lemma 4.5]{BaiX}. From \cite[Theorem 10.1.1]{CM}, we have $\pi_w=\pi_{p(\lambda)}$ since $w$ belongs to the double cell corresponding to the partition $p(\lambda)$.
		
		In general, we suppose $\lambda=(\lambda_1,\dots,\lambda_n)$ is divided into several parts $\lambda_{Y_1},\dots,\lambda_{Y_m}$, with lengths being $r_1, r_2,\dots,r_m$. So we have $$W_{[\lambda]}=S_{r_1}\times S_{r_2}\times\dots\times S_{r_m}. $$
		From the definition of $w$, there exists a unique minimal length element $w_i\in S_{r_i}$ for $1\leq i\leq m$ such that $w=w_1\times w_2\times \dots\times w_m$ and $w_i^{-1}\lambda_{Y_i}$ is anti-dominant.
		Thus we have $$\pi_w=\pi_{w_1}\otimes \dots\otimes \pi_{w_m}.$$
		
	Similarly we have	$\pi_{w_i}=\pi_{p(\lambda_{Y_i})}$ since $p(w_i)=p(\lambda_{Y_i})$.
		From Lemma \ref{j.2}, we have
		$$\tilde{\pi}_w=j_{W_{[\lambda]}}^W(\pi_w)=j_{\prod\limits_{i}S_{r_i}}^W(\pi_{w_1}\otimes \dots\otimes \pi_{w_m})=\pi_{{{\stackrel{c}{\sqcup}}_i{p(\lambda_{Y_i})}}}.$$
		
		Let $ p(\lambda)$ be the partition of the Young tableau $P(\lambda)=\underset{1\leq i\leq m}{\stackrel{c}{\sqcup}} P(\lambda_{Y_i}).$
		Obviously we have $p(\lambda)=\stackrel{c}{\sqcup}_i{p(\lambda_{Y_i})}$ since $p(\lambda_{Y_i})$ is the partition of the Young tableau $P(\lambda_{Y_i})$.
		
		Thus we have completed the proof.
	\end{proof}

	\begin{ex}Let $\mathfrak{g}=\mathfrak{sl}(9, \mathbb{C})$.
		Suppose $\lambda=(2, -0.9, 4, 6, 4.1, 9,2.1, -29, 8)$. Then we can take $\lambda_{Y_1}=(2, 4, 6, 9, -29, 8)$ and $\lambda_{Y_2}=(-0.9,4.1,2.1)$. So we get
		\[
		P(\lambda_{Y_1})=
		\small{\begin{tikzpicture}[scale=\domscale+0.2,baseline=-20pt]
		\hobox{0}{0}{-29}
		\hobox{1}{0}{4}
		\hobox{2}{0}{6}
		\hobox{3}{0}{8}
		\hobox{0}{1}{2}
		\hobox{1}{1}{9}
		\end{tikzpicture}}\quad \text{and }\quad
		P(\lambda_{Y_2})=
		\begin{tikzpicture}[scale=\domscale+0.2,baseline=-20pt]
		\hobox{0}{0}{-0.9}
		\hobox{1}{0}{2.1}
		\hobox{0}{1}{4.1}
		\end{tikzpicture}.
		\]
		Therefore, $p(\lambda_{Y_1})=[4, 2]$,  $p(\lambda_{Y_2})=[2, 1]$ and $p(\lambda)=p(\lambda_{Y_1}){\stackrel{c}{\sqcup}} p(\lambda_{Y_2})=[6, 3]$. So the annihilator variety of $L(\lambda)$ is $V(\Ann (L(\lambda)))=\overline{\mathcal{O}}_{[6,3]}.$
	\end{ex}

	\section{The  annihilator varieties for types $B_n$ and $C_n$}\label{BC}
	
	In this section we suppose  $\mathfrak{g}=\mathfrak{so}(2n+1, \mathbb{C})$ or $\mathfrak{sp}(n,\mathbb{C})$.
	The Weyl group of $\mathfrak{g}$ is $W_n$ which is described in \S  \ref{B}.
	
	A weight $ \mu\in\hs $ is called \textit{anti-dominant} if $ \bil{\mu}{\al} \notin\mathbb{Z}_{>0}$ for all $ \al\in\Phi^+ $. For any $\lambda\in\mathfrak{h}^*$, there exists a unique anti-dominant weight $\mu\in\hs$ and a unique $w\in W_{[\mu]}^J$ such that $\lambda=w\mu$.
	
	Recall that $p({}^-{w})=\sh(P({}^-{w}))$ is the shape of the Young tableau $P({}^-{w})$, which is a partition of $2n$.
	
	The partition types of nilpotent orbits of types $B_n$ and $C_n$ are given in the following propositions.
	
	\begin{Prop}[\cite{ge61,CM}]
		Nilpotent orbits in $\mathfrak{so}{(2n+1,\mathbb{C})}$ are in one-to-one correspondence with the set of partitions of $2n+1$ in which even parts occur with even multiplicity.
	\end{Prop}
	
	\begin{Prop}[\cite{ge61,CM}]
		Nilpotent orbits in $\mathfrak{sp}{(n,\mathbb{C})}$ are in one-to-one correspondence with the set of partitions of $2n$ in which odd parts occur with even multiplicity.
	\end{Prop}

	\begin{Thm}\label{zuheB}
		Suppose the root system $\Phi=B_n$ or $C_n$. For an integral weight $\lambda\in \mathfrak{h}^*$,  we can write $\lambda=w\mu$ for a unique $w\in W^J$ and a unique anti-dominant $\mu\in \mathfrak{h}^*$.
		Then $p(\lambda^-)=p({}^-{w})$. So we have
		$$V(\Ann (L(\lambda)))=V(\Ann (L_w))=V(I_w)=\overline{\mathcal{O}}_w.$$	
		From $p(\lambda^-)=p({}^-{w})$, we can get a $B$-symbol $\Lambda_B$ and a special symbol $\Lambda_{B}^s$. Then we can get the corresponding special nilpotent orbit ${\mathcal{O}}_{w}=\mathcal{O}_{\Lambda_B^s}$.
	\end{Thm}
\begin{proof}
    From \cite[Lemma 5.2]{BXX}, we have $p(\lambda^-)=p({}^-{w})$. The rest results followed from Proposition \ref{trans}.
\end{proof}

	\begin{ex}Let $\mathfrak{g}=\mathfrak{so}(11, \mathbb{C})$ or $\mathfrak{sp}{(5,\mathbb{C})}$.
		Suppose $\lambda=(2, -4, 4, 9,0)$. Then
		\[
		P(\lambda^-)=
		\begin{tikzpicture}[scale=\domscale+0.1,baseline=-27pt]
		\hobox{0}{0}{-9}
		\hobox{1}{0}{-4}
		\hobox{2}{0}{-2}
		\hobox{3}{0}{4}
		\hobox{0}{1}{-4}
		\hobox{1}{1}{0}
		\hobox{2}{1}{0}
		\hobox{0}{2}{2}
		\hobox{1}{2}{4}
		\hobox{2}{2}{9}
		\end{tikzpicture}.
		\]
		Therefore, $p(\lambda^-)=[4, 3, 3]$. We have a  strictly increasing sequence $
		(3,4,6).
		$ Thus we have a $B$-symbol $\Lambda= \begin{pmatrix}
		2~3\\1
		\end{pmatrix} $, which is not special. But from $\Lambda$, we have a unique special $B$-symbol $\Lambda_{B}^s= \begin{pmatrix}
		1~3\\2
		\end{pmatrix} $. The corresponding special partition is ${\bf d}=[5,3,3]$ for $\mathfrak{g}=\mathfrak{so}(11, \mathbb{C})$ and ${\bf d}=[4,4,2]$ for $\mathfrak{g}=\mathfrak{sp}(5, \mathbb{C})$. Thus we have $V(\Ann (L(\lambda)))=\overline{\mathcal{O}}_{\Lambda_{B}^s}=\overline{\mathcal{O}}_{[5,3,3]}$ for  $\mathfrak{g}=\mathfrak{so}(11, \mathbb{C})$ and $V(\Ann (L(\lambda)))=\overline{\mathcal{O}}_{\Lambda_{B}^s}=\overline{\mathcal{O}}_{[4,4,2]}$  for $\mathfrak{g}=\mathfrak{sp}(5, \mathbb{C})$.
		
	\end{ex}

	\subsection{Bipartiton algorithm}
	\hspace{0.5cm}

	When $\lambda\in\hs$ is non-integral,  we write $\lambda=(\lambda_1,\dots,\lambda_n)$. We can associate to $ \lambda $ a set $ P(\lambda) $ of some Young tableaux as follows. 
	
	Recall that for $\mathfrak{g}=\mathfrak{sp}(n,\mathbb{C}), \mathfrak{so}(2n,\mathbb{C}) $ or $ \mathfrak{so}(2n+1,\mathbb{C}) $, we use $[\lambda]=(\lambda)_{(0)}\cup (\lambda)_{(\frac{1}{2})}\cup [\lambda]_3 $ to denote the set of  maximal subsequences $ \lambda_Y $ of  $ \lambda $ such that any two entries of $ \lambda_Y $ have an integral  difference or sum.

	For $\lambda_Y\in (\lambda)_{(0)}\cup (\lambda)_{(\frac{1}{2})}$, we can get a Young tableau $P(\lambda_Y^-)$.
	For $\lambda_Y\in [\lambda]_3$, we can get a Young tableau $P(\tilde{\lambda}_Y)$.

	Then we  have the following theorem.
	
	\begin{Thm}\label{type-B}
		Let $\mathfrak{g}=\mathfrak{so}(2n+1, \mathbb{C})$. Suppose
		$\lambda\in \mathfrak{h}^*$ and
		$[\lambda]=(\lambda)_{(0)}\cup (\lambda)_{(\frac{1}{2})}\cup [\lambda]_3$
		with $[\lambda]_3=\{{\lambda}_{{Y}_1},\dots,{\lambda}_{{Y}_m}\}$. Then
		$$V(\Ann (L(\lambda)))=\overline{\mathcal{O}}_{p_{_B}(\lambda)},$$
		where $ p_{_B}(\lambda)$ is the partition obtained from the following
		bipartition
		$$({\bf d}_{0}{\stackrel{c}{\sqcup}}{\bf d}_{00}{\stackrel{c}{\sqcup}}_i{\bf d}_i, {\bf f}_{0}{\stackrel{c}{\sqcup}}{\bf f}_{00}{\stackrel{c}{\sqcup}}_i{\bf f}_i ).$$
		Here $ ({\bf d}_{0},{\bf f}_{0})$ is the $B$-type bipartition obtained
		from $p((\lambda)_{(0)}^-)$, $({\bf d}_{00},{\bf f}_{00})$ is the $B$-type
		bipartition obtained from $p((\lambda)_{(\frac{1}{2})}^-)$ and
		$ ({\bf d}_{i},{\bf f}_{i})$ is the $B$-type bipartition obtained from
		the partition $p(\tilde{\lambda}_{{Y}_i})$.

	\end{Thm}
	
	\begin{proof}Let  $\mathfrak{g}=\mathfrak{so}(2n+1, \mathbb{C})$.	
		For any $\lambda\in\mathfrak{h}^*$, there exists a unique anti-dominant weight $\mu\in\hs$ and a unique $w\in W_{[\mu]}^J$ such that $\lambda=w\mu$.
		Denote $\tilde{\pi}_w=j_{W_{[\lambda]}}^W(\pi_w)$. From Proposition \ref{AV}, we have $$V(\Ann (L(w\mu)))=\overline{\mathcal{O}}_{\tilde{\pi}_w}.$$
		
		When $\lambda$ is integral,  the result is given in Theorem \ref{zuheB}.
		
		In general, we suppose $\lambda=(\lambda_1,\dots,\lambda_n)$ is divided into several parts $(\lambda)_{(0)}=[\lambda]_1, (\lambda)_{(\frac{1}{2})}=[\lambda]_2,[\lambda]_3=\{{\lambda_{Y_1},\dots,\lambda_{Y_m}}\}$, with lengths being $a$, $b$ and $r_1, r_2,\dots,r_m$. So we have $$W_{[\lambda]}=W_a\times W_{b}\times S_{r_1}\times S_{r_2}\times\dots\times S_{r_m}. $$
		From the definition of $w$, there exists some unique minimal length elements $x_a\in W_a$, $x_b\in W_b$ and $w_i\in S_{r_i}$ for $1\leq i\leq m$ such that $$w=x_a\times x_b \times w_1\times w_2\times \dots\times w_m,$$ and $x_a^{-1}(\lambda)_{(0)}$, $x_b^{-1}(\lambda)_{(\frac{1}{2})}$, $w_i^{-1}\tilde{\lambda}_{Y_i}$ are anti-dominant.
		Thus we have $$\pi_w=\pi_{x_a}\otimes\pi_{x_b}\otimes\pi_{w_1}\otimes \dots\otimes \pi_{w_m}.$$
		Suppose from $(\lambda)_{(0)}$ (note that $p((\lambda)^-_{(0)})$=$p({}^-x_a)$ by Theorem \ref{zuheB}), by using R-S  algorithm, we can get a $B$-symbol $\Lambda_0$. The corresponding special symbol $\Lambda_{0}^s$ will give us a bipartition $({\bf d}_0,{\bf f}_0)$ and a special partition ${\bf d}_{B_a}$  of type $B_a$. Thus we have
		$\pi_{x_a}=\pi_{({\bf d}_0,{\bf f}_0)}$.
		
		Similarly, from $(\lambda)_{(\frac{1}{2})}$, we can get a bipartition $({\bf d}_{00},{\bf f}_{00})$ and a special partition ${\bf d}_{B_b}$ of type $B_b$. Thus we have
		$\pi_{x_b}=\pi_{({\bf d}_{00},{\bf f}_{00})}$.
		
		From $\tilde{\lambda}_{Y_i}$, by using R-S  algorithm, we can get a Young tableau $P(\tilde{\lambda}_{{Y}_i})$, which will give us a partition $p(\tilde{\lambda}_{{Y}_i})={\bf p}_i={[p_1,p_2,\dots,p_{r_i}]}$ of type $A_{r_i}$. Thus we have
		$\pi_{w_i}=\pi_{{\bf p}_i}$
		and from equation (\ref{A-B}) we can get \[  j^{W_{r_i}}_{S_{r_i}}\pi_{ w_i}=j^{W_{r_i}}_{\prod\limits_{k}{S_{p^*_{k}}}} \sgn=j^{W_{r_i}}_{\prod\limits_k{W_{p^*_{k}}}} (\bigotimes_{k}\pi_{([1^{[\frac{p^*_{k}+1}{2}]}],~ [1^{[\frac{p^*_k}{2}]}])})=
		\pi_{({\bf d_i},{\bf f_i})},\]
		with ${\bf d_i}^*=[[\frac{p^*_{1}+1}{2}],\dots,[\frac{p^*_{r_i}+1}{2}]]$ and ${\bf f_i}^*=[[\frac{p^*_{1}}{2}],\dots,[\frac{p^*_{r_i}}{2}]]$.


		

		From Lemma \ref{j.1} and \begin{align*}j_{W_a\times S_b}^{W_{a+b}}(\pi_1\otimes \pi_2)&=j_{W_a\times W_b}^{W_{a+b}}\circ j^{W_a\times W_b}_{W_{a}\times S_b}(\pi_1\otimes \pi_2)\\
		&=j_{W_a\times W_b}^{W_{a+b}} (j^{W_a }_{W_{a}}\pi_1\otimes j^{W_b }_{S_{b}}\pi_2)\\
		&=j_{W_a\times W_b}^{W_{a+b}} (\pi_1 \otimes j^{W_b }_{S_{b}}\pi_2)
		\end{align*}
		for any $
		\pi_1 \in \mathrm{Irr}(W_a)$ and any $\pi_2\in \mathrm{Irr}(S_b)$,
		we have
		\begin{align*}
		\tilde{\pi}_w&=j_{W_{[\lambda]}}^{W_n}(\pi_w)=j_{W_{[\lambda]}}^{W_n}(\pi_{x_a}\otimes\pi_{x_b}\otimes\pi_{w_1}\otimes \dots\otimes \pi_{w_m})\\
		&=\pi_{({\bf d}_{0}{\stackrel{c}{\sqcup}}{\bf d}_{00}{\stackrel{c}{\sqcup}}_i{\bf d}_i, {\bf f}_{0}{\stackrel{c}{\sqcup}}{\bf f}_{00}{\stackrel{c}{\sqcup}}_i{\bf f}_i )}.
		\end{align*}
		
		From the bipartition $({\bf d}_{0}{\stackrel{c}{\sqcup}}{\bf d}_{00}{\stackrel{c}{\sqcup}}_i{\bf d}_i, {\bf f}_{0}{\stackrel{c}{\sqcup}}{\bf f}_{00}{\stackrel{c}{\sqcup}}_i{\bf f}_i )$, we can get a partition $p_{_B}(\lambda)$ for $2n+1$. Then we have	$$V(\Ann (L(\lambda)))=\overline{\mathcal{O}}_{p_{_B}(\lambda)}.$$
		
		Thus we have completed the proof.
	\end{proof}
	\begin{Thm}\label{type-C}
		Let $\mathfrak{g}=\mathfrak{sp}(n, \mathbb{C})$.	Suppose $\lambda\in \mathfrak{h}^*$ and $[\lambda]=(\lambda)_{(0)}\cup (\lambda)_{(\frac{1}{2})}\cup [\lambda]_3$ with $[\lambda]_3=\{{\lambda}_{{Y}_1},\dots,{\lambda}_{{Y}_m}\}$.  Then
		$$V(\Ann (L(\lambda)))=\overline{\mathcal{O}}_{p_{_C}(\lambda)},$$
		where $ p_{_C}(\lambda)$ is the $C$-type partition obtained from the following bipartition $$({\bf d}_{0}{\stackrel{c}{\sqcup}}{\bf d}_{00}{\stackrel{c}{\sqcup}}_i{\bf d}_i, {\bf f}_{0}{\stackrel{c}{\sqcup}}{\bf f}_{00}{\stackrel{c}{\sqcup}}_i{\bf f}_i ).$$
		Here $ ({\bf d}_{0},{\bf f}_{0})$ is the $B$-type bipartition obtained from $p((\lambda)_{(0)}^-)$, $({\bf d}_{00},{\bf f}_{00})$ is the $C$-type bipartition obtained from $p((\lambda)_{(\frac{1}{2})}^-)$ and $ ({\bf d}_{i},{\bf f}_{i})$ is the $B$-type bipartition obtained from the partition $p(\tilde{\lambda}_{{Y}_i})$.

	\end{Thm}
	\begin{proof}Let  $\mathfrak{g}=\mathfrak{sp}(n, \mathbb{C})$. There are some differences with type $B_n$.
		%
		%
		Similarly, we have $$W_{[\lambda]}=W_a\times W'_{b}\times S_{r_1}\times S_{r_2}\times\dots\times S_{r_m}. $$
		From the definition of $w$, there exists some unique minimal length elements $x_a\in W_a$, $x_b\in W'_b$ and $w_i\in S_{r_i}$ for $1\leq i\leq m$ such that $$w=x_a\times x_b \times w_1\times w_2\times \dots\times w_m,$$ and $x_a^{-1}(\lambda)_{(0)}$, $x_b^{-1}(\lambda)_{(\frac{1}{2})}$, $w_i^{-1}\tilde{\lambda}_{Y_i}$ are anti-dominant.
		Thus we have $$\pi_w=\pi_{x_a}\otimes\pi_{x_b}\otimes\pi_{w_1}\otimes \dots\otimes \pi_{w_m}.$$
		Suppose from $(\lambda)_{(0)}$, we can get a $B$-symbol $\Lambda_0$. The corresponding special symbol $\Lambda_{0}^s$ will give us a bipartition $({\bf d}_0,{\bf f}_0)$ and a special partition ${\bf d}_{B_a}$  of type $B_a$. Thus we have
		$\pi_{x_a}=\pi_{({\bf d}_0,{\bf f}_0)}$.
		
		Similarly, from $(\lambda)_{(\frac{1}{2})}$, we can get a bipartition $({\bf d}'_{00},{\bf f}'_{00})$ and a special partition ${\bf d}_{b}$ of type $D_b$. Thus we have
		$\pi_{x_b}=\pi_{({\bf d}'_{00},{\bf f}'_{00})}$. By Proposition \ref{D-C},
		$j^{W_b}_{W'_b}\pi_{x_b}$ will correspond to the metaplectic special partition $(({\bf d}_{b}^*)_D)^*$. Suppose the corresponding $C$-type bipartition is $({\bf d}_{00},{\bf f}_{00})$.
		
		From $\tilde{\lambda}_{Y_i}$, by using R-S  algorithm, we can get a Young tableau $P(\tilde{\lambda}_{{Y}_i})$, which will give us a partition $p(\tilde{\lambda}_{{Y}_i})={\bf p}_i={[p_1,p_2,\dots,p_{r_i}]}$ for $A_{r_i}$. Thus we have
		$\pi_{w_i}=\pi_{{\bf p}_i}$
		and from equation (\ref{A-B}) we can get \[  j^{W_{r_i}}_{S_{r_i}}\pi_{ w_i}=j^{W_{r_i}}_{\prod\limits_{k}{S_{p^*_{k}}}} \sgn=j^{W_{r_i}}_{\prod\limits_k{W_{p^*_{k}}}} (\bigotimes_{k}\pi_{([1^{[\frac{p^*_{k}+1}{2}]}],~ [1^{[\frac{p^*_k}{2}]}])})=
		\pi_{({\bf d_i},{\bf f_i})},\]
		with ${\bf d_i}^*=[[\frac{p^*_{1}+1}{2}],\dots,[\frac{p^*_{r_i}+1}{2}]]$ and ${\bf f_i}^*=[[\frac{p^*_{1}}{2}],\dots,[\frac{p^*_{r_i}}{2}]]$.


		

		From Lemma \ref{j.1} and \begin{align*}j_{W_a\times W'_b}^{W_{a+b}}(\pi_1\otimes \pi_2)&=j_{W_a\times W_b}^{W_{a+b}}\circ j^{W_a\times W_b}_{W_{a}\times W'_b}(\pi_1\otimes \pi_2)\\
		&=j_{W_a\times W_b}^{W_{a+b}} (j^{W_a }_{W_{a}}\pi_1\otimes j^{W_b }_{W'_{b}}\pi_2)\\
		&=j_{W_a\times W_b}^{W_{a+b}} (\pi_1 \otimes j^{W_b }_{W'_{b}}\pi_2)
		\end{align*}
		for any $
		\pi_1 \in \mathrm{Irr}(W_a)$ and any $\pi_2\in \mathrm{Irr}(W'_b)$,
		we have
		\begin{align*}
		\tilde{\pi}_w&=j_{W_{[\lambda]}}^{W_n}(\pi_w)=j_{W_{[\lambda]}}^{W_n}(\pi_{x_a}\otimes\pi_{x_b}\otimes\pi_{w_1}\otimes \dots\otimes \pi_{w_m})\\
		&=\pi_{({\bf d}_{0}{\stackrel{c}{\sqcup}}{\bf d}_{00}{\stackrel{c}{\sqcup}}_i{\bf d}_i, {\bf f}_{0}{\stackrel{c}{\sqcup}}{\bf f}_{00}{\stackrel{c}{\sqcup}}_i{\bf f}_i )}.
		\end{align*}
		
		From the bipartition $({\bf d}_{0}{\stackrel{c}{\sqcup}}{\bf d}_{00}{\stackrel{c}{\sqcup}}_i{\bf d}_i, {\bf f}_{0}{\stackrel{c}{\sqcup}}{\bf f}_{00}{\stackrel{c}{\sqcup}}_i{\bf f}_i )$, we can get a $ C$-type partition $p_{_C}(\lambda)$ for $2n$. Then we have	
		$$V(\Ann (L(\lambda)))=\overline{\mathcal{O}}_{p_{_C}(\lambda)}.$$
		
		Thus we have completed the proof.
	\end{proof}
	
		
		


	\begin{ex}\label{Bex}
		Let $\mathfrak{g}=\mathfrak{so}(25, \mathbb{C})$.
		Suppose we have $$\lambda=(2.5,1.5,3.5,2.6,1.6,0.6,-0.6,-2,6,5,-1,0).$$
		From $(\lambda)_{(0)}=(-2,6,5,-1,0)$, we can get a bipartition $({\bf d}_{0},{\bf f}_{0})=([3],[1,1])$.
		From $(\lambda)_{(\frac{1}{2})}=(2.5,1.5,3.5)$, we can get a bipartition $({\bf d}_{00},{\bf f}_{00})=([1],[1,1])$.
		From $\lambda_{Y_1}=(2.6,1.6,0.6,-0.6)$, we have $\tilde{\lambda}_{Y_1}=(2.6,1.6,0.6,0.6)$. By R-S  algorithm, we can get a partition ${\bf p}_1=[2,1,1]$ and ${\bf p}^*_1=[3,1]$. Then we  can get a bipartition $({\bf d}_{1},{\bf f}_{1})=([2,1],[1])$.
		Thus we can get a bipartition $({\bf d}_{0}{\stackrel{c}{\sqcup}}{\bf d}_{00}{\stackrel{c}{\sqcup}}_i{\bf d}_i, {\bf f}_{0}{\stackrel{c}{\sqcup}}{\bf f}_{00}{\stackrel{c}{\sqcup}}_i{\bf f}_i )=([6,1],[3,2])$.
		The corresponding  partition is ${\bf d}=[13,5,3,3,1]$. Thus we have $V(\Ann (L(\lambda)))=\overline{\mathcal{O}}_{[13,5,3,3,1]}$ for  $\mathfrak{g}=\mathfrak{so}(25, \mathbb{C})$.

		If  we let $\mathfrak{g}=\mathfrak{sp}(12, \mathbb{C})$. Then we consider the same $\lambda$. 	From $(\lambda)_{(\frac{1}{2})}=(2.5,1.5,3.5)$, we can get a $ D$-type special partition ${\bf d}_b=[2,2,1,1]$. From the metaplectic special partition $(({\bf d}^*_b)_D)^*=[2,2,2]$, we can get a $C$-type bipartition $({\bf d}_{00},{\bf f}_{00})=([1,1],[1])$. 	Thus we can get a bipartition $({\bf d}_{0}{\stackrel{c}{\sqcup}}{\bf d}_{00}{\stackrel{c}{\sqcup}}_i{\bf d}_i, {\bf f}_{0}{\stackrel{c}{\sqcup}}{\bf f}_{00}{\stackrel{c}{\sqcup}}_i{\bf f}_i )=([6,2],[3,1])$.
		The corresponding $C$-type partition is ${\bf d}=[12,6,4,2]$. Thus we have $V(\Ann (L(\lambda)))=\overline{\mathcal{O}}_{[12,6,4,2]}$ for  $\mathfrak{g}=\mathfrak{sp}(12, \mathbb{C})$.

	\end{ex}

	\subsection{Partition algorithm}

	Now we introduce another algorithm to compute the annihilator varieties of highest weight modules. In the following, we identify partitions, Young diagrams with nilpotent orbits, and
		bipartitions with $\mathrm{Irr}(W_n)$. 
  
  Firstly we need some lemmas.
	
	\begin{Lem}\label{2A}
		Let  $\bfpp$ be a special partition of $2n$
		in types $C,D$ and $2n+1$ in type $B$.  Let $\bfda$ be a partition of $m$ in type $A$. We use $X$ to denote the type of corresponding partitions.   Then we have
		$$j^{W_{n+m}^X}_{W_{n}\times S_m}(E^X_{\bfdp}\otimes \pi_{\bfda})=E^X_{\bfdp \csqcup 2\bfda},$$
where  $E_{\bf p}$ is the Springer representation for the partition $\bfdp$ in Proposition \ref{som}.
	\end{Lem}	
	
	\begin{proof}
		The  result is obvious by comparing Sommers's result in Proposition \ref{som} and the equations (\ref{A-B}) and (\ref{A-D}).
	\end{proof}	
	
	\begin{Lem}\label{BD}
		Let  $\bfdp$ be a special partition of $2n+1$
		in type $B$ and $\bfpq$ be a special partition of $2m$ in type $D$.    Then we have
		$$j^{W_{m}}_{W'_{m}} E^D_{\bfpq_e}=E^B_{\bfpq_e}    $$ and
		$$j^{W_{n+m}}_{W_{n}\times W'_m}(E^B_{\bfdp_o}\otimes E^D_{\bfpq_e})=E^B_{\bfdp_o \csqcup \bfpq_e},$$
	where $\bfpp_{o}$ is the union of odd columns of $\bfpp$ and $\bfpq_{e}$ is the union of even columns of $\bfpq$.
 \end{Lem}	
	\begin{proof} Since $\bfpp$ is a special partition of $2n+1$ in type $B$, its dual partition $\bfpp^*$ is also a partition of type $B$.  Thus the even columns of $\bfpp$ occur with even multiplicity. We can write $\bfdp=\bfdp_o \csqcup 2\bfdp_e$, where $\bfpp_{o}$ is the union of odd columns of $\bfpp$ and $2\bfpp_{e}$ is the union of even columns of $\bfpp$ with even multiplicity.

 Similarly $\bfpq$ is a special partition of $2m$ in type $D$, its dual partition $\bfpq^*$ is  a partition of type $C$.  Thus the odd columns of $\bfpq$ occur with even multiplicity. We can write $\bfpq=\bfpq_e \csqcup 2\bfpq_o$, where $\bfpq_{e}$ is the union of even columns of $\bfpq$ and $2\bfpq_{o}$ is the union of odd columns of $\bfpq$ with even multiplicity.
		
  The two equalities are very clear by inspecting  Sommers's result in Proposition \ref{som}.	
		
	\end{proof}

	\begin{Thm}\label{b-bipartition}
		Suppose $\fgg = \so(2n+1, \mathbb{C})$, $\lambda\in \fhh^*$ and
		$[\lambda]=(\lambda)_{(0)} \cup (\lambda)_{(\frac{1}{2})}\cup [\lambda]_3$
		with $[\lambda]_3=\{{\lambda}_{{Y}_1},\dots,{\lambda}_{{Y}_m}\}$. Let
		\begin{enumerate}
			\item $\bfpp_{0}$ be the $B$-type special partition associated to
			$(\lambda)_{(0)}$;
			\item ${\bf p}_{\frac{1}{2}}$ be the $C$-type special partition associated to
			$(\lambda)_{(\frac{1}{2})}$;
			\item $\bfpp_{i}$ be the $A$-type partition associated to
			$\tilde{\lambda}_{Y_i}$.
		\end{enumerate}
		Let $\bfpp_{\lambda}$ be the $B$-collapse of
		\begin{equation}\label{eq:rawB}
		\bfdd_{\lambda} := \bfpp_{0} \csqcup {\bf p}_{\frac{1}{2}} \csqcup (\csqcup_i 2 \bfpp_i).
		\end{equation}
		Then we have
		\[
		V(\Ann (L(\lambda)))=\overline{\mathcal{O}}_{\bfpp_{\lambda}}
		\]
	\end{Thm}
	\begin{proof}
		
		By  Theorem \ref{type-B}, $V(\Ann (L(\lambda))) = \overline{\cO}$ such that
		the Springer representation of $\cO$ equals to
		\[
		E_{\cO} := ({\bf d}_{0}\csqcup{\bf d}_{00}\csqcup_i{\bf d}_i, {\bf f}_{0}\csqcup{\bf f}_{00}\csqcup_i{\bf f}_i ).
		\]

		Now  $\bfpp_{0}$ is the $B$-type special partition associated to
		$(\lambda)_{(0)}$, so we can write
		$\bfpp_{0}=\bfpp_{o} \csqcup 2\bfpp_{e}$ where $\bfpp_{o}$ is the union of odd columns of $\bfpp_{0}$ and $2\bfpp_{e}$ is the union of even columns of $\bfpp_{0}$ with even multiplicity.
		
		Similarly, ${\bf p}_{\frac{1}{2}}$ is the $C$-type special partition associated to
		$(\lambda)_{(\frac{1}{2})}$, so we can write
		$\bfpp_{\frac{1}{2}}=\bfpp_{\frac{e}{2}} \csqcup 2\bfpp_{\frac{o}{2}}$ where $\bfpp_{\frac{e}{2}}$ is the union of even columns of $\bfpp_{\frac{1}{2}}$ and $2\bfpp_{\frac{o}{2}}$ is the union of odd columns of $\bfpp_{\frac{1}{2}}$ with even multiplicity.
  
Denote $2p=|{\bf p}_{\frac{1}{2}}|$.
By \S \ref{d-weyl} and the proof in Theorem \ref{type-B}, we also have 	$\pi_{({\bf d}_{0},{\bf f}_{0})}=E^B_{\bfpp_{0}}$, 	$\pi_{({\bf d}_{00},{\bf f}_{00})}=j_{W_p^D}^{W_p^{B}}E^D_{{\bf p}_{\frac{1}{2}}}=E^B_{{\bf p}_{\frac{1}{2}}}$ and $\pi_{({\bf d}_{i},{\bf f}_{i})}=j^{W_{r_i}}_{S_{r_i}}\pi_{ \bfpp_i}$.

	On the other hand, from the proof in Theorem \ref{type-B} and Proposition \ref{d-induciton}, we also have $E_{\cO} =j_{W_{[\lambda]}}^{W_n}(E^B_{\bfpp_{0}}\otimes E^B_{{\bf p}_{\frac{1}{2}}}\otimes\pi_{\bfpp_{1}}\otimes \dots\otimes \pi_{\bfpp_{m}})$.  Recall that $W_{[\lambda]}=W_a\times W_{b}\times S_{r_1}\times S_{r_2}\times\dots\times S_{r_m}$. We denote $W_{[\lambda]'}=W_a\times W'_{b}\times S_{r_1}\times S_{r_2}\times\dots\times S_{r_m}$.

		Thus by Lemma \ref{2A} and Lemma \ref{BD}, we have
		\begin{align*}
		&j_{W_{[\lambda]}}^{W_n}(E^B_{\bfpp_{0}}\otimes E^B_{{\bf p}_{\frac{1}{2}}}\otimes\pi_{\bfpp_{1}}\otimes \dots\otimes \pi_{\bfpp_{m}})\\
		=&j_{W_{[\lambda]'}}^{W_n}(E^B_{\bfpp_{o} \csqcup 2\bfpp_{e}}\otimes E^D_{\bfpp_{\frac{e}{2}} \csqcup 2\bfpp_{\frac{o}{2}}}\otimes\pi_{\bfpp_{1}}\otimes \dots\otimes \pi_{\bfpp_{m}})\\
		=&j_{W_{[\lambda]'}}^{W_n}\left(j^{W_{p_o+p_e}}_{W_{p_o}\times S_{p_e}}(E^B_{\bfpp_{o}}\otimes \pi_{\bfpp_{e}})\otimes j^{W'_{p'_{e}+p'_o}}_{W'_{p'_e}\times S_{p'_o}}(E^D_{\bfpp_{\frac{e}{2}}}\otimes \pi_{\bfpp_{\frac{o}{2}}})\otimes\pi_{\bfpp_{1}}\otimes \dots\otimes \pi_{\bfpp_{m}}\right)\\
  & j_{W_{[\lambda]'}}^{W_n}\left(j^{W_{p_o+p_e}\times W'_{p'_{e}+p'_o}}_{W_{p_o}\times S_{p_e}\times W'_{p'_e}\times S_{p'_o} }(E^B_{\bfpp_{o}}\otimes \pi_{\bfpp_{e}}\otimes E^D_{\bfpp_{\frac{e}{2}}}\otimes \pi_{\bfpp_{\frac{o}{2}}})\otimes \pi_{\bfpp_{1}}\otimes \dots\otimes \pi_{\bfpp_{m}}\right)\\
=& j_{W_{[\lambda]'}}^{W_n}\left(j^{W_{p_o+p_e}\times W'_{p'_{e}+p'_o}}_{W_{p_o}\times W'_{p'_e}\times S_{p_e}\times S_{p'_o} }(E^B_{\bfpp_{o}}\otimes E^D_{\bfpp_{\frac{e}{2}}}\otimes \pi_{\bfpp_{e}}\otimes \pi_{\bfpp_{\frac{o}{2}}})\otimes \pi_{\bfpp_{1}}\otimes \dots\otimes \pi_{\bfpp_{m}}\right)\\	
=&j_{W_{p_o}\times W'_{p'_e}\times S_{p_e}\times S_{p'_o} \times S_{p_1}\times\dots \times S_{p_m}}^{W_n}\left(E^B_{\bfpp_{o}}\otimes E^D_{\bfpp_{\frac{e}{2}}}\otimes\pi_{\bfpp_{e}}\otimes \pi_{\bfpp_{\frac{o}{2}}}  \otimes\pi_{\bfpp_{1}}\otimes \dots\otimes \pi_{\bfpp_{m}}\right)\\
  =&j_{W_{p_o+p'_e}\times S_{p_e+p'_o+p_1+\dots+  p_m}}^{W_n}\left(j_{W_{p_o}\times W'_{p'_e}}^{W_{p_o+p'_e}}(E^B_{\bfpp_{o}}\otimes E^D_{\bfpp_{\frac{e}{2}}})\otimes j^{S_{p_e+p'_o+p_1+\dots + p_m}}_{S_{p_e}\times S_{p'_o} \times S_{p_1}\times\dots \times S_{p_m} }(\pi_{\bfpp_{e}}\otimes \pi_{\bfpp_{\frac{o}{2}}}  \otimes\pi_{\bfpp_{1}}\otimes \dots\otimes \pi_{\bfpp_{m}})\right)\\
  =&j_{W_{p_o+p'_e}\times S_{p_e+p'_o+p_1+\dots+  p_m}}^{W_n}\left(E^B_{{\bfpp_{o}}\csqcup {\bfpp_{\frac{e}{2}}}}\otimes\pi_{\bfpp_{e}\csqcup \bfpp_{\frac{o}{2}}\csqcup_i\bfpp_{i}}\right)\\
		=&E^B_{{\bfpp_{o}}\csqcup {\bfpp_{\frac{e}{2}}}\csqcup 2\bfpp_{e}\csqcup 2\bfpp_{\frac{o}{2}}\csqcup(\csqcup_i2\bfpp_{i}) }\\
		=&E^B_{{\bfpp_0}\csqcup {\bfpp_{\frac{1}{2}}}\csqcup(\csqcup_i 2\bfpp_{i} )}\\
		=&E_{\bfdd_{\lambda}},
		\end{align*}
		where $p_o=\frac{1}{2}|\bfpp_{o}|$, $p_e=|\bfpp_{e}|$,  $p'_o=|\bfpp_{\frac{o}{2}}|$, $p'_e= \frac{1}{2}|\bfpp_{\frac{e}{2}}|$ and $p_i=|\bfpp_{i}|$ for $1\leq i\leq m$.
		
		

		So	following the recipe in Proposition \ref{som},
		we have proved that the Weyl group representation $E_{\bfdd_{\lambda}}$ attached to $\bfdd_{\lambda}$
		equals to $E_{\cO}$.
		
		%
		%
		By Sommers \cite[Lemma~9]{Sommers},
		$E_{\bfdd_{\lambda}} = E_{\bfpp_{\lambda}}$ since $\bfpp_{\lambda}$ and
		$\bfdd_{\lambda}$ have the same $B$-collapse (which is $\bfpp_{\lambda}$).
		Now by the injectivity of Springer correspondence, we conclude that
		$\cO = \bfpp_{\lambda}$.
	\end{proof}
	
	\begin{Thm}
		Suppose $\fgg = \sp(n, \mathbb{C})$, $\lambda\in \fhh^*$ and
		$[\lambda]=(\lambda)_{(0)} \cup (\lambda)_{(\frac{1}{2})}\cup [\lambda]_3$
		with $[\lambda]_3=\{{\lambda}_{{Y}_1},\dots,{\lambda}_{{Y}_m}\}$. Let
		\begin{enumerate}
			\item $\bfpp_{0}$ be the $C$-type special partition associated to
			$(\lambda)_{(0)}$;
			\item ${\bf p}_{\frac{1}{2}}$ be the $C$-type metaplectic special partition associated to
			$(\lambda)_{(\frac{1}{2})}$;
			\item $\bfpp_{i}$ be the $A$-type partition associated to
			$\tilde{\lambda}_{Y_i}$.
		\end{enumerate}
		Let $\bfpp_{\lambda}$ be the $C$-collapse of
		\begin{equation}\label{eq:rawc}
		\bfdd_{\lambda} := \bfpp_{0} \csqcup {\bf p}_{\frac{1}{2}} \csqcup (\csqcup_i 2 \bfpp_i).
		\end{equation}
		Then we have
		\[
		V(\Ann (L(\lambda)))=\overline{\mathcal{O}}_{\bfpp_{\lambda}}
		\]
	\end{Thm}
	
	\begin{proof}

		By our Theorem \ref{type-C}, $V(\Ann (L(\lambda))) = \overline{\cO}$ such that
		the Springer representation of $\cO$ equals to
		\[
		E_{\cO} := ({\bf d}_{0}\csqcup{\bf d}_{00}\csqcup_i{\bf d}_i, {\bf f}_{0}\csqcup{\bf f}_{00}\csqcup_i{\bf f}_i ).
		\]

		Now  $\bfpp_{0}$ is the $C$-type special partition associated to
		$(\lambda)_{(0)}$, so we can write
		$\bfpp_{0}=\bfpp_{e} \csqcup 2\bfpp_{o}$ where $\bfpp_{e}$ is the union of even columns of $\bfpp_{0}$ and $2\bfpp_{o}$ is the union of odd columns of $\bfpp_{0}$ with even multiplicity.

		From $(\lambda)_{(\frac{1}{2})}$, we can get  a special partition ${\bf d}_{\frac{1}{2}}$ of type $D$.  From Proposition \ref{D-C}, we
		know	${\bf p}_{\frac{1}{2}}$ will be the metaplectic special partition $(({\bf d}_{\frac{1}{2}}^*)_D)^*$.
		So $\bfpp_{\frac{1}{2}}^*=({\bf d}_{\frac{1}{2}}^*)_D$ is a partition of type $D$  in which  even rows occur with even multiplicity. Thus we can write
		$\bfpp_{\frac{1}{2}}=\bfpp_{\frac{o}{2}} \csqcup 2\bfpp_{\frac{e}{2}}$ where $\bfpp_{\frac{o}{2}}$ is the union of odd columns of ${\bf p}_{\frac{1}{2}}$ and $2\bfpp_{\frac{e}{2}}$ is the union of even columns of ${\bf p}_{\frac{1}{2}}$ with even multiplicity.

Denote $2p=|{\bf p}_{\frac{1}{2}}|$.
By \S \ref{d-weyl} and the proof in Theorem \ref{type-B}, we also have 	$\pi_{({\bf d}_{0},{\bf f}_{0})}=E^C_{\bfpp_{0}}$, 	$\pi_{({\bf d}_{00},{\bf f}_{00})}=j_{W_p^D}^{W_p^{C}}E^D_{{\bf d}_{\frac{1}{2}}}=E^C_{{\bf p}_{\frac{1}{2}}}$ and $\pi_{({\bf d}_{i},{\bf f}_{i})}=j^{W_{r_i}}_{S_{r_i}}\pi_{ \bfpp_i}$.

		On the other hand, from the proof in Theorem \ref{type-C} and Proposition \ref{d-induciton}, we also have $E_{\cO} =j_{W_{[\lambda]}}^{W_n}(E^C_{\bfpp_{0}}\otimes E^C_{{\bf p}_{\frac{1}{2}}}\otimes\pi_{\bfpp_{1}}\otimes \dots\otimes \pi_{\bfpp_{m}})$.

		Thus by Lemma \ref{2A} and Lemma \ref{BD}, similarly we have
		\begin{align*}
		&j_{W_{[\lambda]}}^{W_n}(E^C_{\bfpp_{0}}\otimes E^C_{{\bf p}_{\frac{1}{2}}}\otimes\pi_{\bfpp_{1}}\otimes \dots\otimes \pi_{\bfpp_{m}})\\
		=&E^C_{{\bfpp_0}\csqcup {\bfpp_{\frac{1}{2}}}\csqcup(\csqcup_i 2\bfpp_{i}) }\\
		=&E_{\bfdd_{\lambda}}.
		\end{align*}

		The rest arguments are the same with Theorem \ref{b-bipartition}.

		
	\end{proof}

	\begin{Rem}
		For  $\mathfrak{g}=\mathfrak{so}(2n+1, \mathbb{C})$, $\lambda\in \mathfrak{h}^*$ and $[\lambda]=(\lambda)_{(0)}\cup (\lambda)_{(\frac{1}{2})}\cup [\lambda]_3$ with $[\lambda]_3=\{{\lambda}_{{Y}_1},\dots,{\lambda}_{{Y}_m}\}$, if $(\lambda)_{(0)}=\emptyset$,  $\bfpp_{\lambda}$ will be the partition obtained from the $B$-collapse of    ${\bf p}_0{\stackrel {c}{\sqcup}}{\bf p}_{\frac{1}{2}}\csqcup(\csqcup_i 2\bfpp_{i})$, where ${\bf p}_0=[1]$ is the trivial partition. 	If $(\lambda)_{(0)}=(\lambda)_{(\frac{1}{2})}=\emptyset$,  $\bfpp_{\lambda}$ will be the partition obtained from the $B$-collapse of    ${\bf p}_0\csqcup(\csqcup_i 2\bfpp_{i})$.

	\end{Rem}

	\begin{ex}
		In the above Example \ref{Bex} for  $\mathfrak{g}=\mathfrak{so}(2n+1, \mathbb{C})$, from $(\lambda)_{(0)}=(-2,6,5,-1,0)$, we can get a $B$-type special partition ${\bf p}_{0}=[7,1,1,1,1]$.
		From $(\lambda)_{(\frac{1}{2})}=(2.5,1.5,3.5)$, we can get a $C$-type special partition ${\bf p}_{\frac{1}{2}}=[2,2,1,1,]$.
		From $\lambda_{Y_1}=(2.6,1.6,0.6,-0.6)$, we can get a $A$-type partition ${\bf p}_{1}=[2,1,1]$. Thus we have ${\bf p}_{0}{\stackrel{c}{\sqcup}}{\bf p}_{00}{\stackrel {c}{\sqcup}}2{\bf p}_1=[13,5,4,2,1]$. Its $B$-collapse is $\bfpp_{\lambda}=[13,5,3,3,1]$. Then we have the same result
		$V(\Ann (L(\lambda)))=\overline{\mathcal{O}}_{[13,5,3,3,1]}$.
		
		For  $\mathfrak{g}=\mathfrak{sp}(n, \mathbb{C})$, from $(\lambda)_{(0)}=(-2,6,5,-1,0)$, we can get a $C$-type special partition ${\bf p}_{0}=[6,2,1,1]$.
		From $(\lambda)_{(\frac{1}{2})}=(2.5,1.5,3.5)$, we can get a $C$-type metaplectic special partition ${\bf p}_{\frac{1}{2}}=[2,2,2]$.
		From $\lambda_{Y_1}=(2.6,1.6,0.6,-0.6)$, we can get a $A$-type partition ${\bf p}_{1}=[2,1,1]$. Thus we have ${\bf p}_{0}{\stackrel{c}{\sqcup}}{\bf p}_{\frac{1}{2}}{\stackrel {c}{\sqcup}}2{\bf p}_1=[12,6,5,1]$. Its $C$-collapse is ${\bfpp_{\lambda}}=[12,6,4,2]$. Then we have the same result
		$V(\Ann (L(\lambda)))=\overline{\mathcal{O}}_{[12,6,4,2]}$.

	\end{ex}

	\section{The  annihilator varieties for type $D_n$}\label{Dchar}
	
	In this section we suppose  $\mathfrak{g}=\mathfrak{so}(2n, \mathbb{C})$.
	The Weyl group of $\mathfrak{g}$ is $W'_n$ which is described in \S  \ref{d-weyl}.
	The partition type of nilpotent orbits of type $D_n$ is given in the following.
	
	\begin{Prop}[\cite{ge61,CM}]
		Nilpotent orbits in $\mathfrak{so}{(2n,\mathbb{C})}$ are in one-to-one correspondence with the set of partitions of $2n$ in which even parts occur with even multiplicity,
		except that each ``very even" partition ${\bf d}$ (consisting of only even parts) correspond to two orbits, denoted by $\mathcal{O}^I_{\bf d}$ and $\mathcal{O}^{II}_{\bf d}$.
	\end{Prop}
\subsection{Lusztig-Spaltenstein induction}	
	
	Firstly we  recall some properties about Lusztig-Spaltenstein induction of nilpotent orbits, some details can be found in \cite[\S 7]{CM}. Let
	${\mathfrak{g}}$ be a complex simple Lie algebra. Let $\mathfrak{p}$ be a parabolic subalgebra with Levi decomposition $\mathfrak{p}=\mathfrak{ l}\oplus\mathfrak{u}$.  Then we have a natural projection
	$P\colon {\mathfrak{g}}^{\ast}\rightarrow \mathfrak{p}^{\ast}$.
	
	\begin{Prop}
		[{\cite[Theorem 7.1.1]{CM}}] Let $\mathfrak{g}$ be a complex simple Lie algebra with a Cartan decomposition $\mathfrak{g}=\mathfrak{n}\oplus\mathfrak{h}\oplus\mathfrak{n}^-$. Let ${\mathcal{O}%
		}_{\mathfrak{l}}\subset\mathfrak{l}^{\ast}$ be a nilpotent orbit. Then there
		exists a unique nilpotent orbit ${\mathcal{O}}_{{\mathfrak{g}}}$ that meets
		$P^{-1}({\mathcal{O}}_{\mathfrak{l}})$ in an open dense subset. We have
		\[
		\dim{\mathcal{O}}_{{\mathfrak{g}}}=\dim{\mathcal{O}}_{\mathfrak{l}}+2\dim
		\mathfrak{n}.
		\]
		
	\end{Prop}

	Note that	the orbit ${\mathcal{O}}_{{\mathfrak{g}}}:=\mathrm{Ind}_{\mathfrak{l}%
	}^{{\mathfrak{g}}}({\mathcal{O}}_{\mathfrak{l}})$ is  called the {\it induced orbit}
	of ${\mathcal{O}}_{\mathfrak{l}}$.
	
	Recall that in \S \ref{j-induction}, for a given $W$-module $\pi$, we use $\mathcal{O}_{{\pi}}$ to denote the nilpotent orbit with trivial local system via the Springer correspondence. For a highest weight module $L(\lambda)$, we write $\lambda=w\mu$ for a unique $w\in W_{[\lambda]}^J$ and a unique anti-dominant $\mu\in \mathfrak{h}^*$.
	From Proposition \ref{intehral-W}, we have $\mathrm{Springer}({\tilde{\pi}_w})={\mathcal{O}}_{\tilde{\pi}_w}$, where $\tilde{\pi}_w=j_{W_{[\lambda]}}^W(\pi_w)$ and $\pi_w$	is the  special representation of $W_{[\lambda]}$ spanned by some Goldie-rank polynomials. From the simple system $\Delta_{[\lambda]}$, we can get a Levi subalgebra $\mathfrak{l}$. Under the Springer correspondence, we have a nilpotent orbit ${\mathcal{O}}_{{\pi}_w}$ in $\mathfrak{l}^*$ (or $\mathfrak{l}$). From \cite{LS} or \cite[Theorem 2.1.3]{GS}, we have the following.

 \begin{Prop}
		\label{ind}Under the above conditions, we have $\mathcal{O}_{\tilde{\pi}}=\mathrm{Ind}_{\mathfrak{l}}^{{\mathfrak{g}}}\left( {\mathcal{O}}_{{\pi}_w}
		\right)  $.
	\end{Prop}

		
	
	\begin{Prop}
		[{\cite[Proposition 7.1.4]{CM}} ]\label{trans-l}
		Let $\mathfrak{g}$ be a complex simple Lie algebra. Let   $\mathfrak{l}_1$ and  $ \mathfrak{ l}_2$ be two Levi subalgebras of $\mathfrak{g}$ with  $\mathfrak{l}_1 \subset \mathfrak{ l}_2$. 
	Then we have
$\mathrm{Ind}_{\mathfrak{l}_2}^{{\mathfrak{g}}}(\mathrm{Ind}_{\mathfrak{l}_1}^{\mathfrak{l}_2}\left(  \mathcal{O}_{\mathfrak{l}_1}
		\right))  $.
	\end{Prop}
Since this transitivity of induction, we may focus on maximal Levi subalgebras. 
From \cite[\S 7.3]{CM}, we know that each maximal Levi subalgebra takes the following form:
$$\mathfrak{l}=\mathfrak{gl}_l\oplus \mathfrak{g}'$$
where $\mathfrak{g}'$ is classical and of the same type as $\mathfrak{g}$. From \cite[Lemma 3.8.1]{CM} and $\mathfrak{l}=\mathfrak{l}_I$ for some subset $I\subset \Delta$, we know that $\Delta\setminus I$  contains a unique simple root.

 Let $\mathcal{O}_{\mathfrak{l}}$ be a nilpotent orbit in a maximal Levi subalgebra of $\mathfrak{g}$. From  \cite[\S 7.3]{CM}, we can write $\mathcal{O}_{\mathfrak{l}}=\mathcal{O}_{\mathbf{d}}\oplus \mathcal{O}_{\mathbf{f}} $,
 where $\mathcal{O}_{\mathbf{d}}$ is a nilpotent orbit in $\mathfrak{sl}_l$ and  $\mathcal{O}_{\mathbf{f}}$ is a nilpotent orbit in $\mathfrak{g}'$. Then we have the following.

 \begin{Prop}[{{\cite[Theorem 7.3.3 and Corollary 7.3.4]{CM}}} and {\cite[p.89]{Lus84}}]\label{typeX}
  Let $\mathfrak{g}$ be a complex simple Lie algebra of type $X$ with $X=B,C$ or $D$.   Let $\mathcal{O}_{\mathfrak{l}}$ be a nilpotent orbit in a maximal Levi subalgebra $\mathfrak{l}=\mathfrak{l}_I$ of $\mathfrak{g}$ and ${\mathcal{O}}=\mathrm{Ind}_{\mathfrak{l}
	}^{{\mathfrak{g}}}({\mathcal{O}}_{\mathfrak{l}})$. Suppose the partition of $\mathcal{O}$ is $\mathbf{p}$. Then we have the following:
 \begin{enumerate}
     \item If $\mathfrak{g}=\mathfrak{so}(2n,\mathbb{C})$ with $n=2m'$ being even and $\mathbf{p}$ is a very even partition of $2n$, but the simple root $\alpha_i~ 
  (\in \Delta \setminus I)\neq \alpha_{n-1}$ or $\alpha_n$, then $\mathcal{O}_{\mathbf{f}}$ is  very even  and the numeral of $\mathcal{O}$ is the same as that of $\mathcal{O}_{\mathbf{f}}$;
     \item If $\mathfrak{g}=\mathfrak{so}(2n,\mathbb{C})$ with $n=2m'$ being even,  $\mathbf{p}$ is a very even partition of $2n$, and the simple root $\alpha_i~ 
  (\in \Delta \setminus I)= \alpha_{n-1}$ or $\alpha_n$, 
   then we have   $\mathcal{O}=\mathcal{O}^I$ if $\alpha_{n}\in \Delta\setminus I$ and $\mathcal{O}=\mathcal{O}^{II}$ if $\alpha_{n-1}\in \Delta\setminus I$.
 \end{enumerate}
 \end{Prop}

\subsection{Very even orbit}
Now we recall a concept called {\it $\tau$-invariant}. Some details can be found in \cite{V79}.
	\begin{Defff}
		For $w \in W_{[\lambda]}$, the $\tau$-invariant of $w$ is defined to be $$\tau (w)=\{\alpha \in \Delta_{[\lambda]}  \mid w\alpha  \notin  \Phi_{[\lambda]}^+     \}.$$
	\end{Defff}

	From  \cite{V79}, $\tau (w)$ depends only on the primitive ideal $I_{w\lambda}$ when $\lambda$ is antidominant and regular. So the Borho-Jantzen-Duflo $\tau$-invariant of $I_{w\lambda}$ is $\tau(I_{w\lambda}):=\tau(w)$.
	
	Recall that we can write $\lambda=w\mu$ for a unique $w\in W^J$ and a unique anti-dominant $\mu\in \mathfrak{h}^*$. 
	The following lemma will tell us how to find such  $w$.


	\begin{Lem}\label{w'}
		Suppose $ \lambda \in \mathfrak{h}^* $ is  integral with $  \lambda=(\lambda_1,\dots,\lambda_n)$. Then we have the following:
		\begin{itemize}
			\item [(1)]There is a unique element $ w\in W_n $ such that
			\begin{itemize}
				\item[(i)] if $ \lambda_i\neq 0 $ then $ \sigma_i' $ has the same sign with $ \lambda_i $;
				\item [(ii)]if $ \lambda_i= 0 $ then $ \sigma_i' <0$;
				\item[(iii)] if $ \lambda_i<\lambda_j $, then $ \sigma_i'<\sigma_j'  $;
				\item[(iv)] if $ \lambda_i=\lambda_j $ and $ i<j $, then $ \sigma_i'<\sigma_j'  $.
			\end{itemize}
			where $ (\sigma_1',\sigma_2',\dots, \sigma_n')=(-w^{-1}(n), -w^{-1}(n-1),\dots, -w^{-1}(1)) $;
			\item [(2)] Let $ w $ be the element in (1) and $t=(1,-1)\in W_n$. Then one of $ w $ and $ wt $ belongs to $ W_n' $, and we denote it by $ w' $. Then $ w' $ is precisely the minimal length element such that $ w'^{-1}\lambda $ is antidominant.
\end{itemize}
\end{Lem}
\begin{proof}
     See the proof in Lemma 5.2 and Lemma 5.3 in \cite{BXX}.
 \end{proof}

\begin{ex}
Let $ \mf{g}=\mathfrak{so}(12, \mathbb{C})$ and $ \lambda=(2,-1,1,-3,1,-1) $. Then from Proposition \ref{anti}, the antidominant weight in $ W\lambda $ is $ \mu= (-3,-2,-1,-1,-1,1) $. We have a string diagram:
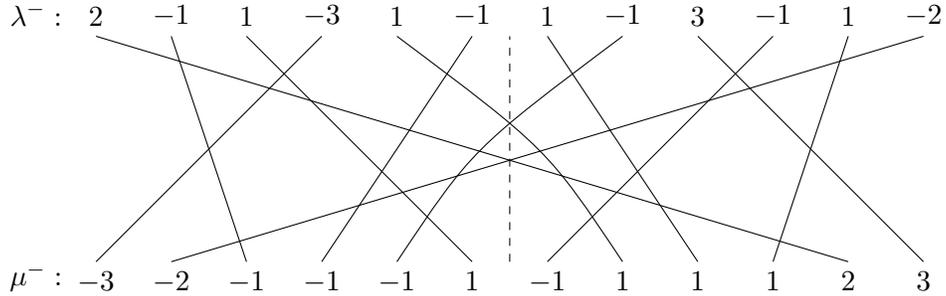
\begin{center}
\begin{tikzpicture}
 \draw(0.2,-0.2) node{$ \mu^-: $};
\filldraw (1,0) \nodd node[below]{$ -3 $};  \draw (1,0)--(4,3);
\filldraw (2,0) \nodd node[below]{$ -2$};  \draw (2,0)--(12,3);
\filldraw (3,0) \nodd node[below]{$ -1 $};	\draw (3,0)--(2,3);
\filldraw (4,0) \nodd node[below]{$ -1 $};	\draw (4,0)--(6,3);
\filldraw (5,0) \nodd node[below]{$ -1 $};	\draw (5,0)..controls (6,1.5)..(8,3);
\filldraw (6,0) \nodd node[below]{$ 1 $};		\draw (6,0)--(3,3);

\filldraw (7,0) \nodd node[below]{$ -1 $}; \draw (7,0)--(10,3);
\filldraw (8,0) \nodd node[below]{$ 1$};	\draw (8,0)..controls (7,1.5)..(5,3);
\filldraw (9,0) \nodd node[below]{$ 1 $};	\draw (9,0)--(7,3);
\filldraw (10,0) \nodd node[below]{$ 1 $};	\draw (10,0)--(11,3);
\filldraw (11,0) \nodd node[below]{$ 2 $};	\draw (11,0)->(1,3);
\filldraw (12,0) \nodd node[below]{$ 3 $}; \draw (12,0)->(9,3);
\draw[dashed](6.5,0)--(6.5,3);
 \draw (0.2,3.3) node{$ \lambda^-: $};
\filldraw (1,3) \nodd node[above]{$2 $};
\filldraw (2,3) \nodd node[above]{$ -1$};
\filldraw (3,3) \nodd node[above]{$ 1 $};
\filldraw (4,3) \nodd node[above]{$ -3 $};
\filldraw (5,3) \nodd node[above]{$ 1 $};
\filldraw (6,3) \nodd node[above]{$ -1 $};

\filldraw (7,3) \nodd node[above]{$ 1 $};
\filldraw (8,3) \nodd node[above]{$ -1$};
\filldraw (9,3) \nodd node[above]{$ 3 $};
\filldraw (10,3) \nodd node[above]{$ -1 $};
\filldraw (11,3) \nodd node[above]{$ 1 $};
\filldraw (12,3) \nodd node[above]{$ -2 $};
\end{tikzpicture}.
\end{center}
By Lemma \ref{w'}, this diagram represents an element $ w'=(4,-2,1,5,-6,3)\in W_6' $.
\end{ex}

 Note that the Weyl group $ (W,S) $ of type $ D_n $ is isomorphic to $ (W_n', S_n') $  via
\begin{equation}\label{isom}
    s_{\ep_i-\ep_{i+1}}\mapsto s_{n-i},1\leq  i\leq n-1,\text{ and } s_{\ep_{n-1}+\ep_n}\mapsto u.
\end{equation}
Then $s_{n-i} (i \geq 1)$ acts on $\mathfrak{h}^*$  by exchanging
 the coefficients of $\varepsilon_i$ and $\ep_{i+1}$, and $u$ acts on $\mathfrak{h}^*$  by  exchanging the coefficients of $\varepsilon_{n-1}$ and $\ep_{n}$, then changing the signs of them.

 \begin{Lem}\label{tau-i}
   Let $w\in W_n'$.   Then we have
   $$\tau(w)=\{\alpha_{n-i}\mid w(i)>w(i+1)\},$$
   where $w(0):=-w(2)$ and $w(n+1):=0$.

 \end{Lem}
 \begin{proof}
     The result follows  from \cite[ Proposition 8.2.2]{BB05} and the above isomorphism (\ref{isom}).
 \end{proof}

\begin{Lem}\label{w12}Let $w\in W'_n$ with $n=2m$ and $\mathcal{O}_w$ be a very even nilpotent orbit of $2n=4m$ with the partition
     ${\bf p}=[2m_1,2m_1,2m_2,2m_2,\dots,2m_k,2m_k]^*$. Then we have the following:
    \begin{enumerate}
        \item $\mathcal{O}_w$ will be type I if and only if $w \stackrel{LR}{\sim} w_I$ with \begin{align*}
          w_I=&(2m_1,2m_1-1,\dots,1,\\
          &\ 2m_2+2m_1,2m_2-1+2m_1,\dots,1+2m_1,\\
          &\ \dots, \ \sum_{1\leq i \leq k}2m_i,\dots,1+\sum_{1\leq i \leq k-1}2m_i);
        \end{align*}
         \item $\mathcal{O}_w$ will be type II if and only if $w \stackrel{LR}{\sim} w_{II}$ with \begin{align*}
          w_{II}=&(-2m_1,2m_1-1,2m_1-2,\dots,2,-1,\\
          &\ 2m_2+2m_1,2m_2-1+2m_1,\dots,1+2m_1,\\
          &\ \dots, \ \sum_{1\leq i \leq k}2m_i,\dots,1+\sum_{1\leq i \leq k-1}2m_i).
          \end{align*}
    \end{enumerate}
\end{Lem}
\begin{proof}
    By the R-S algorithm, we can find that $w_{I}$ and $w_{II}$ will give us the same shape of Young tableaux, i.e., $p({}^-{w_I})=\sh(P({}^-{w_I}))=\sh(P({}^-{w_{II}}))=p({}^-{w_{II}})$. And the numbers of boxes in the columns of $p({}^-{w_{I}})=p({}^-{w_{II}})$ are: $2m_1,2m_1,2m_2,2m_2,\dots,2m_k,2m_k$. Thus we have $p({}^-{w_{I}})=p({}^-{w_{II}})={\bf p}$.

    From Lemma \ref{tau-i}, we find that the simple root $e_{n-1}-e_n\in \tau(w_I)$ and  $e_{n-1}+e_n\in \tau(w_{II})$. So by \cite{V79}, $\mathcal{O}_{w_I}$ is type I and $\mathcal{O}_{w_{II}}$ is type II.

    Thus $\mathcal{O}_w$ will be type I if and only if $w \stackrel{LR}{\sim} w_I$, and vice versa.
\end{proof}

From now on, the Weyl group elements in the form of some $w_I$ and $w_{II}$ are called {\it very even elements of type I and type II} respectively.

 \begin{Rem}In the algorithm of domino tableaux, two elements $w\stackrel{LR}{\sim} x$ and
 their corresponding special domino tableaux  have the same very even shape if and only if the number of vertical dominos in the corresponding domino tableau $T_L(w)$ is congruent to the number of vertical dominoes in $T_L(x)$ modulo $4$, see \cite[\S 3]{Mc96}.
     The domino tableau $T_L(w_I)$ doesn't contain any vertical domino and $T_L(w_{II})$ contains only two vertical dominos. So they correspond to different types of  very even orbits. Some basic definitions and properties about domino tableaux can be found in \cite{BXX} or Garfinkle's work \cite{Garfinkle1,Garfinkle2,Garfinkle3,Ga-D}. There is a web page designed by Garfinkle which can be used to check that two elements of a Weyl group belong to the same left cell, double cell or not: https://devragj.github.io/.
 \end{Rem}

	Now we have the following.
	
	\begin{Thm}\label{zuheD}
		Suppose $\mathfrak{g}=\mathfrak{so}(2n, \mathbb{C})$. For an integral weight $\lambda\in \mathfrak{h}^*$,  we can write $\lambda=w\mu$ for a unique $w\in W_n'$ and a unique anti-dominant $\mu\in \mathfrak{h}^*$.
		Then $p(\lambda^-)=p({}^-{w})$  or $p(\lambda^-)=p({}^-({w}t))$. So we have
		$$V(\Ann (L(\lambda)))=V(\Ann (L_w))=V(I_w)=\overline{\mathcal{O}}_w.$$	
		From $p(\lambda^-)$, we can get a $D$-symbol $\Lambda_D$ and a special symbol $\Lambda_{D}^s$. Then we can get the corresponding special nilpotent orbit ${\mathcal{O}}_{w}=\mathcal{O}_{\Lambda_{D}^s}$.
		When $n=2m'$ is even and ${\mathcal{O}}_{w}=\mathcal{O}_{\Lambda_{D}^s}$ is a very even orbit, it will be type I if  there exists some very even element $ w_I\stackrel{LR}{\sim} w$  and type II if  there exists some very even element $ w_{II}\stackrel{LR}{\sim} w$.
	\end{Thm}

 
 \begin{proof}
		From \cite[Lemma 5.3]{BXX}, we have $p(\lambda^-)=p({}^-{w})$  or $p(\lambda^-)=p({}^-({w}t))$.
  So we only need to prove the case when  $n=2m$ is even and ${\mathcal{O}}_{w}=\mathcal{O}_{\Lambda_{D}^s}$ is a very even orbit.

Thus the result follows from Lemma \ref{w12}.
	\end{proof}

	\begin{ex}\label{type1}
		Let $\mathfrak{g}=\mathfrak{so}(12, \mathbb{C})$.
		Suppose $\lambda=(3,2,1, -5, -6, 7)$. Then
		\[
		P(\lambda^-)=
		\begin{tikzpicture}[scale=\domscale+0.1,baseline=-55pt]
		\hobox{0}{0}{-7}
		\hobox{1}{0}{-3}
		\hobox{0}{1}{-6}
		\hobox{1}{1}{-2}
		\hobox{0}{2}{-5}
		\hobox{1}{2}{-1}
		\hobox{0}{3}{1}
		\hobox{1}{3}{5}
  \hobox{0}{4}{2}
  \hobox{1}{4}{6}
  \hobox{0}{5}{3}
  \hobox{1}{5}{7}
		\end{tikzpicture}.
		\]
		Therefore, $p(\lambda^-)=[2^6]$. The corresponding special partition is ${\bf d}=[2^6]$. Thus we have $$V(\Ann (L(\lambda)))=\overline{\mathcal{O}}_{\mathbf{d}}=\overline{\mathcal{O}}_{[2^6]}.$$
		We can write $\lambda=w\mu$ for a unique $w\in W_6'$ and a unique anti-dominant $\mu\in \mathfrak{h}^*$. In our example, $w=(-4,-5,-6,3,2,-1)$ and $\mu=(-7,-6,-5,-3,-2,-1)$. We  find that $w\stackrel{LR}{\sim} w_{II}=(-6,5,4,3,2,-1)$.
  Thus we have $$V(\Ann (L(\lambda)))=\overline{\mathcal{O}}^{II}_{\mathbf{d}}=\overline{\mathcal{O}}^{II}_{[2^6]}.$$
	\end{ex}

\subsection{Bipartition algorithm}	
 
 In general, we have the following.

	\begin{Thm}\label{typeD}
		Let $\mathfrak{g}=\mathfrak{so}(2n, \mathbb{C})$.	Suppose $\lambda\in \mathfrak{h}^*$ and $[\lambda]=(\lambda)_{(0)}\cup (\lambda)_{(\frac{1}{2})}\cup [\lambda]_3$ with $[\lambda]_3=\{{\lambda}_{{Y}_1},\dots,{\lambda}_{{Y}_m}\}$.  Then
		$$V(\Ann (L(\lambda)))=\overline{\mathcal{O}}_{p_{_D}(\lambda)},$$
		where $ p_{_D}(\lambda)$ is the partition obtained from the following unordered bipartition  $$\{{\bf d}_{0}{\stackrel{c}{\sqcup}}{\bf d}_{00}{\stackrel{c}{\sqcup}}_i{\bf d}_i, {\bf f}_{0}{\stackrel{c}{\sqcup}}{\bf f}_{00}{\stackrel{c}{\sqcup}}_i{\bf f}_i \}.$$
		Here $ \{{\bf d}_{0},{\bf f}_{0}\}$ is the $D$-type unordered bipartition  obtained from $(\lambda)_{(0)}$, $\{{\bf d}_{00},{\bf f}_{00}\}$ is the $D$-type unordered bipartition  obtained from $(\lambda)_{(\frac{1}{2})}$ and $ \{{\bf d}_{i},{\bf f}_{i}\}$ is the $D$-type unordered bipartition  obtained from $\tilde{\lambda}_{{Y}_i}$. 
  When $n=2m'$ is even and ${\mathcal{O}}_{p_{_D}(\lambda)}$ is a very even orbit, it will be type I  if $k(\lambda)\equiv 0 
  ~(\mathrm{mod}~ {2})$ and type II if $k(\lambda)\equiv 1 
  ~(\mathrm{mod}~ {2})$. Here
	 we use $k(\lambda)$ to denote the number of very even unordered bipartitions with numeral II in the set of  $\{\{{\bf d}_{0},{\bf f}_{0}\}, \{{\bf d}_{00},{\bf f}_{00}\},  \{{\bf d}_{i},{\bf f}_{i}\}| 1\leq i\leq m  \}$. 
 \end{Thm}
	\begin{proof}
		For any $\lambda\in\mathfrak{h}^*$, there exists a unique anti-dominant weight $\mu\in\hs$ and a unique $w\in W_{[\mu]}^J$ such that $\lambda=w\mu$.
		Denote $\tilde{\pi}_w=j_{W_{[\lambda]}}^W(\pi_w)$. From Proposition \ref{AV}, we have $$V(\Ann (L(w\mu)))=\overline{\mathcal{O}}_{\tilde{\pi}_w}.$$
		
		When $\lambda$ is integral,  the result is given in Theorem \ref{zuheD}.
		
		In general, we suppose $\lambda=(\lambda_1,\dots,\lambda_n)$ is divided into several parts $(\lambda)_{(0)}=[\lambda]_1, (\lambda)_{(\frac{1}{2})}=[\lambda]_2,[\lambda]_3=\{{\lambda_{Y_1},\dots,\lambda_{Y_m}}\}$, with lengths being $a$, $b$ and $r_1, r_2,\dots,r_m$. So we have $$W_{[\lambda]}=W'_a\times W'_{b}\times S_{r_1}\times S_{r_2}\times\dots\times S_{r_m}. $$
		From the definition of $w$, there exists some unique minimal length elements $x_a\in W'_a$, $x_b\in W'_b$ and $w_i\in S_{r_i}$ for $1\leq i\leq m$ such that $$w=x_a\times x_b \times w_1\times w_2\times \dots\times w_m,$$ and $x_a^{-1}(\lambda)_{(0)}$, $x_b^{-1}(\lambda)_{(\frac{1}{2})}$, $w_i^{-1}\tilde{\lambda}_{Y_i}$ are anti-dominant.
		Thus we have $$\pi_w=\pi_{x_a}\otimes\pi_{x_b}\otimes\pi_{w_1}\otimes \dots\otimes \pi_{w_m}.$$
		Suppose from $(\lambda)_{(0)}$, we can get a $D$-symbol $\Lambda_0$. The corresponding special symbol $\Lambda_{0}^s$ will give us an unordered bipartition $\{{\bf d}_0,{\bf f}_0\}$ and a special partition ${\bf p}_{0}$ of type  $D_a$. Thus we have
		$\pi_{x_a}=\pi_{\{{\bf d}_0,{\bf f}_0\}}$.
		
		Similarly, from $(\lambda)_{(\frac{1}{2})}$, we can get an unordered bipartition  $\{{\bf d}_{00},{\bf f}_{00}\}$ and a special partition ${\bf d}_{\frac{1}{2}}$ of type $D_b$. Thus we have
		$\pi_{x_b}=\pi_{\{{\bf d}_{00},{\bf f}_{00}\}}$.
		
		From $\tilde{\lambda}_{Y_i}$, by using R-S  algorithm, we can get a Young tableau $P(\tilde{\lambda}_{{Y}_i})$, which will give us a partition $p(\tilde{\lambda}_{{Y}_i})={\bf p}_i={[p_1,p_2,\dots,p_{r_i}]}$ for $A_{r_i}$. Thus we have
		$\pi_{w_i}=\pi_{{\bf p}_i}$
		and from equation (\ref{A-D}) we can get \[  j^{W'_{r_i}}_{S_{r_i}}\pi_{ w_i}=j^{W'_{r_i}}_{\prod\limits_{k}{S_{p^*_{k}}}} \sgn=j^{W'_{{r_i}}}_{\prod\limits_k{W'_{p^*_{k}}}} (\bigotimes_{k}\pi_{\{[1^{[\frac{p^*_{k}}{2}]}],~ [1^{[\frac{p^*_k+1}{2}]}]\}})=
		\pi_{\{{\bf d_i},{\bf f_i}\}},\]
		with ${\bf d_i}^*=[[\frac{p^*_{1}}{2}],\dots,[\frac{p^*_{r_i}}{2}]]$ and ${\bf f_i}^*=[[\frac{p^*_{1}+1}{2}],\dots,[\frac{p^*_{r_i}+1}{2}]]$.


		

		From Proposition \ref{d-induciton}
		we have
		\begin{align*}
		\tilde{\pi}_w&=j_{W_{[\lambda]}}^{W_n}(\pi_w)=j_{W_{[\lambda]}}^{W_n}(\pi_{x_a}\otimes\pi_{x_b}\otimes\pi_{w_1}\otimes \dots\otimes \pi_{w_m})\\
		&=\pi_{\{{\bf d}_{0}{\stackrel{c}{\sqcup}}{\bf d}_{00}{\stackrel{c}{\sqcup}}_i{\bf d}_i, {\bf f}_{0}{\stackrel{c}{\sqcup}}{\bf f}_{00}{\stackrel{c}{\sqcup}}_i{\bf f}_i \}}.
		\end{align*}
		
		From the unordered bipartition  $\{{\bf d}_{0}{\stackrel{c}{\sqcup}}{\bf d}_{00}{\stackrel{c}{\sqcup}}_i{\bf d}_i, {\bf f}_{0}{\stackrel{c}{\sqcup}}{\bf f}_{00}{\stackrel{c}{\sqcup}}_i{\bf f}_i \}$, we can get a partition $p_{_D}(\lambda)$ for $2n$. Then we have	$$V(\Ann (L(\lambda)))=\overline{\mathcal{O}}_{p_{_D}(\lambda)}.$$
		
		
 When $n=2m'$ is even and ${\mathcal{O}}_{p_D(\lambda)}$ is a very even orbit,	
the numeral of it follows from Proposition \ref{d-induciton}.		
		
		Thus we have completed the proof.
	\end{proof}
	\subsection{Partition algorithm}	
	
	Similarly we have the following partition algorithm.
	
	\begin{Thm}\label{partition-algorithm-}
		Suppose $\fgg = \so(2n, \mathbb{C})$, $\lambda\in \fhh^*$ and
		$[\lambda]=(\lambda)_{(0)} \cup (\lambda)_{(\frac{1}{2})}\cup [\lambda]_3$
		with $[\lambda]_3=\{{\lambda}_{{Y}_1},\dots,{\lambda}_{{Y}_m}\}$. Let
		\begin{enumerate}
			\item $\bfpp_{0}$ be the $D$-type special partition associated to
			$(\lambda)_{(0)}$;
			\item ${\bf p}_{\frac{1}{2}}$ be the $C$-type metaplectic special partition associated to
			$(\lambda)_{(\frac{1}{2})}$;
			\item $\bfpp_{i}$ be the $A$-type partition associated to
			$\tilde{\lambda}_{Y_i}$.
		\end{enumerate}
		Let $\bfpp_{\lambda}$ be the $D$-collapse of
		\begin{equation}\label{eq:rawd}
		\bfdd_{\lambda} := \bfpp_{0} \csqcup {\bf p}_{\frac{1}{2}} \csqcup (\csqcup_i 2 \bfpp_i).
		\end{equation}
		Then we have
		\[
		V(\Ann (L(\lambda)))=\overline{\mathcal{O}}_{\bfpp_{\lambda}}
		\]
	\end{Thm}
	\begin{proof}
	
		By our Theorem \ref{typeD}, $V(\Ann (L(\lambda))) = \overline{\cO}$ such that
		the Springer representation of $\cO$ equals to
		\[
		E_{\cO} := ({\bf d}_{0}\csqcup{\bf d}_{00}\csqcup_i{\bf d}_i, {\bf f}_{0}\csqcup{\bf f}_{00}\csqcup_i{\bf f}_i ).
		\]

		Now  $\bfpp_{0}$ is the $D$-type special partition associated to
		$(\lambda)_{(0)}$, so we can write
		$\bfpp_{0}=\bfpp_{e} \csqcup 2\bfpp_{o}$ where $\bfpp_{e}$ is the union of even columns of $\bfpp_{0}$ and $2\bfpp_{o}$ is the union of odd columns of $\bfpp_{0}$ with even multiplicity.

		From $(\lambda)_{(\frac{1}{2})}$, we can get  a special partition ${\bf d}_{\frac{1}{2}}$ of type $D$.  From Proposition \ref{D-C}, we
		know	${\bf p}_{\frac{1}{2}}$ will be the metaplectic special partition $(({\bf d}_{\frac{1}{2}}^*)_D)^*$.
		So $\bfpp_{\frac{1}{2}}^*=({\bf d}_{\frac{1}{2}}^*)_D$ is a partition of type $D$  in which  even rows occur with even multiplicity. Thus we can write
		$\bfpp_{\frac{1}{2}}=\bfpp_{\frac{o}{2}} \csqcup 2\bfpp_{\frac{e}{2}}$ where $\bfpp_{\frac{o}{2}}$ is the union of odd columns of ${\bf p}_{\frac{1}{2}}$ and $2\bfpp_{\frac{e}{2}}$ is the union of even columns of ${\bf p}_{\frac{1}{2}}$ with even multiplicity.

Denote $2p=|{\bf p}_{\frac{1}{2}}|$. By \S \ref{d-weyl} and the proof in Theorem \ref{typeD}, we also have 	$\pi_{({\bf d}_{0},{\bf f}_{0})}=E^D_{\bfpp_{0}}$, 	$\pi_{({\bf d}_{00},{\bf f}_{00})}=E^D_{{\bf d}_{\frac{1}{2}}}=j_{W^C_p}^{W^D_p}j_{W^D_p}^{W^C_p}E^D_{{\bf d}_{\frac{1}{2}}}=j_{W^C_p}^{W^D_p}{\bf p}_{\frac{1}{2}}$ and $\pi_{({\bf d}_{i},{\bf f}_{i})}=j^{W_{r_i}}_{S_{r_i}}\pi_{ \bfpp_i}$.

		On the other hand, from the proof in Theorem \ref{typeD} and Proposition \ref{d-induciton}, we also have $E_{\cO} =j_{W_{[\lambda]}}^{W_n}(E^D_{\bfpp_{0}}\otimes E^D_{{\bf d}_{\frac{1}{2}}}\otimes\pi_{\bfpp_{1}}\otimes \dots\otimes \pi_{\bfpp_{m}})$. Recall that $W_{[\lambda]}=W^D_a\times W^D_{b}\times S_{r_1}\times S_{r_2}\times\dots\times S_{r_m}$. We denote $W_{[\lambda]'}=W^D_a\times W^C_{b}\times S_{r_1}\times S_{r_2}\times\dots\times S_{r_m}$.

		The left arguments are similar to type $B$. We omit them here.
	\end{proof}
	
	%
	%
	
	\begin{Thm}\label{I-II}
 Suppose $\fgg = \so(2n, \mathbb{C})$ and $\lambda\in \fhh^*$,  we can write $\lambda=w\mu$ for a unique $w\in W_{[\mu]}^J$ and a unique anti-dominant $\mu\in \mathfrak{h}^*$. Let $\lambda_D$ be the subsequence of $\lambda$ consisting of all the entries in $(\lambda)_{(0)} \cup (\lambda)_{(\frac{1}{2})}$. Let $\lambda_A$ be the subsequence of $\lambda$ consisting of all the entries in $[\lambda]_3$. We use $\mathbf{f}$ to denote the $D$-collapse of the partition $\bfpp_{0} \csqcup {\bf p}_{\frac{1}{2}}$ and $\mathbf{d}$ to denote the  the partition $\csqcup_i \bfpp_i$. We denote $d=|\mathbf{d}|$ and $2f=|\mathbf{f}|$. Thus $d+f=n$. So we can write $\lambda_A=(\lambda_{i_1},\dots,\lambda_{i_d})$ and $\lambda_D=(\lambda_{j_1},\dots,\lambda_{j_{f}})$. We denote $\bar{\lambda}=(\lambda_{A},\lambda_{D})=(\lambda_{i_1},\dots,\lambda_{i_d},\lambda_{j_1},\dots,\lambda_{j_{f}})$.
   Then we have $$V(\Ann (L(\Bar{\lambda})))=V(\Ann (L(\lambda)))=\overline{\mathcal{O}}_{{\bf p}_{\lambda}},$$ 
   where  $\mathcal{O}_{{\bf p}_{\lambda}}=\mathrm{Ind}_{\mathfrak{l}
 }^{{\mathfrak{g}}}({\mathcal{O}}_{\mathfrak{l}})$ with $\mathcal{O}_{\mathfrak{l}}=\mathcal{O}_{\mathbf{d}}\oplus \mathcal{O}_{\mathbf{f}} $ and $\mathfrak{l}$ is the maximal Levi subalgebra corresponding to the simple root $\alpha_d$ when $d\leq n-2$.
  
 When $n=2m'$ is even and $V(\Ann (L(\lambda)))=\overline{\mathcal{O}}_{{\bf p}_{\lambda}}$ is a very even nilpotent orbit, then 
$\mathcal{O}_{\mathbf{f}}$    is a very even nilpotent orbit and 
 $\mathcal{O}_{{\bf p}_{\lambda}}$ is labelled in accordance with Proposition \ref{d-induciton} and Proposition \ref{typeX}.
  
\end{Thm}
 \begin{proof}
The result follows from Proposition \ref{ind}, \ref{trans-l} and \ref{typeX}.

In general  we can rewrite $\lambda_D$ as $$\bar{\lambda}_D=(\lambda_{s_1},\dots,\lambda_{s_{f_1}},\lambda_{s_{f_1+1}},\dots,\lambda_{s_{f}}),$$
where the first $f_1$ elements are all the ones in $(\lambda)_{(0)}$. Thus we can regard $\bar{\lambda}_{D}$ as a weight for a Lie algebra $\mathfrak{g}'$ of type $D_f$. 
Then we also have $$V(\Ann (L(\Bar{\lambda}_D)))=V(\Ann (L(\lambda_D)))=\overline{\mathcal{O}}_{{\bf f}},$$ 
and   $\mathcal{O}_{\bf f}= (\bfpp_{0} \csqcup {\bf p}_{\frac{1}{2}})_{_D}$.
When $n=2m'$ is even and $V(\Ann (L(\lambda)))=\overline{\mathcal{O}}_{{\bf p}_{\lambda}}$ is a very even nilpotent orbit, then 
$\mathcal{O}_{\mathbf{f}}$    is a very even nilpotent orbit. From Proposition \ref{d-induciton} and \ref{typeX}, we can get the numeral of ${\mathcal{O}}_{{\bf p}_{\lambda}}$.

 \end{proof}

 

Recall that for  $ x=(\lambda_{i_1}, \lambda_{i_2},\dots \lambda_{i_r})\in[\lambda]_3 $, let $  y=(\lambda_{j_1}, \lambda_{j_2},\dots, \lambda_{j_p}) $ be the maximal subsequence of $ x $ such that $ j_1=i_1 $ and the difference of any two entries of $ y$ is an integer. Let $ z= (\lambda_{k_1}, \lambda_{k_2},\dots, \lambda_{k_q}) $ be the subsequence obtained by deleting $ y$ from $ x $, which is possible empty.
	We have
	  $\tilde{x}=(\lambda_{j_1}, \lambda_{j_2},\dots, \lambda_{j_p}, -\lambda_{k_q}, -\lambda_{k_{q-1}},\dots ,-\lambda_{k_1}).  $
We call $q$ the \emph{negative index} of $x$, denoted by $q(x)$. Note that the Weyl group will not change the conjugacy classes of Levi subalgebras by \cite[Lemma 3.3.1]{CM}. Then we have the following corollary.
\begin{Cor}Suppose $n=2m'$ is even, $V(\Ann (L(\lambda)))=\overline{\mathcal{O}}_{{\bf p}_{\lambda}}$ and ${\mathcal{O}}_{{\bf p}_{\lambda}}$ is a very even nilpotent orbit. 
   When $[\lambda]=[\lambda]_3$
		with $[\lambda]_3=\{{\lambda}_{{Y}_1},\dots,{\lambda}_{{Y}_m}\}$ (i.e., $(\lambda)_{(0)} = (\lambda)_{(\frac{1}{2})}=\emptyset$), 
  $\mathcal{O}_{{\bf p}_{\lambda}}=\mathcal{O}_{\csqcup_i 2\bfpp_i}$
  will be type I if $\sum_{1\leq i\leq m}q(\lambda_{Y_i})\equiv 0 
  ~(\mathrm{mod}~ {2})$ and   type II if $\sum_{1\leq i\leq m}q(\lambda_{Y_i})\equiv  1  ~(\mathrm{mod}~ {2})$.
		
\end{Cor}

 
From Proposition \ref{d-induciton}, we have the following corollary.

\begin{Cor}\label{nemural-d}
Suppose $n=2m'$ is even, $V(\Ann (L(\lambda)))=\overline{\mathcal{O}}_{{\bf p}_{\lambda}}$ and ${\mathcal{O}}_{{\bf p}_{\lambda}}$  is a very even nilpotent orbit. 
   When $[\lambda]=(\lambda)_{(0)} \cup (\lambda)_{(\frac{1}{2})}$
		with $[\lambda]_3=\emptyset$,  $\bfpp_{0}$ and  $ {\bf p}_{\frac{1}{2}}$ will be very even partitions. Then $\mathcal{O}_{{\bf p}_{\lambda}}=\mathcal{O}_{\bfpp_{0} \csqcup {\bf p}_{\frac{1}{2}}}$ will be type I if $\mathcal{O}_{\bfpp_{0}}$ and $\mathcal{O}_ {{\bf p}_{\frac{1}{2}}}$ have the same numeral, otherwise it will be type II.

   When $[\lambda]=(\lambda)_{(0)} \cup (\lambda)_{(\frac{1}{2})}$
		with $[\lambda]_3\neq \emptyset$, $\bfpp_{0}$ and $ {\bf p}_{\frac{1}{2}}$ will be   very even partitions. We use $\mathbf{f}$ to denote the  $D$-collapse of the partition $\bfpp_{0} \csqcup {\bf p}_{\frac{1}{2}}$ and $\mathbf{d}$ to denote the  the partition $\csqcup_i \bfpp_i$.
  Then  $\mathcal{O}_{{\bf p}_{\lambda}}$ will be type I if  $\mathcal{O}_{\mathbf{f}}$  and $\mathcal{O}_ {2\bf d}$ have the same numeral, otherwise it will be type II.
		
\end{Cor}

Now we use $k(\lambda)$ to denote the number of very even orbits with numeral II in the set of very even orbits of type $D$:    $\{\mathcal{O}_{\bfpp_{0}}, \mathcal{O}_ {{\bf p}_{\frac{1}{2}}},  \mathcal{O}_{2\bfpp_i}| 1\leq i\leq m  \}$. 
From the above two corollaries, we can see that  $\mathcal{O}_{{\bf p}_{\lambda}}$ will be type I if $k(\lambda)\equiv 0 
  ~(\mathrm{mod}~ {2})$ and type II if $k(\lambda)\equiv 1 
  ~(\mathrm{mod}~ {2})$.

	\begin{ex}
		Let $\mathfrak{g}=\mathfrak{so}(20, \mathbb{C})$.
		Suppose we have $$\lambda=(2.5,1.5,3.5,2.6,1.6,0.6,-0.6,-2,6,-5).$$
		From $(\lambda)_{(0)}=(-2,6,-5)$, we can get a $D$-type unordered bipartition $\{{\bf d}_{0},{\bf f}_{0}\}=\{[1,1],[1]\}$.
		From $(\lambda)_{(\frac{1}{2})}=(2.5,1.5,3.5)$, we can get a $D$-type unordered bipartition $\{{\bf d}_{00},{\bf f}_{00}\}=\{[1],[1,1]\}$.
		From $\lambda_{Y_1}=(2.6,1.6,0.6,-0.6)$, we have $\tilde{\lambda}_{Y_1}=(2.6,1.6,0.6,0.6)$. By R-S  algorithm, we can get a  partition ${\bf p}_1=[2,1,1]$ and ${\bf p}^*_1=[3,1]$. Then we  can get a $D$-type unordered bipartition $\{{\bf d}_{1},{\bf f}_{1}\}=\{[1],[2,1]\}$.
		Thus we can get a $D$-type unordered bipartition $\{{\bf d}_{0}{\stackrel{c}{\sqcup}}{\bf d}_{00}{\stackrel{c}{\sqcup}}_i{\bf d}_i, {\bf f}_{0}{\stackrel{c}{\sqcup}}{\bf f}_{00}{\stackrel{c}{\sqcup}}_i{\bf f}_i \}=\{[3],[4,3]\}$.
		The corresponding  partition is ${\bf d}=[7,7,5,1]$. Thus we have $V(\Ann (L(\lambda)))=\overline{\mathcal{O}}_{[7,7,5,1]}$ for  $\mathfrak{g}=\mathfrak{so}(20, \mathbb{C})$.

		From $(\lambda)_{(0)}=(-2,6,-5)$, we can get a $D$-type special partition ${\bf p}_{0}=[2,2,1,1]$.
		From $(\lambda)_{(\frac{1}{2})}=(2.5,1.5,3.5)$, we can get a $C$-type metaplectic special partition ${\bf p}_{\frac{1}{2}}=[2,2,2]$.
		From $\lambda_{Y_1}=(2.6,1.6,0.6,-0.6)$, we can get a $A$-type partition ${\bf p}_{1}=[2,1,1]$. Thus we have ${\bf p}_{0}{\stackrel{c}{\sqcup}}{\bf p}_{\frac{1}{2}}{\stackrel {c}{\sqcup}}2{\bf p}_1=[8,6,5,1]$. Its $D$-collapse is ${\bfpp_{\lambda}}=[7,7,5,1]$. Then we have the same result
		$V(\Ann (L(\lambda)))=\overline{\mathcal{O}}_{[7,7,5,1]}$.
		
	\end{ex}	

	\begin{ex}
		Let $\mathfrak{g}=\mathfrak{so}(24, \mathbb{C})$.
		Suppose we have $$\lambda=(1.1,2,0.1,1.5,4,2.5,-1,7,-3,6,-8,5).$$

		From $(\lambda)_{(0)}=(2,4,-1,7,-3,6,-8,5)$, we can get a $D$-type special partition ${\bf p}_{0}=[4,4,2^4]$, which is a very even partition. We can write $(\lambda)_{(0)}=w_1\mu_0$ with $\mu=(-8,-7,-6,-5,-4,-3,-2,1)$ being antidominant and $w_1=(-6,-8,4,-7,-1,-3,-5,2)$ being the unique minimal length element in $W_{[\mu]}^J$. We have $w_1 \stackrel{LR}{\sim} w_{II}$, where $w_{II}$ is a very even element of type II for the partition ${\bf p}_{0}$, see Lemma \ref{w12}. Thus $\mathcal{O}_{{\bf p}_{0}}$ is a very even nilpotent orbit of type II.
		
  From $(\lambda)_{(\frac{1}{2})}=(1.5,2.5)$, we can get a $C$-type metaplectic special partition ${\bf p}_{\frac{1}{2}}=[2,2]$. 
  We can write $(\lambda)_{(\frac{1}{2})}=w_2\mu_1$ with $\mu_1=(-2.5,-1.5)$ being antidominant and $w_2=(-2,-1)$ being the unique minimal length element in $W_{[\mu_1]}^J$. We have $w_2 \stackrel{LR}{\sim} w_{II}$, where $w_{II}$ is a very even element of type II for the partition ${\bf p}_{\frac{1}{2}}$. Thus $\mathcal{O}_{{\bf p}_{\frac{1}{2}}}$ is a very even nilpotent orbit of type II.

  From $\lambda_{Y_1}=(1.1,0.1)$, we can get a $A$-type partition ${\bf p}_{1}=[1,1]$. Thus we have ${\bf p}_{0}{\stackrel{c}{\sqcup}}{\bf p}_{\frac{1}{2}}{\stackrel {c}{\sqcup}}2{\bf p}_1=[8,8,2^4]$. Its $D$-collapse is the same, i.e., ${\bfpp_{\lambda}}=[8,8,2^4]$. Then we have 
		$V(\Ann (L(\lambda)))=\overline{\mathcal{O}}_{[8,8,2^4]}$.
		Note that ${\bfpp_{\lambda}}=[8,8,2^4]$ is a very even partition. From Corollary \ref{nemural-d}, we know that ${\mathcal{O}}_{\bfpp_{\lambda}}$ is a very even nilpotent orbit of type I and $V(\Ann (L(\lambda)))=\overline{\mathcal{O}}^{I}_{[8,8,2^4]}$.

	\end{ex}

 \section{Hollow diagram algorithm}\label{refi}
	Recall that in the partition algorithm, we used symbols to obtain the partitions ${\bf p}_{0}$ and  ${\bf p}_{\frac{1}{2}}$. In this section we want to give another algorithm to obtain the partitions ${\bf p}_{0}$ and  ${\bf p}_{\frac{1}{2}}$ without  using symbols.
\subsection{Hollow diagram}	
	
	Firstly we recall some notations from \cite{BXX}.
	
	For a Young diagram $P$, we use $ (k,l) $ to denote the box in the $ k $-th row and the $ l $-th column.
	We say that the box $ (k,l) $ is \textit{even} (resp. \textit{odd}) if $ k+l $ is even (resp. odd). Let $ p_i \ev$ (resp. $ p_i\od $) be the number of even (resp. odd) boxes in the $ i $-th row of the Young diagram $ P $.
	Then we have \begin{equation*}
	p_i\ev=\begin{cases}
	\left\lceil \frac{p_i}{2} \right\rceil&\text{ if } i \text{ is odd},\\
	\left\lfloor \frac{p_i}{2} \right\rfloor&\text{ if } i \text{ is even},
	\end{cases}
	\quad p_i\od=\begin{cases}
	\left\lfloor \frac{p_i}{2} \right\rfloor&\text{ if } i \text{ is odd},\\
	\left\lceil \frac{p_i}{2} \right\rceil&\text{ if } i \text{ is even}.
	\end{cases}
	\end{equation*}
	Here for $ s\in \mathbb{R} $, $ \lfloor s \rfloor $ is the largest integer $ n $ such that $ n\leq s $, and $ \lceil s \rceil$ is the smallest integer $n$ such that $ n\geq s $. For convenience, we set
	\begin{equation*}
	p\ev=(p_1\ev,p_2\ev,\dots)\quad\mbox{and}\quad p\od=(p_1\od,p_2\od,\dots).
	\end{equation*}
	
	
	\begin{ex}\label{evenbox}
		Let $P$ be a Young diagram with shape $ p=[5,5,3,3,3] $. The even and odd boxes in $P$ are marked as follows:
		\[
		\tytb{EOEOE,OEOEO,EOE,OEO,EOE}.
		\]
		Then $ p\ev=(3,2,2,1,2)$ and $ p\od=(2,3,1,2,1) $.
	\end{ex}
	
	\begin{Defff}
		By removing all the odd boxes from a Young  diagram  $P$, we obtain a  diagram $ P^{\ev} $  consisting of only even boxes and inheriting the filling from $ P $. We say $ P^{\ev} $ is the \textit{even  diagram} of $ P $. Similarly we can define $ P^{\od} $. These  diagrams are called {\it Hollow  diagrams} of $P$ in \cite{BXX}.
	\end{Defff}
	In the following, we identify a partition and its corresponding Young diagram. If ${\bf p}$ is the partition corresponding to the 	Young diagram $P$, we will identify  $ p\ev$ with $P^{\ev}$, and $ p\od$ with $P^{\od}$.
	\begin{ex}
		When $P $ is the Young diagram in Example \ref{evenbox}, we have
		\[
		P^{\ev}=\tytb{E\none E\none E, \none E\none E, E\none E, \none E\none, E\none E}
		\]

		and
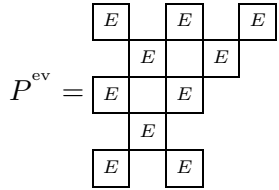
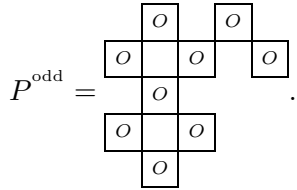
\[
		P^{\od}=\tytb{\none O\none O\none,O\none O\none O,\none O, O\none O,\none O}.
		\]	
	\end{ex}

Recall that in \S \ref{Weyl}, from a partition ${\bf p}$ or nilpotent orbit of type $B$, $C$ or $D$, we can get a symbol $\Lambda$ and a special symbol $\Lambda^s$. From $\Lambda^s$, we can get a special partition ${\bf p}^s$ or special nilpotent orbit.

 \begin{Lem}\label{odev}
 Let ${\bf p}$ be a partition of type $B$, $C$ or $D$.
 We use $P$ and $\bar{P}$ to denote the Young diagrams corresponding to  the partition
	${\bf p}$ and its corresponding special partition
 ${\bf p}^s$. Then we have $P^{\od}=\bar{P}^{\od}$ for types $B$ and $C$, $P^{\ev}=\bar{P}^{\ev}$ for type $D$.
	\end{Lem}
	\begin{proof}
	    From the proof of Proposition $3.3$ in \cite{BXX} or the process of constructing a symbol from a given partition ${\bf p}=[p_1,\dots,p_{2m+1}]$ of type $B$ or $C$, we assume that the corresponding $B$-symbol is
	$\Lambda_B= \begin{pmatrix}
	\lambda_1~\lambda_2~\dots~\lambda_{m+1}\\
	\mu_1~\mu_2~\dots~\mu_{m}
	\end{pmatrix} $. Then
	\[
	\{\lambda_i,\mu_j \mid  1\leq  i\leq m+1, 1\leq j \leq m \}  = \left\{ \alpha_k\mid 1\leq k\leq 2m+1\right \}
	\]
	with $ \alpha_k=p_k\od+ \left\lfloor \frac{2m+1-k  }{2}\right\rfloor $.
	
	Similarly from the proof of Proposition $3.6$ in \cite{BXX} or the process of constructing a symbol from a given partition
	${\bf p}=[p_1,\dots,p_{2m}]$  of type $D$, we assume that the corresponding $D$-symbol is
	$\Lambda_D= \begin{pmatrix}
	\lambda_1~\lambda_2~\dots~\lambda_{m}\\
	\mu_1~\mu_2~\dots~\mu_{m}
	\end{pmatrix} $. Then
	\[
	\{\lambda_i,\mu_j \mid  1\leq  i\leq m, 1\leq j \leq m \}  = \left\{ \beta_k\mid 1\leq k\leq 2m\right \}
	\]
	with $ \beta_k=p_k\ev+ \left\lfloor \frac{2m-k  }{2}\right\rfloor $.

This completed the proof.
	\end{proof}

\subsection{ H-algorithm of classical types}

\def\RSK{\mathrm{RSK}}
\def\bfss{\mathbf{s}}
\def\cF{\mathcal{F}}

From Theorem \ref{zuheB} and \ref{zuheD}, for an integral weight $\lambda\in \mathfrak{h}^*$,  we can write $\lambda=w\mu$ for a unique $w\in W^J$ and a unique anti-dominant $\mu\in \mathfrak{h}^*$.
 Also we have $p(\lambda^-)=p({}^-{w})$ (for types $B$ and $C$) and $p(\lambda^-)=p({}^-{w'})$ ($w'=w$ or $w'=wt$).
 We call it a {\it domino type partition  of $2n$}. In general, a partition $\bf p$ is also called a domino type partition  of $2n$ if it has the same corresponding diagram with some partition $p({}^-{w})$ for some $w\in W_n$.

We consider the type $C$ case first.

For any nilpotent orbit $\cO$  (identified with a Young diagram $P$ and its partition ${\bf p}$) of type $C_n$, let  $P^{\od}$ denote the odd diagram of $P$.
Let $\Lambda_{\cO}=\Lambda_{P}$ be the symbol attached to $\cO$ (or  $P$) by Barbasch-Vogan's algorithm in \cite{BarV82} or the algorithm in \S \ref{Weyl}.

\begin{Lem}
	Let $w\in W_n$ and $P({}^-{w})$ be the diagram obtained by the R-S algorithm.
	Then there exists a unique special Young diagram $\bar{P}$ (or special orbit $\cO^s$) such that
	$\bar{P}^{\od} = P^{\od}$.
	Moreover, $\bar{P}$ is the maximal diagram (under the closure relation of corresponding nilpotent orbits) in
	\[
	\cF_{P} := \{P' \text{ is of type } C \mid P'^{\od} = P^{\od} \}.
	\]
\end{Lem}
\begin{proof}
	By  \cite{BarV82}, from the Young diagram $P$, we have a  symbol $\Lambda_{P}$ attached to $P$, which has rank $n$ and defect $1$. By re-arranging the symbol $\Lambda_P$, we have a unique special symbol $\Lambda_{P}^s$. Let $\cO^s$ be the $C$-type orbit corresponding to $\Lambda_{P}^s$ under the Springer correspondence. Note that $\cO^s$ is special since $\Lambda_{P}^s$ is a special symbol.  We use $\bar{P}$ to denote the Young diagram for $\cO^s$.
	By Lemma \ref{odev}, $\bar{P}^{\od}= P^{\od}$. In particular, we see that $\cF_{P}$ is non-empty.
	Furthermore, by the arguments in Lemma \ref{odev}, for any $C$-type  Young diagram $P'$,
	$P'^{\od}= P^{\od}=\bar{P}^{\od}$  if and only if $\Lambda_{P'}$ and $\Lambda_{P}$ (or $\Lambda_{\bar{P}}$)
	have the same set of entries, if and only if they are in the same double cell by Proposition \ref{dcell}. In other words,
	\[
	\cF_{P} = \{P' \text{ is of type } C \mid \text{$\Lambda_{P'}$ is in the family of $\Lambda_{\bar{P}}$} \}.
	\]
We use $\cO'$ to denote the nilpotent orbit for $P'$.	By \cite[Corollary 7.2]{KT}, $\Lambda_{P'}$ is in the family of  $\Lambda_{\bar{P}}$ implies $\cO' \subset \overline{\cO^s}$.
	Hence $\bar{P}$ is the maximal diagram in $\cF_{P}$.
\end{proof}

Now we give the algorithm to obtain the special partition  ${\bf p}^s$ (partition of $\cO^s$) from a domino type partition ${\bf p}$ (corresponding to a Young diagram $P$) as follows.
\begin{Defff}[H-algorithm of type $C$]
    Let ${\bf p}$ be a domino type partition (whose Young diagram is $P$) of $2n$, then we can get a special partition  ${\bf p}^s$ of type $C_n$ by the following steps:
    \begin{enumerate}
\item Construct  the Hollow diagram $P^{\od}$ consisting of odd boxes;
	\item Label the rows starting from $1$ but avoid all the consecutive rows
	ending with the shape  $\tytb{O,\none O}$ (when  two consecutive rows
	has the shape  $\tytb{E,O}$ in $P$, these two rows will not be labeled);
\item Keep odd labeled rows unchanged and put $\tytb{E}$ on the end of each even labeled row;
\item Fill the holes and we are done.
\end{enumerate}
We call the above algorithm {\it Hollow diagram algorithm} or {\it H-algorithm of type $C$}.
\end{Defff}

\begin{ex}Let ${\bf p}=[5,4,3,3,3,3,2,1]$ be a domino type partition  of $24$. Then we have
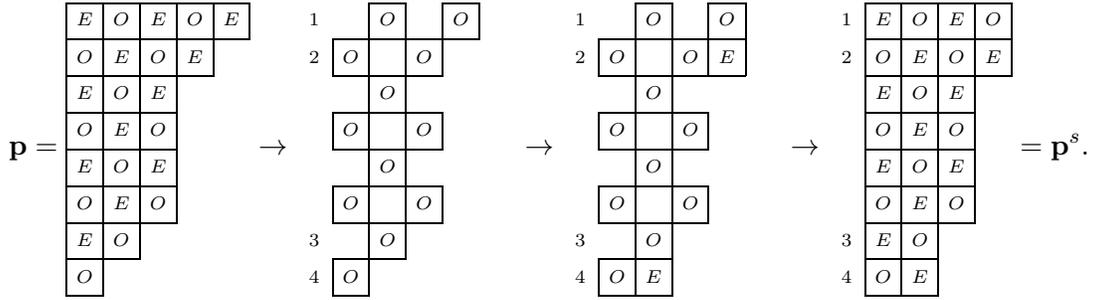
\[{\bf p}=
\tytb{EOEOE,OEOE,EOE,OEO,EOE,OEO,EO,O} \to
\tytb{{\none[1]}\none O\none O\none,{\none[2]}O\none O,\none\none O,\none O\none O,\none\none O,\none O\none O,{\none[3]}\none O,{\none[4]}O}
\to
\tytb{{\none[1]}\none O\none O\none,{\none[2]}O\none OE,\none\none O,\none O\none O,\none\none O,\none O\none O,{\none[3]}\none O,{\none[4]}OE}
\to
\tytb{{\none[1]}EOEO,{\none[2]}OEOE,\none EOE,\none OEO,\none EOE,\none OEO,{\none[3]}EO,{\none[4]}OE}={\bf p}^s.
\]\\
Thus ${\bf p}^s=[4,4,3,3,3,3,2,2]$ is a special partition of type $C_{12}$.
\end{ex}

 \begin{Lem}\label{domino}
 Given a domino type partition ${\bf p}$, the H-algorithm of type $C$ produces the desired special partition ${\bf p}^s$ of type $C$.
 \end{Lem}
 \begin{proof}
 Because of the shape $\tytb{O,\none O}$, we know that the unlabeled rows are always the same for $P$ and $P'\in \cF_{P}$.
Note that deleting (or adding) a pair of the same length of rows does not change the type of the partition (still, type $C$) and is order-preserving.
So we can reduce the problem to the case when every row is labeled by deleting the unlabeled rows (and then putting them back in the end).

We now argue by induction on the number of rows.
In the rest of the proof, we constantly use the fact that $\bar{P}$ (Young diagram of ${\bf p}^s$) is a $C$-type 
diagram such that $\bar{P}^{\od}=P^{\od}$.
If  $\bar{P}$ has only one or two rows, then $P$ has at most $3$ rows. It is clear that the partition or diagram obtained by using the H-algorithm from $P$ is the same with $\bar{P}$, which is maximal in $\cF_{P}$.

Now we assume $\bar{P}$ has at least $3$ rows.
We have the following two cases:
\begin{enumerate}
    \item If the numbers of boxes in the first and second rows of $\bar{P}$ are even, then we can delete the first two rows and reduce the case to the case with a smaller number of rows since the H-algorithm will not change these two rows. If the numbers of boxes in the first and second rows of $\bar{P}$ are odd and equal, then we can delete the first two rows and reduce the case to the case with a smaller number of rows since these two rows are not labeled.
    \item Otherwise, the  number of boxes in the first row of $\bar{P}$ is even and the end box of
the second row is an odd box.
Now the shape of the top four rows in $\bar{P}$ has the following form (the dotted parts may or may not exist)
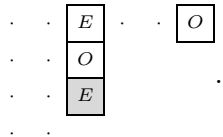
\[
\tytb{{\none[\cdot]}{\none[\cdot]}E{\none[\cdot]}{\none[\cdot]}O,{\none[\cdot]}{\none[\cdot]}O,{\none[\cdot]}{\none[\cdot]}{*(srcol)E},{\none[\cdot]}{\none[\cdot]}}.
\]
Note that below the shaded box $\tytb{{*(srcol)E}}$   there are no other boxes.
So bringing the shaded box $\tytb{{*(srcol)E}}$ to the end of the second row will produce a larger diagram in $\cF_{P}$, with the top four rows of the shape
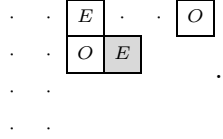
\[
\tytb{{\none[\cdot]}{\none[\cdot]}E{\none[\cdot]}{\none[\cdot]}O,{\none[\cdot]}{\none[\cdot]}O{*(srcol)E},{\none[\cdot]}{\none[\cdot]},{\none[\cdot]}{\none[\cdot]}}.
\]
Therefore this case can not happen since $\bar{P}$ is maximal in $\cF_{P}$.

\end{enumerate}

This finished the proof.
 \end{proof}
 The H-algorithm of type $B$ is similar to type $C$.
 \begin{Defff}[H-algorithm of type $B$]
    Let ${\bf p}$ be a domino type partition (whose Young diagram is $P$) of $2n$, then we can get a special partition  ${\bf p}^s$ of type $B_n$ by the following steps:
    \begin{enumerate}
\item Construct  the Hollow diagram $P^{\od}$ consisting of odd boxes;
	\item Label the rows starting from $1$ but avoid all the consecutive rows
	ending with the shape  $\tytb{O,\none O}$;
\item Keep even labeled rows unchanged and put $\tytb{E}$ on the end of each odd labeled row;
\item Fill the holes. Then if there are only $2n$ boxes in our new Young diagram, we put a box $\tytb{E}$  below the last row and we are done. If there are  $2n+1$ boxes in our new Young diagram, we are done.
\end{enumerate}
We call the above algorithm  {\it H-algorithm of type $B$}.
\end{Defff}

\begin{ex}Let ${\bf p}=[6,4^3,2,2,1,1]$ be a domino type partition  of $24$. Then we have
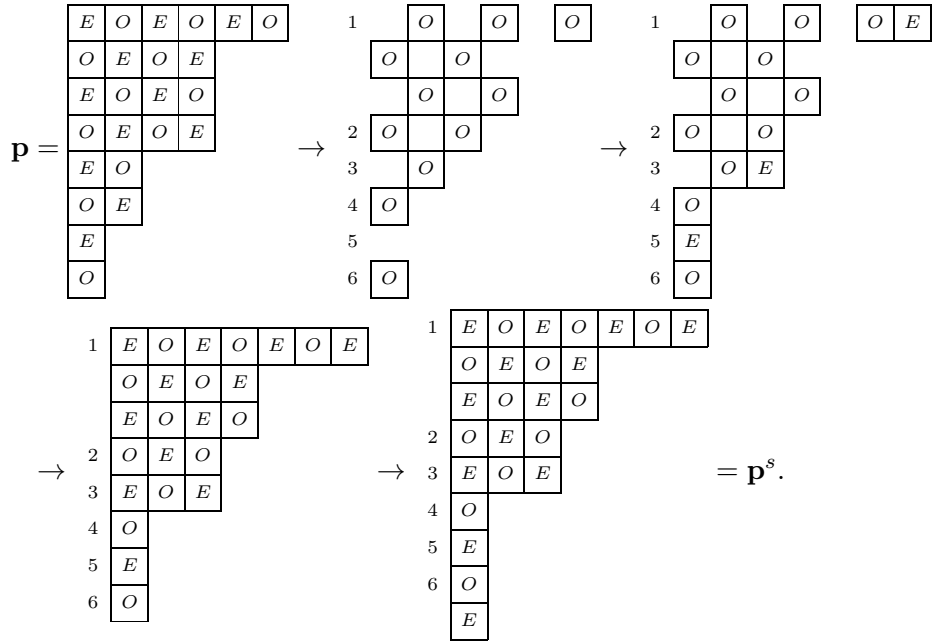
\begin{align*}
    {\bf p}&=\tytb{EOEOEO,OEOE,EOEO,OEOE,EO,OE,E,O} \to
\tytb{{\none[1]}\none O\none O\none O, {\none}O\none O\none,\none\none O\none O, {\none[2]} O\none O,{\none[3]}\none O,{\none[4]} O\none, {\none[5]},{\none[6]}O}
\to
\tytb{{\none[1]}\none O\none O\none O E, {\none}O\none O\none,\none\none O\none O,{\none[2]} O\none O,{\none[3]}\none O E,{\none[4]} O\none,{\none[5]}E, {\none[6]}O}\\
&\to
\tytb{{\none[1]}E OE OE O E ,{\none}OE OE,\none E O E O,{\none[2]} OE O,{\none[3]}E O E,{\none[4]} O\none,{\none[5]}E,{\none[6]}O}
\to \tytb{{\none[1]}E OE OE O E,{\none}OE OE,\none E O E O,{\none[2]} OE O,{\none[3]}E O E,{\none[4]} O\none,{\none[5]}E,{\none[6]}O, \none E}
={\bf p}^s.
\end{align*}

Thus ${\bf p}^s=[7,4,4,3,3,1,1,1,1]$ is a special partition of type $B_{12}$.
\end{ex}

\begin{Defff}[H-algorithm of type $D$]
    Let ${\bf p}$ be a domino type partition (whose Young diagram is $P$) of $2n$, then we can get a special partition  ${\bf p}^s$ of type $D_n$ by the following steps:
    \begin{enumerate}
\item Construct  the Hollow diagram $P^{\ev}$ consisting of even boxes;
	\item Label the rows starting from $1$ but avoid all the consecutive rows
	ending with the shape  $\tytb{E,\none E}$;
\item Keep odd labeled rows unchanged and put $\tytb{O}$ on the end of each even labeled row;
\item Fill the holes. Then if there are only $2n-1$ boxes in our new Young diagram, we put a box $\tytb{O}$  below the last row and we are done. If there are  $2n$ boxes in our new Young diagram, we are done.
\end{enumerate}
We call the above algorithm  {\it H-algorithm of type $D$}.
\end{Defff}

Suppose ${\bf p}=[p_1,\dots, p_N]$ is  a (special) partition of type $X$ and $p_1\geq p_2\geq \dots \geq p_N>0$. Then $N$ is odd (resp. even) for type $B$ (resp. type $D$). So the fourth step in types $B$ and $D$ are different with type $C$.

\begin{ex}Let ${\bf p}=[6, 4, 2^5, 1^4]$ be a domino type partition  of $24$. Then we have
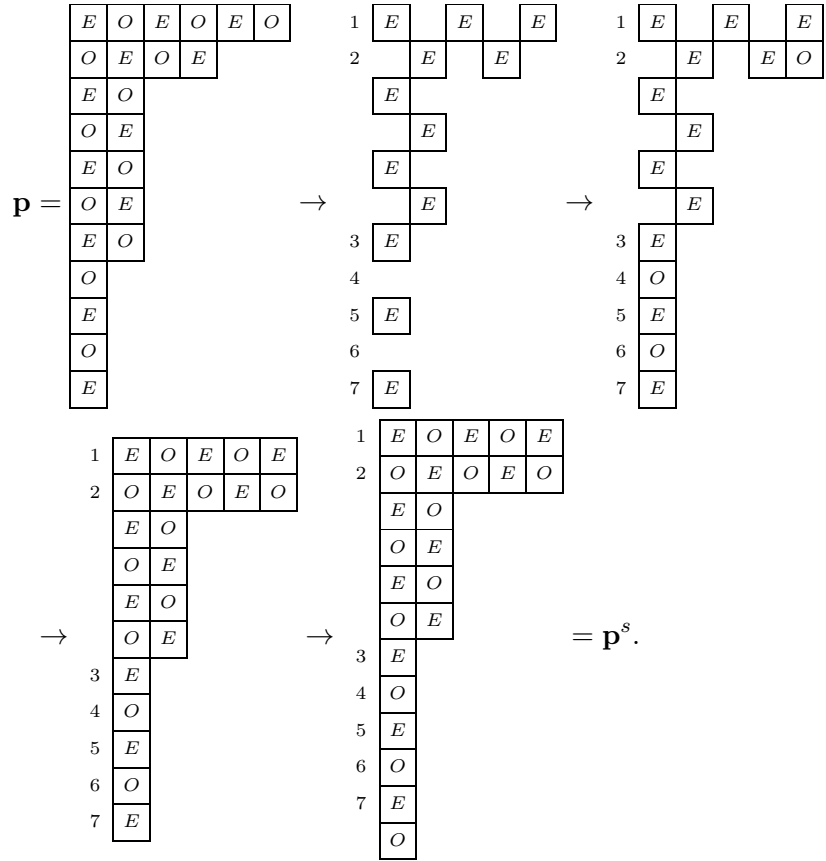
\begin{align*}
    {\bf p}&=\tytb{EOEOEO,OEOE,EO,OE,EO,OE, EO, O, E, O, E} \to
\tytb{{\none[1]}E\none E\none E, {\none[2]}\none E\none E,\none E\none, {\none}\none E,{\none} E\none,{\none} \none E, {\none[3]} E\none,{\none[4]},{\none[5]}E, {\none[6]},{\none[7]}E}
\to
\tytb{{\none[1]}E\none E\none E, {\none[2]}\none E\none EO,\none E\none, {\none}\none E,{\none} E\none,{\none} \none E, {\none[3]} E\none,{\none[4]}O,{\none[5]}E, {\none[6]}O,{\none[7]}E}\\
&\to
\tytb{{\none[1]}EO EO E, {\none[2]}O EO EO,\none EO, {\none}O E,{\none} EO,{\none} O E, {\none[3]} E\none,{\none[4]}O,{\none[5]}E, {\none[6]}O,{\none[7]}E}
\to \tytb{{\none[1]}EO EO E, {\none[2]}O EO EO,\none EO, {\none}O E,{\none} EO,{\none} O E, {\none[3]} E\none,{\none[4]}O,{\none[5]}E, {\none[6]}O,{\none[7]}E, \none O}
={\bf p}^s.
\end{align*}

Thus ${\bf p}^s=[5,5,2^4,1^6]$ is a special partition of type $D_{12}$.
\end{ex}

Similar to Lemma \ref{domino}, for types $B$ and $D$ we have the same result. Thus we have the follows.
\begin{Cor}\label{special-partition}
    For any integral weight $\lambda$ of type $X$ ($X=B,C,D$), we can always get a special partition ${\bf p}^s$ after the $H$-algorithm of type $X$ from the domino type partition ${\bf p}=p(\lambda^-)$.
\end{Cor}
\begin{proof}
		For an integral weight $\lambda\in \mathfrak{h}^*$,  we can write $\lambda=w\mu$ for a unique $w\in W^J$ and a unique anti-dominant $\mu\in \mathfrak{h}^*$.
		Then from Theorem \ref{zuheB} and Theorem \ref{zuheD}, we have   $p(\lambda^-)=p({}^-{w})$ for types $B$ and $C$, and $p(\lambda^-) =p({}^-{w})$ or $p({}^-({wt}))$ for type $D$.   Thus after the H-algorithm, we can get a special partition corresponding to the special symbol obtained from  $p(\lambda^-)=p({}^-{w})$ (for types $B$ and $C$) or $p(\lambda^-)=p({}^-{w'})$ ($w'=w$  or $w'=wt$ for type $D$).
\end{proof}

From now on, we use $H_X({\bf p})$ to denote the special partition of type $X$ after the H-algorithm of type $X$.

Recall that for a given $w\in W_n^X$ (Weyl group of type $X$), we can get a
	Young tableau $ P({}^-{ w}) $ by R-S algorithm. Then we can get a symbol $\Lambda_w$. By some permutation, we can get a special symbol $\Lambda_w^s$ of type $B$ or $D$, then we can get the corresponding special nilpotent orbit $\mathcal{O}_w=\mathcal{O}_{\Lambda_w^s}$ of type $X$. Now from the H-algorithm, we have the following result.
\begin{Cor}
    For any $x,y\in W^X$, we have: $$x \stackrel{LR}{\sim} y \text{~ if~ and~ only~ if ~} H_X(p({}^-{ x}))=H_X(p({}^-{y})),$$
   except that for type $D$, when $n=2m'$ is even and $H_D(p({}^-{ x}))=H_D(p({}^-{y}))$ is a very even partition, we have $x \stackrel{LR}{\sim} y $ 
    if and only if there exists some very even element $ w=w_I$ or $w_{II}$ such that $x\stackrel{LR}{\sim} w$ and $y\stackrel{LR}{\sim}w$.

\end{Cor}
\begin{proof}
From Corollary \ref{special-partition}, we note that  $\mathcal{O}_w$ has the same partition with $H_X(p({}^-{ w}))$  for any $w\in W_n^X$. Then the result follows from Theorem \ref{zuheD}.
\end{proof}

\subsection{H-algorithm of metaplectic type}
The partition $p((\lambda)_{(\frac{1}{2})}^-)$ will give us a special partition  ${\bf b}$ of type $D$ by the $H$-algorithm of type $D$, then  we have  ${\bf p}_{\frac{1}{2}}={\bf b}$ for type $B$, ${\bf p}_{\frac{1}{2}}=(({\bf b}^*)_D)^*$ for types $C$ and $D$. In other words, we have $${\bf p}_{\frac{1}{2}}=\left(\left(\left({\bf b}\right)^*\right)_D\right)^*$$ for types $C$ and $D$.

In the following, we want to give some simpler algorithm to compute this  partition ${\bf p}_{\frac{1}{2}}$.

\begin{Defff}[H-algorithm of  metaplectic type]
    Let ${\bf p}$ be a domino type partition (whose Young diagram is $P$) of $2n$, then we can get a metaplectic special partition  ${\bf p}_{\frac{1}{2}}$ of type $C_n$ by the following steps:
    \begin{enumerate}
\item Construct  the hollow diagram $P^{\ev}$ consisting of even boxes;
	\item Label the rows starting from $1$ but avoid all the consecutive rows
	ending with the shape  $\tytb{E,\none E}$ (when  two consecutive rows
	has the shape  $\tytb{O, E}$  in $P$, these two rows will not be labeled);
\item Keep even labeled rows unchanged and put $\tytb{O}$ on the end of each odd labeled row;
\item Fill the holes. Then if there are only $2n-1$ boxes in our new Young diagram, we put a box $\tytb{O}$  below the last row and we are done. If there are  $2n$ boxes in our new Young diagram, we are done.
\end{enumerate}
We call the above algorithm  {\it H-algorithm of metaplectic type}.
\end{Defff}

Recall the definition of the metaplectic cell, see \cite[Section~6.2]{BMSZ-typeC}.
Two representations $\sigma_1,\sigma_2\in \Irr(W_n)$ are in the same metaplectic cell if and only if
\begin{itemize}
   \item
$\sigma_1|_{W'_n}$ and $\sigma_2|_{W'_n}$  are reducible and $\sigma_1=\sigma_2$, or
\item
$\sigma_1|_{W'_n}$ and $\sigma_2|_{W'_n}$ are irreducible and in the same double cell of $\Irr(W'_n)$.
\end{itemize}
In other words, $\sigma_1$ and $\sigma_2$ are in the same metaplectic double cell if and only if
the $D$-symbols of all irreducible components of $\sigma_1|_{W'_n}$
and $\sigma_2|_{W'_n}$ are the same.




Let $\sigma_P$ be the irreducible $W_n$-representation attached to a given Young diagram $P$ (of domino type) via Barbasch-Vogan's algorithm in \cite{BarV82} or the algorithm in \S \ref{Weyl}.

\begin{Lem}
Let $w\in W'_n$ and $P({}^-{w})$ be the diagram obtained by the R-S algorithm.
Then there is a unique metaplectic special Young diagram $\barP$ such that $\barP\ev=P\ev$.
Moreover, $\barP$ is the unique  maximal dimensional element in
\[
\cF_P:=\Set{\text{$Q$ is of type $C$} | Q\ev = P\ev}.
\]
\end{Lem}
\begin{proof}
Let $\sigma_s$ be the unique metaplectic spacial representation in the metaplectic double cell
of $\sigma_P$, which is the unique representation in the metaplectic double cell with minimal fake degree (see \cite[Lemma~6.2]{BMSZ-typeC}).
Let $\barP$ be the $C$-type orbit (Young diagram) attached to $\sigma_s$ via the Springer correspondence.
Let $Q$ be a $C$-type partition (Young diagram)  such that $Q\ev=P\ev$ and $Q\ne P$.
By the arguments in Lemma \ref{odev}, the condition $Q\ev=P\ev$ is equivalent to that $\sigma_{Q}$ is in the metaplectic double cell of $\sigma_P$, see \cite[p419~Type~$C$]{Ca85}.
Therefore by Proposition \ref{intehral-W}, the fake degree of $\sigma_{Q}$ is strictly greater than the fake degree of $\sigma_s$.
Using the fact that  $2$ times the fake degree of $\sigma_Q$ equals the codimension of $Q$ in the nilpotent cone, we proved the lemma.
\end{proof}

Similar to Lemma \ref{domino}, we have the following result.

 \begin{Lem}\label{domino-meta}
 Given a domino type partition ${\bf p}$, the H-algorithm of metaplectic type produces the desired partition ${\bf p}_{\frac{1}{2}}$.
 \end{Lem}
\begin{Cor}
    For any weight $\lambda$ of type $X$ ($X=C,D$), we can always get a metaplectic special partition ${\bf p}_{\frac{1}{2}}$ after the H-algorithm of metaplectic type from the domino type partition $p((\lambda)_{(\frac{1}{2})}^-)$.
\end{Cor}

\begin{ex}Let ${\bf p}=[6, 4, 3,3]$ be a domino type partition  of $16$. Then we have
\begin{align*}
    {\bf p}&=\tytb{EOEOEO,OEOE,EOE,OEO} \to
    \tytb{{\none[1]}E\none E\none E, {\none[2]}\none E\none E, {\none[3]} E\none E, {\none[4]}\none E} \to
\tytb{{\none[1]}E\none E\none EO, {\none[2]}\none E\none E, {\none[3]} E\none EO, {\none[4]}\none E}
\to
\tytb{{\none[1]}EO EO EO, {\none[2]}O EOE, {\none[3]}EO EO,{\none[4]} OE}
={\bf p}^s.
\end{align*}

Thus ${\bf p}^s=[6,4,4,2]$ is a metaplectic special partition of type $C_{8}$.
\end{ex}

\subsection*{Acknowledgments}
	 Z. Bai was supported  by the National Natural Science Foundation of
	China (No. 12171344) and the National Key $\textrm{R}\,\&\,\textrm{D}$ Program of China (No. 2018YFA0701700 and No. 2018YFA0701701). J.-J. Ma was supported by the National Natural Science Foundation of China (No. 11701364 and No. 11971305) and Xiamen University Malaysia Research Fund (No. XMUMRF/2022-C9/IMAT/0019).	Y. Wang was supported by Hui-Chun Chin and Tsung-Dao Lee Chinese Undergraduate Research Endowment (CURE). 	We would like to thank  D. A. Vogan for helpful discussions about $\tau$-invariants. We also would like to thank D. Garfinkle Johnson for helpful discussions about domino tableaux.

\printbibliography

@ARTICLE{BXX,
  author = {Bai, Z. and Xiao, W. and Xie, X.},
  title = {Gelfand-Kirillov dimensions and associated varieties of highest weight
	modules},
  journal = {Int. Math. Res. Not. IMRN},
  year = {2023},
pages = {8101-8142}
}

@ARTICLE{BaiX,
  author = {Bai, Z. and Xie, X.},
  title = {Gelfand-{K}irillov dimensions of highest weight {H}arish-{C}handra
	modules for {$SU(p,q)$}},
  journal = {Int. Math. Res. Not. IMRN},
  year = {2019},
  pages = {4392--4418}
}

@ARTICLE{BMSZ-ABCD,
  author = {Barbasch, D. and Ma, J.-J. and Sun, B. and Zhu, C.-B.},
  title = {Special unipotent representations of real classical groups: construction and unitarity},
  journal = {arXiv:1712.05552v5,},
year = {2023}
}

@ARTICLE{BMSZ-counting,
  author = {Barbasch, D. and Ma, J.-J. and Sun, B. and Zhu, C.-B.},
  title = {Special unipotent representations of real classical groups: counting
	and reduction to good parity},
  journal = { arXiv:2205.05266v3,},
year = {2023}
}

@ARTICLE{BMSZ-typeC,
  author = {Barbasch, D. and Ma, J.-J. and Sun, B. and Zhu, C.-B.},
  title = {On the notion of metaplectic Barbasch-Vogan duality},
  journal = {arXiv:2010.16089v2,},
year = {2022}
}

@ARTICLE{BarV82,
  author = {Barbasch, D. and Vogan, D. A.},
  title = {Primitive ideals and orbital integrals in complex classical groups},
  journal = {Math. Ann.},
  year = {1982},
  volume = {259},
  pages = {153--199},
  number = {2},
 
  fjournal = {Mathematische Annalen}
}

@ARTICLE{BV83,
  author = {Barbasch, D. and Vogan, D.},
  title = {Primitive ideals and orbital integrals in complex exceptional groups},
  journal = {J. Algebra},
  year = {1983},
  volume = {80},
  pages = {350--382},
  number = {2},
  fjournal = {Journal of Algebra}
}

@BOOK{BB05,
  title = {Combinatorics of {C}oxeter groups},
  publisher = {Springer, New York},
  year = {2005},
  author = {Bj\"{o}rner, A. and Brenti, F.},
  volume = {231},
  pages = {xiv+363},
  series = {Graduate Texts in Mathematics},
  isbn = {978-3540-442387; 3-540-44238-3},
  mrclass = {05-01 (05E15 20F55)},
  mrnumber = {2133266},
  mrreviewer = {Jian-yi Shi}
}

@ARTICLE{BoB1,
  author = {Borho, W. and Brylinski, J.-L.},
  title = {Differential operators on homogeneous spaces. {I}. {I}rreducibility
	of the associated variety for annihilators of induced modules},
  journal = {Invent. Math.},
  year = {1982},
  volume = {69},
  pages = {437--476},
  number = {3},
  
  fjournal = {Inventiones Mathematicae}
}

@BOOK{BBM,
  title = {Nilpotent orbits, primitive ideals, and characteristic classes},
  publisher = {Birkh\"{a}user Boston, Inc., Boston, MA},
  year = {1989},
  author = {Borho, W. and Brylinski, J.-L. and MacPherson, R.},
  volume = {78},
  pages = {vi+131},
  series = {Progress in Mathematics},
  note = {A geometric perspective in ring theory}
 
}

@ARTICLE{BJ77,
  author = {Borho, W. and Jantzen, J. C.},
  title = {\"{U}ber primitive {I}deale in der {E}inh\"{u}llenden einer halbeinfachen
	{L}ie-{A}lgebra},
  journal = {Invent. Math.},
  year = {1977},
  volume = {39},
  pages = {1--53},
  number = {1},
  
  fjournal = {Inventiones Mathematicae}
}

@article {BM,
    AUTHOR = {Borho, W. and MacPherson, R.},
     TITLE = {Repr\'{e}sentations des groupes de {W}eyl et homologie
              d'intersection pour les vari\'{e}t\'{e}s nilpotentes},
   JOURNAL = {C. R. Acad. Sci. Paris S\'{e}r. I Math.},
  FJOURNAL = {Comptes Rendus des S\'{e}ances de l'Acad\'{e}mie des Sciences. S\'{e}rie
              I. Math\'{e}matique},
    VOLUME = {292},
      YEAR = {1981},
    NUMBER = {15},
     PAGES = {707--710},
      ISSN = {0249-6291},
   MRCLASS = {14B05 (17B05 17B35 20G05)},
  MRNUMBER = {618892},
}

@BOOK{Ca85,
  title = {Finite groups of {L}ie type},
  publisher = {John Wiley \& Sons, Inc., New York},
  year = {1985},
  author = {Carter, R. W.},
  pages = {xii+544},
  series = {Pure and Applied Mathematics (New York)},
  note = {Conjugacy classes and complex characters, A Wiley-Interscience Publication},
  isbn = {0-471-90554-2},
  mrclass = {20G40 (20-02 20C15)},
  mrnumber = {794307},
  mrreviewer = {David B. Surowski}
}

@BOOK{CM,
  title = {Nilpotent orbits in semisimple {L}ie algebras},
  publisher = {Van Nostrand Reinhold Co., New York},
  year = {1993},
  author = {Collingwood, D. H. and McGovern, W. M.},
  pages = {xiv+186},
  series = {Van Nostrand Reinhold Mathematics Series},
  isbn = {0-534-18834-6},
  mrclass = {17-02 (17B20 17B25 22E60)},
  mrnumber = {1251060},
  mrreviewer = {Stephen Slebarski}
}

@ARTICLE{Du,
  author = {Duflo, M.},
  title = {Sur la classification des id\'{e}aux primitifs dans l'alg\`ebre enveloppante
	d'une alg\`ebre de {L}ie semi-simple},
  journal = {Ann. of Math. (2)},
  year = {1977},
  volume = {105},
  pages = {107--120},
  number = {1},
 
  fjournal = {Annals of Mathematics. Second Series}
}

@ARTICLE{GT,
  author = {Gao, F. and Tsai, W.-Y.},
  title = {On the wavefront sets associated with theta representations},
  journal = {Math. Z.},
  year = {2022},
  volume = {301},
  pages = {1--40},
  number = {1},
  
  fjournal = {Mathematische Zeitschrift}
}

@ARTICLE{Garfinkle3,
  author = {Garfinkle, D.},
  title = {On the classification of primitive ideals for complex classical {L}ie
	algebras. {III}},
  journal = {Compositio Math.},
  year = {1993},
  volume = {88},
  pages = {187--234},
  number = {2},
  fjournal = {Compositio Mathematica}
}

@ARTICLE{Garfinkle2,
  author = {Garfinkle, D.},
  title = {On the classification of primitive ideals for complex classical {L}ie
	algebras. {II}},
  journal = {Compositio Math.},
  year = {1992},
  volume = {81},
  pages = {307--336},
  number = {3},
  fjournal = {Compositio Mathematica}
}

@ARTICLE{Garfinkle1,
  author = {Garfinkle, D.},
  title = {On the classification of primitive ideals for complex classical {L}ie
	algebras. {I}},
  journal = {Compositio Math.},
  year = {1990},
  volume = {75},
  pages = {135--169},
  number = {2},
  fjournal = {Compositio Mathematica}
}

@ARTICLE{ge61,
  author = {Gerstenhaber, M.},
  title = {Dominance over the classical groups},
  journal = {Ann. of Math. (2)},
  year = {1961},
  volume = {74},
  pages = {532--569},
 
  fjournal = {Annals of Mathematics. Second Series}
}

@ARTICLE{GS,
  author = {Gourevitch, D. and Sahi, S.},
  title = {Annihilator varieties, adduced representations, {W}hittaker functionals,
	and rank for unitary representations of {${\rm GL}(n)$}},
  journal = {Selecta Math. (N.S.)},
  year = {2013},
  volume = {19},
  pages = {141--172},
  number = {1},
 
  fjournal = {Selecta Mathematica. New Series}
}

@ARTICLE{GSK,
  author = {Gourevitch, D. and Sayag, E. and Karshon, I.},
  title = {Annihilator varieties of distinguished modules of reductive {L}ie
	algebras},
  journal = {Forum Math. Sigma},
  year = {2021},
  volume = {9},
  pages = {Paper No. e52, 30},
  
  fjournal = {Forum of Mathematics. Sigma}
}

@BOOK{Hum08,
  title = {Representations of semisimple {L}ie algebras in the {BGG} category
	{$\mathscr{O}$}},
  publisher = {American Mathematical Society, Providence, RI},
  year = {2008},
  author = {Humphreys, J. E.},
  volume = {94},
  pages = {xvi+289},
  series = {Graduate Studies in Mathematics}
 
}

@BOOK{Ja79,
  title = {Moduln mit einem h\"ochsten {G}ewicht},
  publisher = {Springer, Berlin},
  year = {1979},
  author = {Jantzen, J. C.},
  volume = {750},
  pages = {ii+195},
  series = {Lecture Notes in Mathematics}
}

@ARTICLE{Ga-D,
  author = {Johnson, D. G.},
  title = {Edge Transport from Parabolic Subgroups of Type $D_4$},
 journal = { arXiv:1907.09717v1,},
year = {2019}
}

@ARTICLE{Jo85,
  author = {Joseph, A.},
  title = {On the associated variety of a primitive ideal},
  journal = {J. Algebra},
  year = {1985},
  volume = {93},
  pages = {509--523},
  number = {2},
  
  fjournal = {Journal of Algebra}
}

@INPROCEEDINGS{Jo83ICM,
  author = {Joseph, A.},
  title = {Primitive ideals in enveloping algebras},
  booktitle = {Proceedings of the {I}nternational {C}ongress of {M}athematicians,
	{V}ol. 1, 2 ({W}arsaw, 1983)},
  year = {1984},
  pages = {403--414},
  publisher = {PWN, Warsaw}
}

@ARTICLE{Jo84,
  author = {Joseph, A.},
  title = {On the variety of a highest weight module},
  journal = {J. Algebra},
  year = {1984},
  volume = {88},
  pages = {238--278},
  number = {1},
  
  fjournal = {Journal of Algebra}
}

@INCOLLECTION{Jo82,
  author = {Joseph, A.},
  title = {On the classification of primitive ideals in the enveloping algebra
	of a semisimple {L}ie algebra},
  booktitle = {Lie group representations, {I} ({C}ollege {P}ark, {M}d., 1982/1983)},
  publisher = {Springer, Berlin},
  year = {1983},
  volume = {1024},
  series = {Lecture Notes in Math.},
  pages = {30--76}
}

@ARTICLE{J81-1,
  author = {Joseph, A.},
  title = {Goldie rank in the enveloping algebra of a semisimple {L}ie algebra.
	{III}},
  journal = {J. Algebra},
  year = {1981},
  volume = {73},
  pages = {295--326},
  number = {2},
 
  fjournal = {Journal of Algebra}
}

@ARTICLE{J80-1,
  author = {Joseph, A.},
  title = {Goldie rank in the enveloping algebra of a semisimple {L}ie algebra.
	{I}},
  journal = {J.  Algebra},
  year = {1980},
  volume = {65},
  pages = {269--283},
  
  fjournal = {Journal of Algebra}
}

@ARTICLE{J80-2,
  author = {Joseph, A.},
  title = {Goldie rank in the enveloping algebra of a semisimple {L}ie algebra.
	{II}},
  journal = {J. Algebra},
  year = {1980},
  volume = {65},
  pages = {284--306},
 
  fjournal = {Journal of Algebra}
}

@ARTICLE{KT,
  author = {Kashiwara, M. and Tanisaki, T.},
  title = {The characteristic cycles of holonomic systems on a flag manifold
	related to the {W}eyl group algebra},
  journal = {Invent. Math.},
  year = {1984},
  volume = {77},
  pages = {185--198},
  number = {1},
  
  fjournal = {Inventiones Mathematicae}
}

@ARTICLE{KL,
  author = {Kazhdan, D. and Lusztig, G.},
  title = {Representations of {C}oxeter groups and {H}ecke algebras},
  journal = {Invent. Math.},
  year = {1979},
  volume = {53},
  pages = {165--184},
  number = {2},
  
  fjournal = {Inventiones Mathematicae}
}

@article {Lo15,
    AUTHOR = {Losev, I.},
     TITLE = {Dimensions of irreducible modules over {W}-algebras and
              {G}oldie ranks},
   JOURNAL = {Invent. Math.},
  FJOURNAL = {Inventiones Mathematicae},
    VOLUME = {200},
      YEAR = {2015},
    NUMBER = {3},
     PAGES = {849--923},
      
}

@ARTICLE{LY,
  author = {Losev, I. and Yu, S.},
  title = {On Harish-Chandra modules over quantizations of nilpotent orbits},
  journal = {arXiv:2309.11191,},
year = {2023}
}

@BOOK{Lus84,
  title = {Characters of reductive groups over a finite field},
  publisher = {Princeton University Press, Princeton, NJ},
  year = {1984},
  author = {Lusztig, G.},
  volume = {107},
  pages = {xxi+384},
  series = {Annals of Mathematics Studies}}

@ARTICLE{Lu79,
  author = {Lusztig, G.},
  title = {A class of irreducible representations of a {W}eyl group},
  journal = {Nederl. Akad. Wetensch. Indag. Math.},
  year = {1979},
  volume = {41},
  pages = {323--335},
  number = {3},
  fjournal = {Koninklijke Nederlandse Akademie van Wetenschappen. Indagationes Mathematicae}
}

@ARTICLE{lusztig1977symbol,
  author = {Lusztig, G.},
  title = {Irreducible representations of finite classical groups},
  journal = {Invent. Math.},
  year = {1977},
  volume = {43},
  pages = {125--175},
  number = {2},
 
  fjournal = {Inventiones Mathematicae}
}

@article {LS,
    AUTHOR = {Lusztig, G. and Spaltenstein, N.},
     TITLE = {Induced unipotent classes},
   JOURNAL = {J. London Math. Soc. (2)},
  FJOURNAL = {Journal of the London Mathematical Society. Second Series},
    VOLUME = {19},
      YEAR = {1979},
    NUMBER = {1},
     PAGES = {41--52}
      
}

@article {Mc72,
    AUTHOR = {Macdonald, I. G.},
     TITLE = {Some irreducible representations of {W}eyl groups},
   JOURNAL = {Bull. London Math. Soc.},
  FJOURNAL = {The Bulletin of the London Mathematical Society},
    VOLUME = {4},
      YEAR = {1972},
     PAGES = {148--150},
     
}

@ARTICLE{Mc96,
  author = {McGovern, W. M.},
  title = {Left cells and domino tableaux in classical {W}eyl groups},
  journal = {Compositio Math.},
  year = {1996},
  volume = {101},
  pages = {77--98},
  number = {1},
  fjournal = {Compositio Mathematica}
}

@INCOLLECTION{NOT,
  author = {Nishiyama, K. and Ochiai, H. and Taniguchi, K.},
  title = {Bernstein degree and associated cycles of {H}arish-{C}handra modules---{H}ermitian
	symmetric case},
  year = {2001},
  number = {273},
  pages = {13--80},
  note = {Nilpotent orbits, associated cycles and Whittaker models for highest
	weight representations},
  fjournal = {Ast\'{e}risque},
  issn = {0303-1179},
  journal = {Ast\'{e}risque},
  mrclass = {22E46 (14L30 32M15)},
  mrnumber = {1845714}
}

@BOOK{Sagan,
  title = {The symmetric group},
  publisher = {Springer-Verlag, New York},
  year = {2001},
  author = {Sagan, B. E.},
  volume = {203},
  pages = {xvi+238},
  series = {Graduate Texts in Mathematics},
  edition = {Second},
  note = {Representations, combinatorial algorithms, and symmetric functions}
}

@article {Sh79,
    AUTHOR = {Shoji, T.},
     TITLE = {On the {S}pringer representations of the {W}eyl groups of
              classical algebraic groups},
   JOURNAL = {Comm. Algebra},
  FJOURNAL = {Communications in Algebra},
    VOLUME = {7},
      YEAR = {1979},
    NUMBER = {16},
     PAGES = {1713--1745},
 
}

@article {Sommers,
    AUTHOR = {Sommers, E.},
     TITLE = {Lusztig's canonical quotient and generalized duality},
   JOURNAL = {J. Algebra},
  FJOURNAL = {Journal of Algebra},
    VOLUME = {243},
      YEAR = {2001},
    NUMBER = {2},
     PAGES = {790--812}
}

@INCOLLECTION{Ta,
  author = {Tanisaki, T.},
  title = {Characteristic varieties of highest weight modules and primitive
	quotients},
  booktitle = {Representations of {L}ie groups, {K}yoto, {H}iroshima, 1986},
  publisher = {Academic Press, Boston, MA},
  year = {1988},
  volume = {14},
  series = {Adv. Stud. Pure Math.},
  pages = {1--30}
}

@ARTICLE{V79,
  author = {Vogan, D. A.},
  title = {A generalized {$\tau $}-invariant for the primitive spectrum of a
	semisimple {L}ie algebra},
  journal = {Math. Ann.},
  year = {1979},
  volume = {242},
  pages = {209--224},
  number = {3},
  
  fjournal = {Mathematische Annalen}
}
 
\end{document}